\documentclass[11pt]{article}
\usepackage{amsmath} \usepackage{amssymb} \usepackage{latexsym} \usepackage{enumitem} \usepackage{mathrsfs} \usepackage{comment}
\usepackage{color}

\usepackage{cite}

\usepackage[colorlinks=true,urlcolor=blue,
citecolor=red,linkcolor=blue,linktocpage,pdfpagelabels,bookmarksnumbered,bookmarksopen]{hyperref}

\newtheorem{assumption}{Assumption}

\def\qed{ \ \vrule width.2cm height.2cm depth0cm\smallskip}
\newenvironment{proof}{\noindent {\bf Proof.\/}}{$\qed$\vskip 0.1in}

\newcommand{\ba}{\begin{array}}
\newcommand{\ea}{\end{array}}
\newcommand{\be}{\begin{equation}}
\newcommand{\ee}{\end{equation}}
\newcommand{\bea}{\begin{eqnarray}}
\newcommand{\eea}{\end{eqnarray}}
\newcommand{\beaa}{\begin{eqnarray*}}
\newcommand{\eeaa}{\end{eqnarray*}}

\def\dbE{\mathbb{E}}
\def\dbF{\mathbb{F}}

\def\dbL{\mathbb{L}}

\def\dbP{\mathbb{P}}
\def\dbR{\mathbb{R}}

\def\dbQ{\mathbb{Q}}
\def\dbZ{\mathbb{Z}}

\def\rme{\mathrm{e}}
\def\rmd{\mathrm{d}}
\def\mcu{\mathcal{U} }
\def\mcv{\mathcal{V} }

\def\si{\sigma}

\def\Th{\Theta}

\def\O{\Omega}
\def\cF{{\cal F}}
\def\cG{{\cal G}}

\def\cK{{\cal K}}
\def\cL{{\cal L}}

\def\cP{{\cal P}}

\def\cW{{\cal W}}

\def\no{\noindent}

\def\q{\quad}

\def\pa{\partial}

\def\qed{ \hfill \vrule width.25cm height.25cm depth0cm\smallskip}

\newcommand{\basa}{\begin{assumption}}
\newcommand{\easa}{\end{assumption}}

\newcommand{\bas}{\begin{assum}}
\newcommand{\eas}{\end{assum}}

\def\pa{\partial}

\def\1{{\bf 1}}

\def\:{\!:\!}
\def\reff{\eqref}
\def \proof{{\noindent \bf Proof.\quad}}

\def\gxX{\nabla X^{x,\xi}}
\def\gxY{\nabla Y^{x,\xi}}
\def\gxZ{\nabla Z^{x,\xi}}

\def\bU{\bar{U}}
\def\bX{\bar{X}^{\xi}} 
\def\bY{\bar{Y}^{\xi}} 
\def\bZ{\bar{Z}^{\xi}}

at 9pt

\definecolor{alp}{rgb}{0.0, 0.5, 0.0}

\newtheorem{thm}{Theorem}[section]
\newtheorem{lem}[thm]{Lemma}
\newtheorem{cor}[thm]{Corollary}
\newtheorem{prop}[thm]{Proposition}

\newtheorem{defn}[thm]{Definition}
\newtheorem{assum}[thm]{Assumption}

\begin{document}

\title{\bf Infinite-Time Mean Field FBSDEs and the Associated  Elliptic Master Equations    } 
\author{ Yongsheng Song$^*$ and Zeyu Yang\thanks{State Key Laboratory of Mathematical Sciences, Academy of Mathematics and Systems Science, Chinese Academy of Sciences, Beijing 100190, China, and School of Mathematical Sciences, University of Chinese Academy of Sciences, Beijing 100049, China. E-mails: yssong@amss.ac.cn (Y. Song), yangzeyu@amss.ac.cn (Z. Yang).}  }
\date{\today}

\maketitle

\begin{abstract}
We study discounted infinite-time mean field forward-backward stochastic
differential equations (FBSDEs) and the associated elliptic master equations,
which were introduced in~\cite{yang2025discounted} as mathematical tools for solving discounted
infinite-time mean field games.
We first establish stability
and flow properties, identify the backward decoupling field with the
spatial derivative of the equilibrium value function, and prove that
the value function is a viscosity solution of the master equation.
Under additional smoothness and quantitative stability assumptions, we
construct its Lions-derivative and the required higher-order derivatives
using variational results from a companion
paper~\cite{songyangregularity}. These results show that the value function
is a classical solution of the elliptic master equation. We also prove
that the value function is the unique displacement-monotone regular solution.
\end{abstract}

\no{\bf Keywords.}  discounted infinite-time 
  mean field games, infinite-time FBSDEs,  elliptic 
master equations, Lions-derivative.

\vfill\eject


\section{Introduction}
\label{sect-Introduction}
\setcounter{equation}{0} 

The study of mean field games was initiated independently by Lasry-Lions 
\cite{lasry2006jeux,lasry2006mean,lasry2007mean} and Huang-Malhamé-Caines \cite{huang2006large}, which is an analysis
of limit models for symmetric weakly interacting $(N+1)-$player differential games.
The parabolic master equation plays a crucial role in the analysis of mean field games, which was
 introduced by Lions in lectures \cite{lions2007cours}.
It describes the strategic interaction between a representative player and the collective environment. 
We refer the
reader to \cite{carmona2018probabilistic,cardaliaguet2019master,gangbo2022mean} for a comprehensive exposition on the subject.

Forward-backward stochastic differential equations (FBSDEs) also serve as a powerful tool for the study of mean field games.
The theory of general nonlinear backward stochastic differential equations (BSDEs) 
was pioneered by Pardoux and Peng \cite{pardoux1990adapted,pardoux2005backward} in the
early 1990s.
It's well known that BSDEs jointly with forward classical 
stochastic differential equations (SDEs) can provide a probabilistic interpretation of 
the solution of a class of second order semilinear parabolic and elliptic partial differential equations (PDEs),
which generalizes the Feynman-Kac formula to nonlinear cases and plays an important role in the field of stochastic optimization.
Here we refer the reader to \cite{peng1991probabilistic,pardoux1998backward,fuhrman2004infinite}.
For fully coupled infinite-time FBSDEs, \cite{peng2000infinite} established the well-posedness, 
and \cite{shi2020forward} investigated the problem in a broader solution space and connected it 
to quasilinear elliptic PDEs. Subsequently, \cite{bayraktar2023solvability} extended this framework to 
McKean-Vlasov FBSDEs and used it to solve infinite-time mean field control problems and mean field games.
This paper investigates the properties of  infinite-time mean field FBSDEs, a class of fully coupled  McKean-Vlasov infinite-time
FBSDEs which were introduced in \cite{yang2025discounted},
and uses them to construct solutions of the associated elliptic master
equation.  Our main objective is a complete passage from the equilibrium
FBSDE to the value function: viscosity existence under weighted
\(L_r^2\) assumptions, classical existence under explicit additional
regularity conditions, and uniqueness in a specified class of regular
solutions.

In the recent work \cite{yang2025discounted}, we proposed the discounted infinite-time  mean field games
and elliptic master equations, which extends the traditional framework to infinite-time case.
In our framework, a representative player interacts with a continuum of other players 
(also referred to as the population or social equilibrium).
Let $\mu_0\in \cP(\dbR)$ be the initial state distribution in the society.
For the player($x$), a representative player with initial state $x$,
the dynamic of its state  is given by
\begin{equation}
    \label{eq: constraint}
  X^{x,\beta}_t=x+\int_0^tb(X^{x,\beta}_s,\mu_s,\beta_s)\rmd s+B_t,
\end{equation}
where  $\beta$, an $\dbR$-valued progressively measurable stochastic process,
is the strategy of the player($x$)
and $(\mu_t)_{t\geq 0}$ is the population distribution flow start at $\mu_0$.
Denote by $\beta^x$ the optimal strategy of the representative player($x$) 
derived from the  stochastic optimization problem:
\begin{equation}
  \min_{\beta} J^\mu(\beta)\triangleq\mathbb{E}\left[\int_0^\infty \rme^{-rt}f(X_t^{x,\beta},\mu_t,\beta_t)\rmd t\right].
\end{equation}
Since the population consists of a multitude of individuals, the macroscopic distribution should satisfy
\begin{equation}
  \label{eq: relationship}
  \mu_t(\cdot)=\int \dbP(X_t^{x,\beta^x}\in \cdot)\mu_0(\rmd x).
\end{equation}
To solve this fixed-point problem,
 we introduce the following two systems of FBSDEs. 
The first one represents the state process of the social equilibrium, 
while the second denotes the state process of the representative player with initial state $x$.
\begin{equation}
  \label{eq: 11}
  \begin{cases}
    \rmd X_{t}^{\xi} = \partial_yH(X_{t}^{\xi},\mathcal{L} _{X_{t}^{\xi}},Y_{t}^{\xi}) \rmd t +\rmd B_{t}, \\
    \rmd Y_{t}^{\xi} = -\left[\partial_x {H} (X_{t}^{\xi},\mathcal{L} _{X_{t}^{\xi}}, Y_{t}^{\xi})-rY_{t}^{\xi}\right] \rmd t + Z_{t}^{\xi}\rmd B_{t}, \\
    X_0^{\xi}=\xi.
    \end{cases}
\end{equation}
\begin{equation}
  \label{eq: 12}
  \begin{cases}
    \rmd X_{t}^{x,\xi} = \partial_yH(X_{t}^{x,\xi},\mathcal{L} _{X_{t}^{\xi}},Y_{t}^{x,\xi}) \rmd t +\rmd B_{t}, \\
    \rmd Y_{t}^{x,\xi} = -\left[\partial_x {H} (X_{t}^{x,\xi},\mathcal{L} _{X_{t}^{\xi}}, Y_{t}^{x,\xi})-rY_{t}^{x,\xi}\right] \rmd t + Z_{t}^{x}\rmd B_{t}, \\
    X_0^x=x.
    \end{cases}
\end{equation}
Here $r>0$ is the discount factor and
\begin{equation}
\label{eq: H}
  H(x,\mu,y)\triangleq \min_{a\in \dbR} \left[ b(x,\mu,a)\cdot y+f(x,\mu,a) \right].
\end{equation}
Denote by \(\hat{\alpha}(x,\mu,y)\) the unique minimizer in
\eqref{eq: H}.
After further assuming that $f$ and $b$ possess good smoothness properties, 
we can obtain the following relationship:
\begin{equation}
  \pa_y H(x,\mu,y)=b(x,\mu,\hat{\alpha}(x,\mu,y)).
\end{equation}

In \cite{yang2025discounted}, 
we define $\mcv(x,\mu)\triangleq Y^{x,\xi}_0 $ with $\cL_\xi=\mu$ and
 demonstrate that \(\mcv(x,\mu)\) serves as a viscosity solution to the distribution-dependent PDE below:
\begin{equation}
  \label{eq: intro-pa}
  \begin{split}
  r \mcu(x,\mu)=&\partial_xH(x,\mu,\mcu(x,\mu))+ \partial_yH(x,\mu,\mcu(x,\mu))\cdot \partial_x\mcu(x,\mu) +
  \frac{1}{2}\partial_{xx}\mcu(x,\mu)\\&+\tilde{\mathbb{E} }
  \left[\frac{1}{2}\partial_{\tilde{x}}\partial_\mu \mcu(x,\mu,\tilde{\xi })+\partial_\mu \mcu(x,\mu,\tilde{\xi })
  \partial_y H(\tilde{\xi },\mu,\mcu(\tilde{\xi},\mu))\right].
\end{split}
\end{equation}
Here $\partial_x,\partial_{xx}$ are standard spatial derivatives, $\partial_{\mu},\partial_{\tilde{x}\mu}$ are $\cW_2$-Wasserstein derivatives, $\tilde \xi$ is a random variable with law $\mu$ and $\tilde \dbE$ is the expectation with respect to its law. 
We also introduce the elliptic master equation:
\begin{equation}
\label{eq: master}
\begin{split}
  r U(x,\mu)=&H(x,\mu,\partial_xU(x,\mu))+\frac{1}{2}\partial_{xx}U(x,\mu)\\&+\tilde{\mathbb{E} }
  \left[\frac{1}{2}\partial_{\tilde{x}}\partial_\mu U(x,\mu,\tilde{\xi })+\partial_\mu U(x,\mu,\tilde{\xi })
  \partial_y H(\tilde{\xi },\mu,\partial_xU(\tilde{\xi},\mu))\right].
\end{split}
\end{equation}
We have proved that if this master equation  admits a classical solution
with sufficient regularity, we can derive the following representation for FBSDEs (\ref{eq: 11}) and (\ref{eq: 12}):
\begin{equation}
  Y_t^\xi=\partial_xU(X_{t}^{\xi}, \mathcal{L} _{X_{t}^{\xi}}),\quad
  Z_t^\xi=\partial_{xx}U(X_{t}^{\xi}, \mathcal{L} _{X_{t}^{\xi}}),
\end{equation}
\begin{equation}
  Y_t^{x,\xi}=\partial_xU(X_{t}^{x,\xi}, \mathcal{L} _{X_{t}^{\xi}}),\quad
  Z_t^{x,\xi}=\partial_{xx}U(X_{t}^{x,\xi}, \mathcal{L} _{X_{t}^{\xi}}).
\end{equation}

In this paper, we make  above relation rigorous and proceed to study the properties of the  elliptic master equation.
We prove that solutions to FBSDEs (\ref{eq: 11}) and (\ref{eq: 12}) possess the flow property:
\begin{equation}
  X_t^{x,\xi}|_{x=\xi}=X_t^{\xi},\quad Y_t^{x,\xi}|_{x=\xi}=Y_t^{\xi},\quad ~\mbox{for $\rmd t\times \rmd \dbP$-a.e. $(t,\omega)$},
\end{equation}
and
\begin{equation}
  Y_0^{x,\xi}|_{x=\xi}=Y_0^{\xi},\quad ~\mbox{for $\rmd \dbP$-a.e. $\omega$}.
\end{equation}
Furthermore, we prove that the value function of the representative player
\begin{equation}
  \label{eq: value}
V(x,\mu)=  \mathbb{E}\bigg[\int_{0}^{+\infty}\rme^{-rt}f\big(X_{t}^{x,\xi},\mathcal{L}_{X_{t}^{\xi}},
\hat{\alpha}(X_{t}^{x,\xi},\mathcal{L}_{X_{t}^{\xi}}, Y_{t}^{x,\xi})        \big)\rmd t\bigg]
\end{equation}
satisfies the relationship
\begin{equation}
  \pa_xV(x,\mu)=Y_0^{x,\xi},
\end{equation}
and constitutes a 
viscosity solution to elliptic master equation (\ref{eq: master}).
To obtain a classical solution, we impose additional
smoothness and quantitative stability assumptions and use the
variational results for \(\partial_xV(x,\mu)\) from the companion
paper \cite{songyangregularity}. The usual finite-time approach to
Lions differentiability does not apply in our \(L_r^2\) framework.
We therefore directly construct
\(\partial_\mu\partial_xV(x,\mu,\tilde x)\), following
\cite{gangbo2022mean,mou2024wellposedness}.
Higher-order variational equations require
higher-moment estimates, obtained using the Riccati results of
Guatteri and Tessitore \cite{guatteri2008backward}.
Using these results and the discounted representation of \(V(x,\mu)\),
we construct \(\partial_\mu V(x,\mu,\tilde x)\), differentiate it
in \(x\) and \(\tilde x\), and establish joint continuity and
growth bounds.
These properties allow us to take \(V(x,\mu)\) as a test function
in the viscosity inequalities and verify that it is a classical
solution of \eqref{eq: master}.

Finally, we address the uniqueness of classical solutions to the
elliptic master equation. An example in \cite{yang2025discounted}
illustrates that classical solutions to the elliptic master equation
are not unique; however, at most one of them can be used to construct
a Nash equilibrium. We introduce the displacement monotonicity
condition (see \cite{gangbo2022mean}) for the solutions, which states
that for any square integrable random variables \(\xi_1,\xi_2\),
the following inequality holds:
        \begin{equation}
            \mathbb{E}\left[ 
                (\xi_1 - \xi_2) \left( 
                    \partial_x U(\xi_1, \mathcal{L}_{\xi_1}) - \partial_x U(\xi_2, \mathcal{L}_{\xi_2}) 
                \right) 
            \right] \geq 0.
        \end{equation}
We establish uniqueness in the class of regular solutions, which
combines this monotonicity with additional differentiability and
growth conditions. The value function \(V(x,\mu)\) belongs to this class
and is therefore the unique regular solution.

Section~\ref{sec:setting} introduces the weighted spaces and measure
calculus.  Section~\ref{sec: Properties} establishes stability and
flow properties of the mean field FBSDEs.
Section~\ref{sec:viscosity} proves the spatial derivative identity
and viscosity existence.
Section~\ref{sec: main higher regularity} recalls the variational
results from \cite{songyangregularity} and proves the regularity
of \(V(x,\mu)\).
Section~\ref{sec: main classical master} verifies that \(V(x,\mu)\)
is a classical solution of the master equation and establishes
uniqueness in the regular class.

\section{Preliminaries} 
\label{sec:setting}
\setcounter{equation}{0}

We will use a filtered probability space $(\Omega, \mathcal{F} ,\mathbb{P},\mathbb{F} )$
endowed
with a Brownian motion $B$.
Its filtration $\mathbb{F} \triangleq (\mathcal{F} _t)_{t\ge 0}$ is
augmented
by all $\mathbb{P}$-null sets and a sufficiently rich sub-$\sigma$-algebra $\mathcal{F}_0$ independent
of $B$, such that it
can support any measure on $ \mathbb{R} $ with finite second moment.

Let $( \O', \mathcal{F}' , \dbP',\dbF')$ be a copy of the filtered probability space $(\Omega, \mathcal{F} ,\mathbb{P},\mathbb{F})$
 with corresponding Brownian motion $B'$, define the larger filtered probability space by
\begin{equation}
\tilde \O \triangleq \O\times \O' ,\q \tilde{ \mathcal{F}} \triangleq \mathcal{F}\otimes \mathcal{F}'\q \tilde\dbF = \{\tilde \cF_t\}_{t\ge 0} \triangleq \{\cF_t \otimes  \cF'_t\}_{t\ge 0},\q \tilde \dbP \triangleq \dbP\otimes \dbP',\q \tilde \dbE\triangleq \dbE^{\tilde \dbP}.
\end{equation}
Throughout the paper we will use the probability space $(\Omega, \mathcal{F} ,\mathbb{P},\mathbb{F})$. However, when we deal with the distribution-dependent master equation, independent copies of random variables or processes are needed. Then we will tacitly use their extensions to the larger space $(\tilde \O,\tilde{\mathcal{F}},  \tilde \dbP,\tilde \dbF)$.

Let $\cP \triangleq\cP(\dbR)$ be the set of all probability measures on $\mathbb R$ and let $\cP_p(p\ge1)$ denote the set of $\mu\in \cP$ with finite $p$-th moment.
For any sub-$\si$-field $\cG\subset \cF$ and $\mu\in \cP_p$, we define $\dbL^p(\cG)$ to be the set of $\dbR$-valued, $\cG$-measurable, and $p$-integrable random variables $\xi$ , and $\dbL^p(\cG;\mu)$ to be the set of $\xi\in \dbL^p(\cG)$ such that the law $\cL_\xi=\mu$ .
For any $\mu,\nu\in \cP_p$, we define the $\cW_p$--Wasserstein distance between them as follows: 
\beaa
\cW_p(\mu, \nu) \triangleq \inf\Big\{\big(\dbE[|\xi-\eta|^p]\big)^{1/p}: \mbox{for all $\xi\in \dbL^p(\cF; \mu)$, $\eta\in \dbL^p(\cF; \nu)$}\Big\}.
\eeaa

Due to our interest in discounted infinite-time mean field games, for any $K\in \mathbb{R} $, we denote by $L^2_K(t_0, \infty, \mathbb{R} )$
 the Hilbert space of all $\mathbb{R} -$valued adapted stochastic process ($v_t$) start from $t_0$ such that
\begin{equation}\mathbb{E}\left[\int_{t_0}^\infty \rme^{-Kt}|v_t|^2dt\right]<+\infty.\end{equation}
To simplify, we set $L^2_K \triangleq L^2_K(0, \infty,\mathbb{R} )$ and
define the exponentially weighted $L^2$ norm $\lVert \cdot \rVert_K$ by
\begin{equation}
  \Vert v\Vert _K^2\triangleq \dbE\left[\int_0^\infty \rme^{-Kt}|v_t|^2\rmd t \right].
\end{equation}
For each $\mathcal{F} _0$-measurable square integrable random
variable $\xi$ , we consider the following infinite-time FBSDE:
\begin{equation}
  \label{eq: fbsde}
  \begin{cases}
    \rmd X_t=G(t,X_t,Y_t,\mathcal{L}_{X_t})dt+ \rmd B_t, \\
    \rmd Y_t=-F(t,X_t,Y_t,\mathcal{L}_{X_t})dt+Z_t\rmd B_t, \\
    X_0=\xi.
    \end{cases}
\end{equation}
Here $G,F:\dbR^+\times \dbR^2\times \cP_2\to\dbR$ are two measurable functions, and satisfy
 the following assumptions:
\begin{assum}
  \label{assum: fbsde}
  (i) There exists a positive constant $\ell$ such that for any $x,x',y,y'\in \mathbb{R}, \mu,\mu'\in \cP_2$
  \begin{equation}
    \begin{aligned}
      |G(t,x,y,\mu) & -G(t,x^{\prime},y^{\prime},\mu^{\prime})|+|F(t,x,y,\mu)-F(t,x^{\prime},y^{\prime},\mu^{\prime})| \\
       & \leq \ell(|x-x^{\prime}|+|y-y^{\prime}|+\mathcal{W}_{2}(\mu,\mu^{\prime})).\quad\mathrm{a.s.}
      \end{aligned}
  \end{equation}

(ii) There exist constants $0<K<2\kappa$ such that for any $t\geq0$ and any square
integrable random variables $X,X',Y,Y'$
\begin{equation}
  \label{eq: condition}
  \begin{aligned}
    &\mathbb{E}\left[-K\hat{X}\hat{Y}-\hat{X}(F(t,U)-F(t,U^{\prime}))+\hat{Y}(G(t,U)-G(t,U^{\prime}))\right] \\
      &\leq-\kappa\mathbb{E}\left[\hat{X}^{2}+\hat{Y}^{2}\right],
    \end{aligned}
\end{equation}
where $ \hat{X}=X-X^{\prime},\hat{Y}=Y-Y^{\prime}$ and $U=(X,Y,\mathcal{L}_X),U'=(X',Y',\mathcal{L}_{X'}) .$
\end{assum}
The following lemma states the existence and uniqueness of a solution to FBSDE(\ref{eq: fbsde}).
For the detailed proof, we refer the reader to \cite[Theorem~2.1]{bayraktar2023solvability}.
\begin{lem}
  \label{lem: fbsde}
  Under Assumption \ref{assum: fbsde}, the FBSDE (\ref{eq: fbsde}) admits a unique solution in $L_K^2$.
\end{lem}

We introduce the Wasserstein space and differential calculus on Wasserstein space. 
For a $\cW_2$-continuous functions $U: \cP_2 \to \dbR$, its $\cW_2$-Wasserstein derivatives\cite{carmona2018probabilistic}(also called Lions-derivative), takes the form $\pa_\mu U: (\mu,\tilde x)\in \cP_2\times \dbR\to \dbR$ and satisfies: 
\bea
\label{pamu}
U(\cL_{\xi +  \eta}) - U(\mu) = \dbE\big[\langle \pa_\mu U(\mu, \xi), \eta \rangle \big] + o(\|\eta\|_2), \ \forall\ \xi\in\mathbb L^2(\mathcal{F};\mu),\eta\in\mathbb L^2(\mathcal{F}).
\eea
Let $C^0(\cP_2)$ denote the set of $\cW_2$-continuous functions $U:\cP_2\to\dbR$.
For $C^1(\cP_2)$, we mean the space of functions $U\in C^0(\cP_2)$ such that $\pa_\mu U$ exists and is continuous on $ \cP_2\times \dbR$, which is uniquely determined by \reff{pamu}.
Let $C^{2,1}(\dbR\times\cP_2)$ denote the set of continuous functions $U:\dbR\times\cP_2\to\dbR$ such that $\pa_xU,\pa_{xx}U$ exist and are joint continuous on $\dbR\times \cP_2$, $\pa_\mu U,\pa_{\tilde x\mu}U$ exist and are continuous on $\dbR\times\cP_2\times\dbR$.

Finally, we consider the space 
$\Th\triangleq [0, T]\times \dbR \times \cP_2$ for some $T>0$,
 and let $C^{1,2,1}(\Th)$ denote the set of  continuous functions $U: \Th\to \dbR$ which has the following continuous derivatives: 
$\pa_t U$, $\pa_x U$, $\pa_{xx} U$, $\pa_\mu U$, $\pa_{\tilde x\mu} U.$ 
One crucial property of functions $U\in C^{1,2,1}(\Th)$ is  the  It\^{o}'s  formula\cite[Vol.~I, Proposition~5.102]{carmona2018probabilistic}. For $i=1,2$, let $\rmd X^i_t \triangleq b^i_t \rmd t + \si^i_t \rmd B_t ,$ where $b^i:[0,T]\times\Omega\to\mathbb R$ and $\sigma^i:[0,T]\times\Omega\to\mathbb R$ are $\dbF$-progressively 
measurable and bounded (for simplicity), and $\rho_t\triangleq \cL_{X^2_t}$.
Suppose that for every compact subset $\cK\subset \dbR\times \cP_2$, it holds that:
\begin{equation}
  \sup_{(t,x,\mu)\in [0,T]\times\cK} 
    \int_\dbR \left( \left\lvert \pa_\mu U(t,x,\mu,\tilde{x})\right\rvert^2 
    +\left\lvert \pa_{\tilde{x}\mu}U(t,x,\mu,\tilde{x})\right\rvert^2 \right) \rmd \mu(\tilde{x})
  <\infty.
\end{equation}
We have
\begin{equation}
  \label{eq:ito}
\begin{split}
  \rmd U(t, X^1_t, \rho_t) = & \Big[\pa_t U + \pa_x U\cdot b^1_t + \frac{1}{2} \pa_{xx} U (\si_t^1 )^2\Big](t, X^1_t, \rho_t) \rmd t \\
&+\Big(\tilde \dbE_{\cF_t}\big[\pa_\mu U(t,X^1_t,\rho_t,\tilde X^2_t) (\tilde b^{2}_t)+\frac{1}{2} \pa_{\tilde x}\pa_\mu U(t, X^1_t, \rho_t, \tilde X^2_t)(\tilde \si_t^2 )^2\big]\Big) \rmd t\\
&+\pa_xU(t,X^1_t,\rho_t)\si_t^1\rmd B_t  .
\end{split}
\end{equation}
Here  $\tilde \dbE_{\cF_t}$ is the conditional expectations given $\cF_t$ corresponding to the probability measure $\tilde \dbP$
and the process $(\tilde X^2_t,\tilde b^{2}_t,\tilde\si_t^2)_{0\leq t\leq T}$
is a copy of the process $( X^2_t, b^{2}_t,\si_t^2)_{0\leq t\leq T}$,
on a copy $( \O', \mathcal{F}' , \dbP',\dbF')$ of the probability space $(\Omega, \mathcal{F} ,\mathbb{P},\mathbb{F} )$.

\section{Properties of mean field FBSDEs}
\label{sec: Properties}
\setcounter{equation}{0}

In this section, for some $\xi\in\dbL^2(\cF_0)$, we investigate the following mean field FBSDEs:
\begin{equation}
  \label{eq: mf1}
  \begin{cases}
    \rmd X_{t}^{\xi} = \partial_yH(X_{t}^{\xi},\mathcal{L} _{X_{t}^{\xi}},Y_{t}^{\xi}) \rmd t +\rmd B_{t}, \\
    \rmd Y_{t}^{\xi} = -\left[\partial_x {H} (X_{t}^{\xi},\mathcal{L} _{X_{t}^{\xi}}, Y_{t}^{\xi})-rY_{t}^{\xi}\right] \rmd t + Z_{t}^{\xi}\rmd B_{t}, \\
    X_0^{\xi}=\xi,
    \end{cases}
\end{equation}
\begin{equation}
  \label{eq: mf2}
  \begin{cases}
    \rmd X_{t}^{x,\xi} = \partial_yH(X_{t}^{x,\xi},\mathcal{L} _{X_{t}^{\xi}},Y_{t}^{x,\xi}) \rmd t +\rmd B_{t}, \\
    \rmd Y_{t}^{x,\xi} = -\left[\partial_x {H} (X_{t}^{x,\xi},\mathcal{L} _{X_{t}^{\xi}}, Y_{t}^{x,\xi})-rY_{t}^{x,\xi}\right] \rmd t + Z_{t}^{x}\rmd B_{t}, \\
    X_0^{x,\xi}=x.
    \end{cases}
\end{equation}

In \cite{yang2025discounted}, 
we provided sufficient conditions on $b(x,\mu,a)$ and $f(x,\mu,a)$ to ensure the above two FBSDEs admit  unique solutions.
To obtain further properties of these solutions , we require the following assumptions on $H$:

\begin{assum}
  \label{assum: H}
\noindent(i)
$H(x,\mu,y)$ is jointly continuous and all second-order partial derivatives exist.

\noindent  (ii) $\pa_yH(x,\mu,y)$ and $\pa_xH(x,\mu,y)$ are Lipschitz continuous in $(x,\mu,y)$.
  More specifically, there exists a constant $\ell>0$ such that
  \begin{equation}
    \begin{aligned}
      & |\pa_yH(x,\mu,y)-\pa_yH(x',\mu',y')|\le \ell\left(|x-x'|+|y-y'|+\cW_2(\mu,\mu')\right),\\
      & |\pa_xH(x,\mu,y)-\pa_xH(x',\mu',y')|\le \ell\left(|x-x'|+|y-y'|+\cW_2(\mu,\mu')\right).
    \end{aligned}
  \end{equation}

\noindent (iii)There exist constants $\kappa,C_0>0$ , such that $C_0+r/2<\kappa$ and
\begin{equation}
  \begin{aligned}
  &-(x-x')\left[\pa_xH(x,\mu,y)-\pa_xH(x',\mu',y') \right]
  +(y-y')\left[\pa_yH(x,\mu,y)-\pa_yH(x',\mu',y') \right]\\
  &\le -\kappa\left(|x-x'|^2+|y-y'|^2\right)+C_0\cW_2^2(\mu,\mu').    
  \end{aligned}
\end{equation}
\end{assum}

\begin{thm}
  Under Assumption \ref{assum: H}, both FBSDE (\ref{eq: mf1}) and (\ref{eq: mf2})
  have unique solutions in $L_r^2$.
\end{thm}
\proof
For FBSDE (\ref{eq: mf1}), it's clear that the conditions in Assumption~\ref{assum: fbsde}(i)--(ii) are satisfied.
Taking four arbitrary square
integrable random variables $X,X',Y,Y'$, we have
\begin{equation}
\begin{aligned}
\mathbb{E}\biggl[
  & -r\hat{X}\cdot \hat{Y}\\
  & - \hat{X}\left[\partial_x H(X,\cL_{X},Y) - \partial_x H(X',\cL_{X'},Y') - r\hat{Y} \right] \\
  & + \hat{Y}\left[\partial_y H(X,\cL_{X},Y) - \partial_y H(X',\cL_{X'},Y') \right] \biggr] \\
\leq& -\kappa\mathbb{E}\left[\hat{X}^{2} + \hat{Y}^{2}\right] + C_0\mathcal{W}_2^2(\cL_{X},\cL_{X'}) \\
\leq& -(\kappa - C_0)\mathbb{E}\left[\hat{X}^{2} + \hat{Y}^{2}\right],
\end{aligned}
\end{equation}
where $ \hat{X}\triangleq X-X^{\prime},\hat{Y}\triangleq Y-Y^{\prime}$.
Since $\kappa-C_0>r/2$, the condition in Assumption~\ref{assum: fbsde}(ii) is satisfied.
Then, FBSDE (\ref{eq: mf1}) has a unique solution in $L_r^2$.

After solving FBSDE (\ref{eq: mf1}), we substitute its solution $\cL_{X_t^{\xi}}$ into FBSDE (\ref{eq: mf2}), 
and it can similarly be shown that there exists a unique solution to FBSDE (\ref{eq: mf1}).

\qed

The following proposition informs us that the solution of the mean field FBSDE
 exhibits favorable continuous dependence on the initial value, 
which is of paramount importance for our subsequent research.

\begin{prop}
  \label{prop: con}
  For FBSDEs (\ref{eq: mf1}) and (\ref{eq: mf2}), assuming all conditions in Assumption \ref{assum: H} are satisfied,
   we have, for any $x,x'\in\mathbb{R} $ and $\xi,\xi'\in \dbL^2(\cF_0)$, there exists a constant $C>0$, such that
\begin{equation}
  \dbE\left[|Y_0^{\xi}-Y_0^{\xi'} |^2\right]+\left\lVert X^{\xi}-X^{\xi'}\right\rVert _r^2 
  +\left\lVert Y^{\xi}-Y^{\xi'}\right\rVert _r^2\leq C \dbE\left[|{\xi}-{\xi'} |^2\right]
\end{equation}
and
\begin{equation}
 |Y_0^{x,\xi}-Y_0^{x',\xi'} |^2+\left\lVert X^{x,\xi}-X^{x',\xi'}\right\rVert _r^2 
  +\left\lVert Y^{x,\xi}-Y^{x',\xi'}\right\rVert _r^2\leq C \left(|x-x'|^2+\dbE\left[|{\xi}-{\xi'} |^2\right]\right) .
\end{equation}
\end{prop}
\proof
Set 
\begin{equation}
  \begin{split}
    &\hat{X}^{\xi}=X^{\xi}-X^{\xi'},\quad \hat{Y}^{\xi}=Y^{\xi}-Y^{\xi'},
    \quad \hat{Z}^{\xi}=Z^{\xi}-Z^{\xi'},\\
    &\hat{X}^{x,\xi}=X^{x,\xi}-X^{x',\xi'},\quad \hat{Y}^{x,\xi}=Y^{x,\xi}-Y^{x',\xi'}
    ,\quad \hat{Z}^{x,\xi}=Z^{x,\xi}-Z^{x',\xi'}.
  \end{split}
\end{equation}
$C_1,C_2,C_3,C_4,C_5,C_6,C_7$ appeared in the following proof are all positive constants.

Applying It\^{o}'s formula to $\rme^{-rt}|Y^{\xi}_t-Y^{\xi'}_t|^2$,
and taking a sequence of $T_i\to\infty$ such that
\begin{equation}
  \dbE\left[\rme^{-rt}|Y^{\xi}_{T_i}-Y^{\xi'}_{T_i}|^2\right]\rightarrow 0,
\end{equation}
we get that
\begin{equation}
\label{eq: 4.1}
\begin{aligned}
\mathbb{E}\left[\left|Y_0^{\xi} - Y_0^{\xi'}\right|^2\right] 
=&\mathbb{E}\int_0^\infty \mathrm{e}^{-rt} \biggl[ r\left|Y_t^{\xi} - Y_t^{\xi'}\right|^2 \\
& \quad + 2\hat{Y}^{\xi} 
 \cdot \left( \partial_x H\left(X_t^{\xi},\mathcal{L}_{X_t^{\xi}},Y_t^{\xi}\right) - \partial_x H\left(X_t^{\xi'},\mathcal{L}_{X_t^{\xi'}},Y_t^{\xi'}\right) - r\hat{Y}^{\xi} \right) \\
& \quad - \left|Z_t^{\xi} - Z_t^{\xi'}\right|^2 \biggr] \mathrm{d}t\\
\leq & C_1\left(\left\lVert X^{\xi}-X^{\xi'}\right\rVert _r^2 
  +\left\lVert Y^{\xi}-Y^{\xi'}\right\rVert _r^2\right).
\end{aligned}
\end{equation}
Applying It\^{o}'s formula to $\rme^{-rt}\hat{X}^{\xi}\hat{Y}^{\xi}$, we get
\begin{equation}
  \label{eq: 4.2}
\begin{aligned}
-\mathbb{E}\left[\hat{X}^{\xi}_0\hat{Y}^{\xi}_0\right] 
= & \ \mathbb{E}\int_0^\infty \mathrm{e}^{-rt} \biggl[ -r\hat{X}^{\xi}_t\hat{Y}^{\xi}_t \\
& + \hat{Y}^{\xi}_t \left( \partial_y H(X^{\xi}_t,\mathcal{L}_{X^{\xi}_t},Y^{\xi}_t) - \partial_y H(X^{\xi'}_t,\mathcal{L}_{X^{\xi'}_t},Y^{\xi'}_t) \right) \\
& - \hat{X}^{\xi}_t \left( \partial_x H(X^{\xi}_t,\mathcal{L}_{X^{\xi}_t},Y^{\xi}_t) - \partial_x H(X^{\xi'}_t,\mathcal{L}_{X^{\xi'}_t},Y^{\xi'}_t) - r\hat{Y}^{\xi}_t \right) \biggr] \mathrm{d}t\\
\leq &-\frac{r}{2}\left(\left\lVert X^{\xi}-X^{\xi'}\right\rVert _r^2 
  +\left\lVert Y^{\xi}-Y^{\xi'}\right\rVert _r^2\right).
\end{aligned}
\end{equation}
Then, we have
\begin{equation}\label{eq: 4.3}\begin{aligned}
  \left\lVert X^{\xi}-X^{\xi'}\right\rVert _r^2 
  +\left\lVert Y^{\xi}-Y^{\xi'}\right\rVert _r^2\leq& \frac{2}{r}\mathbb{E}\left[\hat{X}^{\xi}_0\hat{Y}^{\xi}_0\right]\\
\leq& \frac{1}{2C_1}\mathbb{E}\left[\left|Y_0^{\xi} - Y_0^{\xi'}\right|^2\right] +\frac{32C_1}{r^2}\dbE\left[|{\xi}-{\xi'} |^2\right].
\end{aligned}
\end{equation}
Combining (\ref{eq: 4.1}) and (\ref{eq: 4.3}), we get that
\begin{equation}
  \mathbb{E}\left[\left|Y_0^{\xi} - Y_0^{\xi'}\right|^2\right]\le C_2\dbE\left[|{\xi}-{\xi'} |^2\right].
\end{equation}
Substituting this into (\ref{eq: 4.3}), we obtain
\begin{equation}
  \left\lVert X^{\xi}-X^{\xi'}\right\rVert _r^2 
  +\left\lVert Y^{\xi}-Y^{\xi'}\right\rVert _r^2\le C_3\dbE\left[|{\xi}-{\xi'} |^2\right].
\end{equation}

Next, we employ the same approach to  FBSDE (\ref{eq: mf2}).
Applying It\^{o}'s formula to $\rme^{-rt}|Y^{x,\xi}_t-Y^{x',\xi'}_t|^2$, we get that
\begin{equation}
  \label{eq: 4.4}
\begin{aligned}
\mathbb{E}\left[\left|Y_0^{x,\xi} - Y_0^{x',\xi'}\right|^2\right] 
=&\mathbb{E}\int_0^\infty \mathrm{e}^{-rt} \biggl[ r\left|Y_t^{x,\xi} - Y_t^{x',\xi'}\right|^2 \\
& \quad + 2\hat{Y}^{x,\xi} 
 \cdot \left( \partial_x H\left(X_t^{x,\xi},\mathcal{L}_{X_t^{\xi}},Y_t^{x,\xi}\right) - \partial_x H\left(X_t^{x',\xi'},\mathcal{L}_{X_t^{\xi'}},Y_t^{x',\xi'}\right) - r\hat{Y}^{x,\xi} \right) \\
& \quad - \left|Z_t^{x,\xi} - Z_t^{x',\xi'}\right|^2 \biggr] \mathrm{d}t\\
\leq & C_4\left(\left\lVert X^{x,\xi}-X^{x',\xi'}\right\rVert _r^2 
  +\left\lVert Y^{x,\xi}-Y^{x',\xi'}\right\rVert _r^2
  +\left\lVert X^{\xi}-X^{\xi'}\right\rVert _r^2\right).
\end{aligned}
\end{equation}
Applying It\^{o}'s formula to $\rme^{-rt}\hat{X}^{x,\xi}\hat{Y}^{x,\xi}$, we can get that
\begin{equation}
\begin{aligned}
-\mathbb{E}\left[\hat{X}^{x,\xi}_0\hat{Y}^{x,\xi}_0\right] 
= & \ \mathbb{E}\int_0^\infty \mathrm{e}^{-rt} \biggl[ -r\hat{X}^{x,\xi}_t\hat{Y}^{x,\xi}_t \\
& + \hat{Y}^{x,\xi}_t \left( \partial_y H(X^{x,\xi}_t,\mathcal{L}_{X^{\xi}_t},Y^{x,\xi}_t) - \partial_y H(X^{x',\xi'}_t,\mathcal{L}_{X^{\xi'}_t},Y^{x',\xi'}_t) \right) \\
& - \hat{X}^{x,\xi}_t \left( \partial_x H(X^{x,\xi}_t,\mathcal{L}_{X^{\xi}_t},Y^{x,\xi}_t) - \partial_x H(X^{x',\xi'}_t,\mathcal{L}_{X^{\xi'}_t},Y^{x',\xi'}_t) - r\hat{Y}^{x,\xi}_t \right) \biggr] \mathrm{d}t\\
\leq &-\frac{r}{2}\left(\left\lVert X^{x,\xi}-X^{x',\xi'}\right\rVert _r^2 
  +\left\lVert Y^{x,\xi}-Y^{x',\xi'}\right\rVert _r^2\right)
  +C_0 \left\lVert X^{\xi}-X^{\xi'}\right\rVert _r^2.
\end{aligned}
\end{equation}
Then, we have 
\begin{equation}\label{eq: 4.5}\begin{aligned}
  \left\lVert X^{x,\xi}-X^{x',\xi'}\right\rVert _r^2 
  +\left\lVert Y^{x,\xi}-Y^{x',\xi'}\right\rVert _r^2
\leq& \frac{1}{2C_4}\mathbb{E}\left[\left|Y_0^{x,\xi} - Y_0^{x',\xi'}\right|^2\right] +C_5\left(|x-x'|^2+\dbE\left[|{\xi}-{\xi'} |^2\right]\right).
\end{aligned}
\end{equation}
Substituting it into (\ref{eq: 4.4}), we have
\begin{equation}
\begin{aligned}
\mathbb{E}\left[\left|Y_0^{x,\xi} - Y_0^{x',\xi'}\right|^2\right] 
\leq& C_6\left(|x-x'|^2+\dbE\left[|{\xi}-{\xi'} |^2\right]\right).
\end{aligned}
\end{equation}
At last we come back to (\ref{eq: 4.5}), we get that
\begin{equation}\begin{aligned}
  \left\lVert X^{x,\xi}-X^{x',\xi'}\right\rVert _r^2 
  +\left\lVert Y^{x,\xi}-Y^{x',\xi'}\right\rVert _r^2
\leq& C_7\left(|x-x'|^2+\dbE\left[|{\xi}-{\xi'} |^2\right]\right).
\end{aligned}
\end{equation}
Now, we complete the proof.
\qed

Benefiting from the continuous dependence on initial values of the solutions to the mean field FBSDEs, 
we can demonstrate the flow property. 
This precisely reflects the essence of the mean field game: 
the aggregation of innumerable representative players constitutes the equilibrium state of the system.

\begin{thm}
  Let Assumption \ref{assum: H} hold. For any $\xi\in\dbL^2(\cF_0)$, we have
\begin{equation}
  X_t^{x,\xi}|_{x=\xi}=X_t^{\xi},\quad Y_t^{x,\xi}|_{x=\xi}=Y_t^{\xi},\quad ~\mbox{for $\rmd t\times \rmd \dbP$-a.e. $(t,\omega)$},
\end{equation}
and
\begin{equation}
  Y_0^{x,\xi}|_{x=\xi}=Y_0^{\xi},\quad ~\mbox{for $\rmd \dbP$-a.e. $\omega$}.
\end{equation}
\end{thm}
\proof
We prove this theorem in two steps.

\emph{Step 1.} We first assume $\xi$ is discrete, that is 
\begin{equation}
  \xi=\sum_{i=1}^{n}x_i I_{A_i},
\end{equation}
where $x_i\in \dbR$ are constants and $\{A_i\}\in \cF_0$ is a partition of $\Omega$.
For each $i\in\{1,2,\cdots,n\}$, we can solve the following FBSDE in $L_r^2$:
\begin{equation}
  \begin{cases}
    \rmd X_{t}^{x_i,\xi} = \partial_yH(X_{t}^{x_i,\xi},\mathcal{L} _{X_{t}^{\xi}},Y_{t}^{x_i,\xi}) \rmd t +\rmd B_{t}, \\
    \rmd Y_{t}^{x_i,\xi} = -\left[\partial_x {H} (X_{t}^{x_i,\xi},\mathcal{L} _{X_{t}^{\xi}}, Y_{t}^{x_i,\xi})-rY_{t}^{x_i,\xi}\right] \rmd t + Z_{t}^{x_i}\rmd B_{t}, \\
    X_0^{x_i}=x_i.
    \end{cases}
\end{equation}
Now, we define 
\begin{equation}
  X_t(\omega)\triangleq \sum_{i=1}^{n}X_t^{x_i,\xi}(\omega)\cdot I_{A_i}(\omega),
\end{equation}
 and similarly define $Y$ and $Z$ in the same manner.
Then, we know $(X,Y,Z)\in L_r^2$. 

Multiplying both sides of the individual integral equations by \( I_{A_i} \) and summing over \( i \), we derive the aggregated equations. For the forward process:
\begin{equation}
\begin{split}
X_t & = \sum_{i=1}^n X_t^{x_i,\xi} I_{A_i} \\
& = \sum_{i=1}^n \left( x_i + \int_0^t \partial_y H(X_s^{x_i,\xi}, \mathcal{L}_{X_s^{\xi}}, Y_s^{x_i,\xi})  \mathrm{d}s + B_t \right) I_{A_i} \\
& = \sum_{i=1}^n x_i I_{A_i} + \int_0^t \sum_{i=1}^n \partial_y H(X_s^{x_i,\xi}, \mathcal{L}_{X_s^{\xi}}, Y_s^{x_i,\xi}) I_{A_i}  \mathrm{d}s + \sum_{i=1}^n I_{A_i} B_t \\
& = \xi + \int_0^t \partial_y H(X_s, \mathcal{L}_{X_s^{\xi}}, Y_s)  \mathrm{d}s + B_t.
\end{split}
\end{equation}
 For the backward process:
\begin{equation}
\begin{split}
Y_t & = \sum_{i=1}^n Y_t^{x_i,\xi} I_{A_i} \\
& = \sum_{i=1}^n \left( Y_T^{x_i,\xi} + \int_t^T \left[ \partial_x H(X_s^{x_i,\xi}, \mathcal{L}_{X_s^{\xi}}, Y_s^{x_i,\xi}) - r Y_s^{x_i,\xi} \right]  \mathrm{d}s - \int_t^T Z_s^{x_i,\xi}  \mathrm{d}B_s \right) I_{A_i} \\
& = \sum_{i=1}^n Y_T^{x_i,\xi} I_{A_i} + \int_t^T \sum_{i=1}^n \left[ \partial_x H(X_s^{x_i,\xi}, \mathcal{L}_{X_s^{\xi}}, Y_s^{x_i,\xi}) - r Y_s^{x_i,\xi} \right] I_{A_i}  \mathrm{d}s - \int_t^T \sum_{i=1}^n Z_s^{x_i,\xi} I_{A_i}  \mathrm{d}B_s \\
& = Y_T + \int_t^T \left[ \partial_x H(X_s, \mathcal{L}_{X_s^{\xi}}, Y_s) - r Y_s \right]  \mathrm{d}s - \int_t^T Z_s  \mathrm{d}B_s.
\end{split}
\end{equation}
Thus, the aggregated processes \( (X, Y, Z) \) satisfy the following FBSDE:
\begin{equation}
  \begin{cases}
    \mathrm{d} X_{t} = \partial_y H(X_{t},\mathcal{L} _{X_{t}^{\xi}},Y_{t}) \mathrm{d} t + \mathrm{d} B_{t}, \\
    \mathrm{d} Y_{t} = -\left[\partial_x {H} (X_{t},\mathcal{L} _{X_{t}^{\xi}}, Y_{t})-rY_{t}\right] \mathrm{d} t + Z_{t} \mathrm{d} B_{t}, \\
    X_0 = \xi.
    \end{cases}
\end{equation}
We have proved that, given $\mathcal{L} _{X_{t}^{\xi}}$, the above FBSDE has a unique solution in $L_r^2$.
So 
\begin{equation}
  \left\lVert X-X^{\xi}\right\rVert _r =\left\lVert Y-Y^{\xi}\right\rVert _r =0
\end{equation}
and
\begin{equation}
  \dbE\left[|Y_0-Y_0^\xi|^2\right]=0.
\end{equation}
Then, we have
\begin{equation}
  X_t^{x,\xi}(\omega)|_{x=\xi(\omega)}=X_t^{\xi}(\omega),
  \quad Y_t^{x,\xi}(\omega)|_{x=\xi(\omega)}=Y_t^{\xi}(\omega),\quad ~\mbox{for $\rmd t\times \rmd \dbP$-a.e. $(t,\omega)$},
\end{equation}
and
\begin{equation}
  Y_0^{x,\xi}|_{x=\xi}=Y_0^{\xi},\quad ~\mbox{for $\rmd \dbP$-a.e. $\omega$}.
\end{equation}

\emph{Step 2.} In the general case, let $\xi_n\in \dbL^2(\cF_0)$ be a discrete approximation of $\xi$
such that
 $ \dbE\left[ \left\lvert \xi_n-\xi\right\rvert ^2  \right]\rightarrow 0$.
By Step 1, we have 
\begin{equation}
  X_t^{x,\xi_n}(\omega)|_{x=\xi_n(\omega)}=X_t^{\xi_n}(\omega),
  \quad Y_t^{x,\xi_n}(\omega)|_{x=\xi_n(\omega)}=Y_t^{\xi_n}(\omega)\quad ~\mbox{for $\mathrm{d}t\times \mathrm{d}\mathbb{P}$-a.e. $(t,\omega)$}.
\end{equation}
Applying Proposition \ref{prop: con}, we have
\begin{equation}
  \left\lVert X^{\xi_n}-X^{\xi}\right\rVert _r\to 0,\quad
  \left\lVert Y^{\xi_n}-Y^{\xi}\right\rVert _r\to 0.
\end{equation}
On the other hand,
\begin{equation}
\begin{aligned}
  \dbE\left[\int_0^\infty \rme^{-rt} \left(X_t^{\xi(\omega),\xi}-X_t^{\xi_n(\omega),\xi_n}\right)^2\rmd t \right]
\leq& \dbE\left[  \dbE\left[\int_0^\infty \rme^{-rt} \left(X_t^{\xi(\omega),\xi}-X_t^{\xi_n(\omega),\xi_n}\right)^2\rmd t\vert \cF_0 \right]   \right]\\
\leq& C\dbE\left[|\xi(\omega)-\xi_n(\omega) |^2+\dbE\left[|\xi-\xi_n|^2\right]   \right]
\rightarrow 0.
\end{aligned}
\end{equation}
So we have $X_t^{\xi}$ and $X_t^{x,\xi}|_{x=\xi}$ are the same in $L_r^2$,
The same holds for $Y_t^{\xi}$ and $Y_t^{x,\xi}|_{x=\xi}$.
Then we conclude that
\begin{equation}
  X_t^{x,\xi}|_{x=\xi}=X_t^{\xi},\quad Y_t^{x,\xi}|_{x=\xi}=Y_t^{\xi},\quad ~\mbox{for $\rmd t\times \rmd \dbP$-a.e. $(t,\omega)$}.
\end{equation}
In addition, by Proposition \ref{prop: con}, we also have
\begin{equation}
  Y_0^{x,\xi_n}|_{x=\xi_n}=Y_0^{\xi_n},\quad ~\mbox{for $\rmd \dbP$-a.e. $\omega$}.
\end{equation}
and
\begin{equation}\begin{aligned}
  &\dbE\left[|Y_0^\xi-Y_0^{\xi_n}|^2\right] \leq C \dbE\left[|\xi-\xi_n|^2\right]\to 0,\\
  &\dbE\left[|Y_0^{\xi(\omega),\xi } -Y_0^{\xi_n(\omega),\xi_n }|^2\right]
  =\dbE\left[\dbE\left[|Y_0^{\xi(\omega),\xi } -Y_0^{\xi_n(\omega),\xi_n }|^2|{\cF_0}\right]\right]\\
 & \leq C\dbE\left[|\xi(\omega)-\xi_n(\omega) |^2+\dbE\left[|\xi-\xi_n|^2\right]   \right]
\rightarrow 0.
  \end{aligned}
\end{equation}
Thus, we have
\begin{equation}
  Y_0^{x,\xi}|_{x=\xi}=Y_0^{\xi},\quad ~\mbox{for $\rmd \dbP$-a.e. $\omega$}.
\end{equation}

\qed

\section{ Solutions to elliptic master equations} 
\label{sec:viscosity}
\setcounter{equation}{0}

This section establishes the viscosity solution of the elliptic
master equation \eqref{eq: master} in the basic weighted \(L_r^2\)
framework.  We first prove the spatial differentiability of the
FBSDE and identify \(\partial_xV\) with its backward initial value.
The distribution-independent equation is already classical at this
level.  The higher-order constructions needed in the measure-dependent
case are developed in Section~\ref{sec: main higher regularity};
classical verification and uniqueness are then proved in
Section~\ref{sec: main classical master}.

\subsection{Differentiability of mean field FBSDEs}
In order to prove the main theorem, we first establish the differentiability of the solutions to 
the mean field FBSDEs. For this purpose, 
we require the following assumptions:
\begin{assum}
  \label{assum: Ha}

\noindent (i) $H(x,\mu,y)$ has at most quadratic growth.  

\noindent (ii) There exist a constant $\lambda_1>r/2$ such that 
  \begin{equation}
    \pa_{yy}H(x,\mu,y)\leq-\lambda_1,\quad \pa_{xx}H(x,\mu,y)\geq\lambda_1.
  \end{equation}

\noindent (iii) There exist a constant $\lambda_2>0$ such that
  \begin{equation}
    |\pa_{xx}H(x,\mu,y)|\leq\lambda_2,\quad
    |\pa_{yy}H(x,\mu,y)|\leq \lambda_2 ,\quad
    |\pa_{xy}H(x,\mu,y)|\leq \lambda_2.
  \end{equation}

\noindent (iv) $\pa_{xx}H(x,\mu,y),\pa_{xy}H(x,\mu,y),\pa_{yy}H(x,\mu,y)$ are Lipschitz continuous.

\noindent (v) $\pa_\mu H(x,\mu,y,\tilde{x})$ satisfies the linear growth condition
and $\pa_{x\mu} H(x,\mu,y,\tilde{x}),\pa_{y\mu} H(x,\mu,y,\tilde{x})$ are bounded.
\end{assum}

We introduce the following FBSDE:
\begin{equation}
  \label{eq: dx}
\begin{cases}
\mathrm{d}  \gxX_t = \left[ \gxX_t \partial_{xy}H(X_{t}^{x,\xi},\mathcal{L} _{X_{t}^{\xi}},Y_{t}^{x, \xi}) + \gxY_t \partial_{yy}H(X_{t}^{x,\xi},\mathcal{L} _{X_{t}^{\xi}},Y_{t}^{x,\xi}) \right] \mathrm{d} t, \\
\begin{aligned}
  \mathrm{d} \gxY_t = & -\left[ \gxX_t \partial_{xx} {H} (X_{t}^{x,\xi},\mathcal{L} _{X_{t}^{\xi}}, Y_{t}^{x,\xi}) + \gxY_t \partial_{xy} {H} (X_{t}^{x,\xi},\mathcal{L} _{X_{t}^{\xi}}, Y_{t}^{x,\xi}) - r\gxY_t \right] \mathrm{d} t \\
  & + \gxZ_t \mathrm{d} B_t,
\end{aligned} \\
\gxX_0 = 1.
\end{cases}
\end{equation}
 Under Assumption \ref{assum: H} and \ref{assum: Ha},
it's clear that  the above FBSDE satisfies all conditions in Assumption \ref{assum: fbsde},
so it admits a unique solution in $L_r^2$.
The following theorem tells us that $(\gxX,\gxY)$ can be 
viewed as the derivative of $(X^{x,\xi},Y^{x,\xi})$ with respect to $x$.

\begin{thm}
  \label{thm: deri}
  For any $x\in\dbR$, we have
  \begin{equation}
    \begin{aligned}
    &\lim_{\delta\to 0}\left\lVert \frac{1}{\delta} \left(X^{x+\delta,\xi}-X^{x,\xi} \right)-\gxX \right\rVert _r= 0,\\
    &\lim_{\delta\to 0}\left\lVert \frac{1}{\delta} \left(Y^{x+\delta,\xi}-Y^{x,\xi} \right)-\gxY \right\rVert _r= 0,\\
    &\lim_{\delta\to 0}\left\lvert \frac{1}{\delta} \left(Y^{x+\delta,\xi}_0-Y^{x,\xi}_0 \right)-\gxY_0          \right\rvert =0.
    \end{aligned}\end{equation}
\end{thm}

\proof
Denote
\begin{equation}
  \Delta X^\delta_t= X^{x+\delta,\xi}_t-X^{x,\xi}_t,\quad
  \Delta Y^\delta_t= Y^{x+\delta,\xi}_t-Y^{x,\xi}_t,\quad
  \Delta Z^\delta_t= Z^{x+\delta,\xi}_t-Z^{x,\xi}_t.
\end{equation}
By Proposition \ref{prop: con}, we have
\begin{equation}
  \left\lVert \left(\Delta X^\delta,\Delta Y^\delta,\Delta Z^\delta \right)\right\rVert _r\leq C\delta,
\end{equation}
where $C$ is a positive constant. 
Therefore we can define the $L^2_r$ processes
\begin{equation}
  \nabla X^\delta \triangleq \frac{\Delta X^\delta}{\delta},\quad 
  \nabla Y^\delta \triangleq \frac{\Delta Y^\delta}{\delta},\quad 
  \nabla Z^\delta \triangleq \frac{\Delta Z^\delta}{\delta},
\end{equation}
and they satisfy the following FBSDE:
\begin{equation}
\begin{cases}
\mathrm{d}  \nabla X^\delta_t = \left[ \nabla X^\delta_t \cdot H_{xy}^\delta + \nabla Y^\delta_t \cdot H_{yy}^\delta \right] \mathrm{d} t, \\
\begin{aligned}
  \mathrm{d} \nabla Y^\delta_t = & -\left[ \nabla X^\delta_t \cdot H_{xx}^\delta + \nabla Y^\delta_t \cdot H_{xy}^\delta - r \nabla Y^\delta_t \right] \mathrm{d} t \\
  & + \nabla Z^\delta_t \mathrm{d} B_t,
\end{aligned} \\
\nabla X^\delta_0 = 1,
\end{cases}
\end{equation}
where
\begin{equation}
  H_{xy}^\delta= \int_0^1
  \pa_{xy}H\left( X_{t}^{x,\xi}+\theta \Delta X^\delta_t,\mathcal{L} _{X_{t}^{\xi}},
  Y_{t}^{x, \xi} +\theta \Delta Y^\delta_t \right)\rmd \theta,
\end{equation}
and $H_{xx}^\delta,H_{yy}^\delta$ are defined similarly.
Define
\begin{equation}
  X^\delta_t\triangleq \nabla X^\delta_t-\nabla X^{x,\xi}_t,\quad
  Y^\delta_t\triangleq \nabla Y^\delta_t-\nabla Y^{x,\xi}_t,\quad
  Z^\delta_t\triangleq \nabla Z^\delta_t-\nabla Z^{x,\xi}_t.
\end{equation}
Then, we have $( X^\delta_t, Y^\delta_t, Z^\delta_t)$ are uniformly bounded in $L_r^2$,
and they satisfy the following FBSDE:
\begin{equation}
\begin{cases}
\mathrm{d} X_t^\delta = \left[ H_{xy}^\delta X_t^\delta + H_{yy}^\delta Y_t^\delta + R_t^x \right] \mathrm{d}t, \\
\begin{aligned}
\mathrm{d} Y_t^\delta &= -\left[ H_{xx}^\delta X_t^\delta + H_{xy}^\delta Y_t^\delta - r Y_t^\delta + R_t^y \right] \mathrm{d}t \\
&\quad + Z_t^\delta \mathrm{d}B_t,
\end{aligned} \\
X_0^\delta = 0,
\end{cases}
\end{equation}
where the remainder terms $R_t^x$ and $R_t^y$ are defined as:
\begin{equation}
\begin{aligned}
R_t^x \triangleq& \nabla X_t^{x,\xi} \left( H_{xy}^\delta - \partial_{xy} H(X_t^{x,\xi}, \mathcal{L}_{X_t^\xi}, Y_t^{x,\xi}) \right) \\
  &+ \nabla Y_t^{x,\xi} \left( H_{yy}^\delta - \partial_{yy} H(X_t^{x,\xi}, \mathcal{L}_{X_t^\xi}, Y_t^{x,\xi}) \right), \\
R_t^y \triangleq& \nabla X_t^{x,\xi} \left(H_{xx}^\delta- \partial_{xx} H(X_t^{x,\xi}, \mathcal{L}_{X_t^\xi}, Y_t^{x,\xi})  \right) 
\\&+ \nabla Y_t^{x,\xi} \left(  H_{xy}^\delta-\partial_{xy} H(X_t^{x,\xi}, \mathcal{L}_{X_t^\xi}, Y_t^{x,\xi})  \right).
\end{aligned}
\end{equation}
Applying It\^{o}'s formula to $\rme^{-rt}X^\delta_t Y^\delta_t$,
and noticing that $H_{yy}^\delta\leq -\lambda_1,  H_{xx}^\delta \geq\lambda_1 ,|H_{xy}^\delta|\leq \lambda_2$,
we get that
\begin{equation}
   \dbE\left[\int_0^\infty \rme^{-rt} \left((X_t^{\delta})^2+(Y_t^{\delta})\right)^2\rmd t \right]
   \leq \frac{2}{r}\dbE\left[\int_0^\infty \rme^{-rt} \left(Y_t^{\delta}\cdot R_t^x+X_t^{\delta}\cdot R_t^y\right)\rmd t \right].
\end{equation}
For simplicity, we only estimate the $Y_t^{\delta}\cdot \nabla X_t^{x,\xi} \left(H_{xx}^\delta- \partial_{xx} H(X_t^{x,\xi}, \mathcal{L}_{X_t^\xi}, Y_t^{x,\xi})  \right)$
term.
Denote 
\begin{equation}
  A_t^\delta\triangleq \left(H_{xx}^\delta- \partial_{xx} H(X_t^{x,\xi}, \mathcal{L}_{X_t^\xi}, Y_t^{x,\xi})  \right).
\end{equation}
From the a priori estimate, we know that $Y_t^{\delta}$ is uniformly bounded in $L_r^2$,
and $|A_t^\delta|\leq 2\lambda_2$, then there exists a constant $C>0$, such that
\begin{equation}\begin{aligned}
  &\dbE\left[\int_0^\infty \rme^{-rt} \left(Y_t^{\delta}\cdot \nabla X_t^{x,\xi}\cdot A_t^\delta\right)\rmd t \right]
\\&\leq  \left( \dbE\left[\int_0^\infty \rme^{-rt} \left(Y_t^{\delta}\right)^2\rmd t \right]
  \dbE\left[\int_0^\infty \rme^{-rt} \left( \nabla X_t^{x,\xi}\cdot A_t^\delta\right)^2\rmd t \right]\right)^{\frac{1}{2}}\\
  &\leq C  \left(\dbE\left[\int_0^\infty \rme^{-rt} \left( \nabla X_t^{x,\xi}\cdot A_t^\delta\right)^2\rmd t \right]\right)^{\frac{1}{2}}.
\end{aligned}\end{equation}
From the Lipschitz continuity of $\pa_{xx}H$, we know that there exists a constant $\ell>0$, such 
that
\begin{equation}
  |A_t^\delta|\leq \ell \left( |\Delta X^\delta_t|+|\Delta Y^\delta_t|\right).
\end{equation}
We consider the following finite measure on $\dbR_+\times \Omega$:
\begin{equation}
  \rmd \dbQ=\rme^{-rt} \rmd t\times \rmd \dbP.
\end{equation}
Then, we know 
\begin{equation}
  \int |A^\delta|^2 \rmd\dbQ \to 0
\end{equation}
as $\delta\to 0$.
This shows that under measure $\dbR$, $|A^\delta|^2$ converges in measure to 0.
Since $(\nabla X^{x,\xi})^2$ is integrable under $\dbQ$,
we know that $( \nabla X^{x,\xi}\cdot A^\delta)^2$ converges in measure to 0.
Finally, we know that $( \nabla X^{x,\xi}\cdot A^\delta)^2 $ is dominated by an integrable function
 $4\lambda_2^2 ( \nabla X^{x,\xi})^2$ under measure Q, 
and by the Dominated Convergence Theorem, we have
\begin{equation}
  \lim_{\delta\to 0}\int \left( \nabla X^{x,\xi}\cdot A^\delta\right)^2 \rmd \dbQ =0.
\end{equation}
This shows that
\begin{equation}
  \lim_{\delta\to 0} \dbE\left[\int_0^\infty \rme^{-rt} \left(Y_t^{\delta}\cdot \nabla X_t^{x,\xi}\cdot A_t^\delta\right)\rmd t \right]=0.
\end{equation}
Then, we have
\begin{equation}
  \lim_{\delta\to 0}\dbE\left[\int_0^\infty \rme^{-rt} \left((X_t^{\delta})^2+(Y_t^{\delta})\right)^2\rmd t \right]=0.
\end{equation}
At last, we apply It\^{o}'s formula to $\rme^{-rt} (Y^\delta_t)^2$.
It's easy to derive that
\begin{equation}
  \lim_{\delta\to 0} \left\lvert Y_0^{\delta}\right\rvert ^2=0.
\end{equation}
Now, we finish the proof.

\qed

Now, we consider a function $F(x,\mu,y):\dbR\times \cP_2\times \dbR\rightarrow \dbR$ satisfying 
the following conditions:
\begin{itemize}
  \item $F$ is of at most quadratic growth and jointly continuous in $(x,\mu,y)$.
  \item $\pa_x F(x,\mu,y)$ and $\pa_y F(x,\mu,y)$ exist and  are Lipschitz continuous in $(x,\mu,y)$.
  \item $\pa_\mu F(x,\mu,y,\tilde{x})$ exist and satisfies the linear growth condition. 
\end{itemize}
 We can define a function $U$ on $\dbR\times \cP_2$ by
\begin{equation}
  U(x,\mu)=\dbE\left[\int_0^\infty \rme^{-rt} F(X_t^{x,\xi},\cL_{X_t^{\xi}},Y_t^{x,\xi})\rmd t \right],
\end{equation}
where $\xi \in\dbL^2(\cF_0;\mu) $.
Since we have already established the uniqueness in distribution for the FBSDEs (see \cite{yang2025discounted}),  $U$ is well-defined.

We now introduce a lemma about $U$ that will be useful for studying the  continuity and differentiability  of the value function (\ref{eq: value}).
\begin{lem}
\label{lem: F}
\noindent (i) $\pa_x U(x,\mu)$ exists and
\begin{equation}
\begin{aligned}
  \pa_xU(x,\mu) &= \dbE\left[\int_0^\infty \rme^{-rt}\left( \pa_x F(X_t^{x,\xi},\cL_{X_t^{\xi}},Y_t^{x,\xi})\gxX_t \right. \right. \\
  &\qquad + \left. \left. \pa_y F(X_t^{x,\xi},\cL_{X_t^{\xi}},Y_t^{x,\xi})\gxY_t\right) \rmd t \right].
\end{aligned}
\end{equation}

\noindent (ii) $U(x,\mu)$ is jointly continuous in $\dbR\times \cP_2$.
\end{lem}
\proof
(i)
We still use the notations in the proof of Theorem \ref{thm: deri}.
Fix $x\in \dbR, \xi\in\dbL^2(\cF_0)$ and assume $|\delta|<1$.
We have $X_t^{x,\xi},{X_t^{\xi}},Y_t^{x,\xi}$ and 
$X_t^{x+\delta,\xi},Y_t^{x+\delta,\xi}$ are all uniformly bounded in $L_r^2$.
From the definition of $U(x,\mu)$, we have
\begin{equation}
\begin{aligned}
  \frac{U(x+\delta,\mu)-U(x,\mu)}{\delta}= \dbE\left[\int_0^\infty \rme^{-rt}\left(  F_x^\delta \cdot \nabla X^\delta_t 
   +   F_y^\delta\cdot \nabla Y^\delta_t\right) \rmd t \right],
\end{aligned}
\end{equation}
where
\begin{equation}
  F_{x}^\delta= \int_0^1
  \pa_{x}F\left( X_{t}^{x,\xi}+\theta \Delta X^\delta_t,\mathcal{L} _{X_{t}^{\xi}},
  Y_{t}^{x, \xi} +\theta \Delta Y^\delta_t \right) \rmd \theta,
\end{equation}
and $F_y^\delta$ is defined similarly. For simplicity, we let
\begin{equation}
\pa_x F=  \pa_x F(X_t^{x,\xi},\cL_{X_t^{\xi}},Y_t^{x,\xi}),\quad
\pa_yF=  \pa_y F(X_t^{x,\xi},\cL_{X_t^{\xi}},Y_t^{x,\xi}).
\end{equation}
Since $\pa_x F(x,\mu,y)$ and $\pa_y F(x,\mu,y)$ are Lipschitz continuous,
we know that
\begin{equation}
   \left\lVert F_x^\delta -F_x \right\rVert _r\to 0,\quad
    \left\lVert F_y^\delta -F_y \right\rVert _r \to 0
\end{equation}
as $\delta\to 0$.
Now, we have
\begin{equation}
\begin{aligned}
& \left\lvert    \frac{U(x+\delta,\mu)-U(x,\mu)}{\delta}- \dbE\left[\int_0^\infty \rme^{-rt}\left( \pa_x F\cdot \gxX_t 
 +  \pa_y F\cdot \gxY_t\right) \rmd t \right]    \right\rvert\\
 \leq &
 \dbE\left[\int_0^\infty \rme^{-rt}\left\lvert F_x^\delta \cdot \nabla X^\delta_t- \pa_x F\cdot \gxX_t 
 \right\rvert \rmd t \right]  
 +\dbE\left[\int_0^\infty \rme^{-rt}\left\lvert F_y^\delta \cdot \nabla Y^\delta_t- \pa_y F\cdot \gxY_t 
 \right\rvert \rmd t \right] \\
 \triangleq&\quad  I+II .
\end{aligned}
\end{equation}
For the first term,
\begin{equation}
\begin{aligned}
I
 \leq &
 \dbE\left[\int_0^\infty \rme^{-rt}\left\lvert F_x^\delta \cdot \left( \nabla X^\delta_t- \gxX_t \right) 
 \right\rvert \rmd t \right]  
 +\dbE\left[\int_0^\infty \rme^{-rt}\left\lvert \left( F_x^\delta -F_x\right) \cdot \gxX_t
 \right\rvert \rmd t \right] \\
 \leq& \left\lVert F_x^\delta\right\rVert _r \cdot  \left\lVert \nabla X^\delta_t- \gxX_t\right\rVert _r
 + \left\lVert F_x^\delta -F_x \right\rVert _r \cdot  \left\lVert \gxX_t\right\rVert _r.
\end{aligned}
\end{equation}
Since $\left\lVert F_x^\delta\right\rVert _r$ and $\left\lVert \gxX_t\right\rVert _r$ are uniformly bounded,
by Theorem \ref{thm: deri}, we know that I tends to 0 as $\delta\rightarrow 0$.
The second term can be treated similarly,
then we finish the proof.

(ii) Fix $x\in \dbR, \mu\in\cP_2$, and consider another $x'\in \dbR, \mu'\in\cP_2$
such that $|x-x'|+\cW_2(\mu,\mu')<\delta$. And then take $\xi,\xi'\in \dbL^2(\cF_0)$,
such that $\cL_\xi=\mu,\cL_{\xi'}=\mu'$.
Without loss of generality, we can assum $|\delta|<1$
and $|x-x'|+(\dbE[\xi-\xi']^2)^{1/2}\leq \delta$
, then all processes appeared in the following
proof are uniformly bounded in $L_r^2$ with respect to $\delta$ .
We have
\begin{equation}
  \begin{aligned}
&F(X_t^{x',\xi'},\cL_{X_t^{\xi'}},Y_t^{x',\xi'}) - F(X_t^{x,\xi},\cL_{X_t^{\xi}},Y_t^{x,\xi})\\
=&F(X_t^{x',\xi'},\cL_{X_t^{\xi'}},Y_t^{x',\xi'}) - F(X_t^{x,\xi},\cL_{X_t^{\xi'}},Y_t^{x',\xi'})\\
&+F(X_t^{x,\xi},\cL_{X_t^{\xi'}},Y_t^{x',\xi'}) - F(X_t^{x,\xi},\cL_{X_t^{\xi}},Y_t^{x',\xi'})\\
&+F(X_t^{x,\xi},\cL_{X_t^{\xi}},Y_t^{x',\xi'}) - F(X_t^{x,\xi},\cL_{X_t^{\xi}},Y_t^{x,\xi})\\ 
\triangleq&\quad I+II+III.   
  \end{aligned}
\end{equation}
For the first term,
\begin{equation}
  \begin{aligned}
    &F(X_t^{x',\xi'},\cL_{X_t^{\xi'}},Y_t^{x',\xi'}) - F(X_t^{x,\xi},\cL_{X_t^{\xi'}},Y_t^{x',\xi'})\\
  =&\int_0^1\pa_x F(X_t^{x,\xi}+\theta (X_t^{x',\xi'}-X_t^{x,\xi}),\cL_{X_t^{\xi'}},Y_t^{x',\xi'})\rmd \theta \cdot (X_t^{x',\xi'}-X_t^{x,\xi})\\
  \triangleq& A_t\cdot  (X_t^{x',\xi'}-X_t^{x,\xi}).
  \end{aligned}
\end{equation}
Since $\pa_xF(x,\mu,y)$ is of at most linear growth,
we know that $A_t$ is uniformly bounded in $L_r^2$.
By Cauchy's Inequality and Proposition \ref{prop: con},  We have
\begin{equation}
  \dbE\left[\rme^{-rt}|I| \rmd t \right]
  \leq \left\lVert A_t \right\rVert _r \cdot \left\lVert X_t^{x',\xi'}-X_t^{x,\xi} \right\rVert _r
  \leq C_1\delta.
\end{equation}
For the second term,
\begin{equation}
\begin{aligned}
  &F(X_t^{x,\xi},\cL_{X_t^{\xi'}},Y_t^{x',\xi'}) - F(X_t^{x,\xi},\cL_{X_t^{\xi}},Y_t^{x',\xi'})\\
=&\int_0^1 \tilde{\dbE}_{\cF_t}
\left[   
  \pa_\mu F\left( X_t^{x,\xi},\cL_{X_t^{\xi}+\theta (X_t^{\xi'}-X_t^{\xi})  },Y_t^{x',\xi'},
  {\tilde{X}_t^{\xi}+\theta (\tilde{X}_t^{\xi'}-\tilde{X}_t^{\xi}) } 
  \right) \cdot \left( \tilde{X}_t^{\xi'}-\tilde{X}_t^{\xi}\right)
\right]\rmd \theta.
\end{aligned}
\end{equation}
Similarly, from the linear growth property of $\pa_\mu F(x,\mu,y,\tilde{x})$, we obtain
\begin{equation}
  \dbE\left[\rme^{-rt}|II| \rmd t \right]
  \leq C_2\delta.
\end{equation}
The approach for III is identical to that for I. Integrating the analyses of I, II, and III, we obtain
there exists a constant $C>0$, such that
\begin{equation}
  |U(x',\mu')-U(x,\mu) |\leq C\delta.
\end{equation}
This completes the proof.
\qed

\subsection{Viscosity solution to elliptic master equation}

\begin{thm}
  \label{thm: deri eq Y}
  Under Assumption \ref{assum: H} and \ref{assum: Ha}, the value function is continuous and satisfies
  \begin{equation}
    \pa_x V(x,\mu)=Y_0^{x,\xi}.
  \end{equation}
\end{thm}
\proof
By the definition of $f(x,\mu,\hat{\alpha}(x,\mu,y))$, we have the relationship:
\begin{equation}
  f(x,\mu,\hat{\alpha}(x,\mu,y))=H(x,\mu,y)-\pa_y H(x,\mu,y)\cdot y.
\end{equation}
It's easy to verify that $H(x,\mu,y)$ satisfies
all conditions in Lemma \ref{lem: F},
whereas $\pa_y H(x,\mu,y)\cdot y$ fails to do so. We therefore require a more delicate analysis.
Fix $x\in \dbR, \xi\in\dbL^2(\cF_0)$,  we continue to use the notation from the proof of Theorem \ref{thm: deri}.
For any $|\delta|<1$,
set
\begin{equation}
\begin{aligned}
U^\delta \triangleq& 
\frac{1}{\delta} \dbE\left[\int_0^\infty \rme^{-rt}
\left( \pa_yH(X_t^{x+\delta,\xi},\cL_{X_t^{\xi}},Y_t^{x+\delta,\xi})\cdot Y^{x+\delta,\xi}_t 
- \pa_yH(X_t^{x,\xi},\cL_{X_t^{\xi}},Y_t^{x,\xi})\cdot Y^{x,\xi}_t
\right)\rmd t \right] \\
=& \dbE\left[\int_0^\infty \rme^{-rt} \nabla Y^{\delta}_t \pa_yH(X_t^{x+\delta,\xi},\cL_{X_t^{\xi}},Y_t^{x+\delta,\xi}) \rmd t \right] \\
& + \dbE\left[\int_0^\infty \rme^{-rt} Y_t^{x,\xi}\left(H_{xy}^\delta \cdot
\nabla X_t^\delta + H_{yy}^\delta\cdot Y_t^\delta \right) \rmd t \right]
\end{aligned}
\end{equation}
and 
\begin{equation}
\begin{aligned}
U \triangleq& \dbE\bigg[ \int_0^\infty \rme^{-rt} 
\nabla Y_t^{x,\xi} \pa_yH(X_t^{x,\xi},\cL_{X_t^{\xi}},Y_t^{x,\xi}) \\
&+ Y_t^{x,\xi} \left(
\pa_{xy}H(X_t^{x,\xi},\cL_{X_t^{\xi}},Y_t^{x,\xi}) \cdot \nabla X_t^{x,\xi}
+ \pa_{yy}H(X_t^{x,\xi},\cL_{X_t^{\xi}},Y_t^{x,\xi}) \cdot \nabla Y_t^{x,\xi}
\right) \rmd t \bigg].
\end{aligned}
\end{equation}
On the one hand, by an application of  H\"{o}lder's inequality, we can obtain that
\begin{equation}
  \begin{aligned}
  &\lim_{\delta\to 0}\dbE\left[\int_0^\infty \rme^{-rt} \nabla Y^{\delta}_t \pa_yH(X_t^{x+\delta,\xi},\cL_{X_t^{\xi}},Y_t^{x+\delta,\xi}) \rmd t \right] 
\\=&\dbE\bigg[ \int_0^\infty \rme^{-rt} 
\nabla Y_t^{x,\xi} \pa_yH(X_t^{x,\xi},\cL_{X_t^{\xi}},Y_t^{x,\xi})\rmd t \bigg].
  \end{aligned}
\end{equation}
On the other hand, employing a method similar to that used in the proof of Theorem \ref{thm: deri}, 
which utilizes the Dominated Convergence Theorem, we can obtain
\begin{equation}
  \begin{aligned}
  &\lim_{\delta\to 0}\dbE\left[\int_0^\infty \rme^{-rt} Y_t^{x,\xi}\left(H_{xy}^\delta \cdot
\nabla X_t^\delta + H_{yy}^\delta\cdot Y_t^\delta \right) \rmd t \right] 
\\=&\dbE\bigg[ \int_0^\infty \rme^{-rt} 
Y_t^{x,\xi} \left(
\pa_{xy}H(X_t^{x,\xi},\cL_{X_t^{\xi}},Y_t^{x,\xi}) \cdot \nabla X_t^{x,\xi}
+ \pa_{yy}H(X_t^{x,\xi},\cL_{X_t^{\xi}},Y_t^{x,\xi}) \cdot \nabla Y_t^{x,\xi}
\right) \rmd t \bigg].
  \end{aligned}
\end{equation}
Therefore
\begin{equation}
  \lim_{\delta\to 0}U^\delta=U.
\end{equation}

Noting that all the first-order partial derivatives of $f(x,\mu,\hat{\alpha}(x,\mu,y))$ have at most linear growth, 
we obtain from Lemma \ref{lem: F} that the value function $V(x,\mu)$ is continuous and
 $\pa_x V(x,\mu)$ exists, which satisfies
\begin{equation}
\begin{aligned}
  \partial_x V(x,\mu) =\dbE \int_0^\infty \mathrm{e}^{-rt} &\left\{
    \left[
      \partial_{x} H(X_t^{x,\xi}, \mathcal{L}_{X_t^\xi}, Y_t^{x,\xi}) - \partial_{xy} H(X_t^{x,\xi}, \mathcal{L}_{X_t^\xi}, Y_t^{x,\xi}) \cdot Y_t^{x,\xi}
    \right] \cdot \gxX_t \right. \\
  &\left. - \partial_{yy} H(X_t^{x,\xi}, \mathcal{L}_{X_t^\xi}, Y_t^{x,\xi}) \cdot Y_t^{x,\xi} \cdot \gxY_t
  \right\} \mathrm{d} t.
\end{aligned}
\end{equation}
Applying It\^{o}'s formula to $\rme^{-rt}\gxX_t\cdot Y_t^{x,\xi}$, we have
\begin{equation}
\begin{aligned}
\mathrm{d}\left( \mathrm{e}^{-rt} \gxX_t Y_t^{x,\xi} \right) &= \mathrm{e}^{-rt} \left[ -\gxX_t \partial_x H(X_{t}^{x,\xi},\mathcal{L} _{X_{t}^{\xi}}, Y_{t}^{x,\xi}) + Y_t^{x,\xi} \gxX_t \partial_{xy}H(X_{t}^{x,\xi},\mathcal{L} _{X_{t}^{\xi}}, Y_t^{x,\xi}) \right. \\
&\quad \left. + Y_t^{x,\xi} \gxY_t \partial_{yy}H(X_{t}^{x,\xi},\mathcal{L} _{X_{t}^{\xi}}, Y_t^{x,\xi}) \right] \mathrm{d}t \\
&\quad + \mathrm{e}^{-rt} \gxX_t Z_t^{x} \mathrm{d}B_t.
\end{aligned}
\end{equation}
 Note that $\gxX_0 = 1$, and take a sequence $T_i \to \infty$ such that 
 \begin{equation}
  \mathbb{E}\left[\mathrm{e}^{-rT_i} \gxX_{T_i}\cdot Y_{T_i}^{x,\xi}\right] \to 0.
 \end{equation}
 Integrating  from $0$ to $T_i$ and taking expectation, after letting $T_i\to \infty$ we have:
\begin{equation}
\begin{aligned}
  Y_0^{x,\xi} =\dbE \int_0^\infty \mathrm{e}^{-rt} &\left\{
    \left[
      \partial_{x} H(X_t^{x,\xi}, \mathcal{L}_{X_t^\xi}, Y_t^{x,\xi}) - \partial_{xy} H(X_t^{x,\xi}, \mathcal{L}_{X_t^\xi}, Y_t^{x,\xi}) \cdot Y_t^{x,\xi}
    \right] \cdot \gxX_t \right. \\
  &\left. - \partial_{yy} H(X_t^{x,\xi}, \mathcal{L}_{X_t^\xi}, Y_t^{x,\xi}) \cdot Y_t^{x,\xi} \cdot \gxY_t
  \right\} \mathrm{d} t.
\end{aligned}
\end{equation}
Now, we get the desired result.

\qed

In preparation for the definition of a viscosity solution
of the elliptic master equation, we first define the class of test functions used for that purpose.
\begin{defn}
\label{def: test function}
A function \(\Psi\in C^{2,1}(\dbR\times\cP_2)\) is said to be
a test function if, for every compact subset
\(\cK\subset\dbR\times\cP_2\),
\begin{equation}
\label{eq: main test integrability}
 \sup_{(x,\mu)\in\cK}\int_{\dbR}
 \left(|\pa_\mu\Psi(x,\mu,\tilde x)|^2
       +|\pa_{\tilde x}\pa_\mu\Psi(x,\mu,\tilde x)|^2\right)
 \rmd\mu(\tilde x)<\infty.
\end{equation}
\end{defn}

Now, we present the definition of the viscosity solution for elliptic master equation (\ref{eq: master}).
\begin{defn}
  Suppose that $U\in C(\mathbb{R}\times \mathcal{P} _2 )$ and its partial derivative
   $\pa_x U\in C(\mathbb{R}\times \mathcal{P} _2 )$. Then $U$ is called a viscosity subsolution (resp. supersolution)
  of PDE (\ref{eq: master}) if, for any $(x^0,\mu^0)\in \mathbb{R} \times \mathcal{P} _2$
  and any test function $\Psi$, such that
  $(x^0,\mu^0)$ is a local maximum (resp. minimum) of $U-\Psi$, we have
\begin{equation}
    \begin{split}
  r U(x^0,\mu^0)\leq &H(x^0,\mu^0,\partial_x U(x^0,\mu^0))+\frac{1}{2}\partial_{xx}\Psi(x^0,\mu^0)\\&+\tilde{\mathbb{E} }
  \left[\frac{1}{2}\partial_{\tilde{x}}\partial_\mu \Psi(x^0,\mu^0,\tilde{\xi }^0)+\partial_\mu \Psi(x^0,\mu^0,\tilde{\xi }^0)
 \cdot \partial_y H(\tilde{\xi }^0,\mu^0,\partial_xU(\tilde{\xi}^0,\mu^0))\right].
\end{split}
\end{equation}
  (respectively,
\begin{equation}
    \begin{split}
  r U(x^0,\mu^0)\geq &H(x^0,\mu^0,\partial_xU(x^0,\mu^0))+\frac{1}{2}\partial_{xx}\Psi(x^0,\mu^0)\\&+\tilde{\mathbb{E} }
  \left[\frac{1}{2}\partial_{\tilde{x}}\partial_\mu \Psi(x^0,\mu^0,\tilde{\xi }^0)+\partial_\mu \Psi(x^0,\mu^0,\tilde{\xi }^0)
 \cdot \partial_y H(\tilde{\xi }^0,\mu^0,\partial_xU(\tilde{\xi}^0,\mu^0))\right],
\end{split}
\end{equation}
  ).

  The function $U$ is called a viscosity solution of PDE (\ref{eq: master}) if it is 
  both a viscosity subsolution and a viscosity supersolution.
\end{defn}

\begin{thm}
  \label{thm: value viscosity}
  Under Assumption \ref{assum: H} and \ref{assum: Ha},
  the value function (\ref{eq: value})
is the viscosity solution to the elliptic master equation (\ref{eq: master}).
\end{thm}

\proof
We have proved that $V(x,\mu)$ is a continuous function
and $\pa_x V(x,\mu)=Y_0^{x,\xi}$. Since $\pa_x V(x,\mu) $ is also a continuous function,
by the flow property, we have 
\begin{equation}
  \label{eq: flow}
  \pa_x V(X_t^{x,\xi},\cL_{X_t^{\xi}})=Y_t^{x,\xi},\quad \pa_x V(X_t^{\xi},\cL_{X_t^{\xi}})=Y_t^{\xi}.
\end{equation}
We only show that $V$ is a viscosity subsolution of PDE (\ref{eq: master}). 
A similar argument will show that it is also a viscosity supersolution of (\ref{eq: master}).

Let $\Psi\in C^{2,1}(\mathbb{R} \times \mathcal{P} _2)$ be a test function
and  $(x^0,\mu^0)\in \mathbb{R} \times \mathcal{P} _2$ be a local maximum  of $V-\Psi$. 
It's natural to get that $\pa_xV(x^0,\mu^0)=\pa_x\Psi(x^0,\mu^0)$.
We assume without loss of generality that $V(x^0,\mu^0)=\Psi(x^0,\mu^0)$.
And we suppose that
\begin{equation}
      \begin{split}
  r \Psi(x^0,\mu^0)> &H(x^0,\mu^0,\partial_xV(x^0,\mu^0))+\frac{1}{2}\partial_{xx}\Psi(x^0,\mu^0)\\
  &+\tilde{\mathbb{E} }
  \left[\frac{1}{2}\partial_{\tilde{x}}\partial_\mu \Psi(x^0,\mu^0,\tilde{\xi }^0)+\partial_\mu \Psi(x^0,\mu^0,\tilde{\xi }^0)
  \partial_y H(\tilde{\xi },\mu^0,\partial_xV(\tilde{\xi}^0,\mu^0))\right].
\end{split}
\end{equation}
Notice that
\begin{equation}
  \begin{split}
   &H(x^0,\mu^0,\pa_x V(x^0,\mu^0))\\=&
   f(x^0,\mu^0,\hat{\alpha}(x^0,\mu^0,\pa_x V(x^0,\mu^0)))+\pa_x V(x^0,\mu^0)\cdot \pa_yH(x^0,\mu^0,\pa_x V(x^0,\mu^0))\\
 =&f(x^0,\mu^0,\hat{\alpha}(x^0,\mu^0,\pa_x V(x^0,\mu^0)))+\pa_x \Psi(x^0,\mu^0)\cdot \pa_yH(x^0,\mu^0,\pa_x V(x^0,\mu^0)).
  \end{split}
\end{equation}
We can get
\begin{equation}
      \begin{split}
  r \Psi(x^0,\mu^0)> &f(x^0,\mu^0,\hat{\alpha}(x^0,\mu^0,\pa_x V(x^0,\mu^0)))+\pa_x \Psi(x^0,\mu^0)\cdot \pa_yH(x^0,\mu^0,\pa_x V(x^0,\mu^0))
  \\&+\frac{1}{2}\partial_{xx}\Psi(x^0,\mu^0)\\
  &+\tilde{\mathbb{E} }
  \left[\frac{1}{2}\partial_{\tilde{x}}\partial_\mu \Psi(x^0,\mu^0,\tilde{\xi }^0)+\partial_\mu \Psi(x^0,\mu^0,\tilde{\xi }^0)
  \partial_y H(\tilde{\xi },\mu^0,\partial_xV(\tilde{\xi}^0,\mu^0))\right].
\end{split}
\end{equation}

The measure integrals are continuous. Indeed, under an \(L^2\)
coupling, the derivatives of the test function converge in probability
and are bounded in \(L^2\), while the drift converges in \(L^2\).
Truncation and the Cauchy--Schwarz inequality give convergence of
the products in \(L^1\).
It follows from the above that there exists an open subset $O\subset \mathbb{R} \times \mathcal{P} _2$ that contains
$(x^0,\mu^0)$, such that for all $(x,\mu)\in O$,
\begin{equation}
  \label{eq: contradic}
  \begin{cases}
    V(x,\mu)\le \Psi(x,\mu),\\
    \begin{aligned}
       r \Psi(x,\mu)>&f(x,\mu,\hat{\alpha}(x,\mu,\pa_x V(x,\mu)))+\pa_x \Psi(x,\mu)\cdot \pa_yH(x,\mu,\pa_x V(x,\mu))
       \\&+\frac{1}{2}\partial_{xx}\Psi(x,\mu)\\&+\tilde{\mathbb{E} }
  \left[\frac{1}{2}\partial_{\tilde{x}}\partial_\mu \Psi(x,\mu,\tilde{\xi })+\partial_\mu \Psi(x,\mu,\tilde{\xi })
  \partial_y H(\tilde{\xi },\mu,\partial_xV(\tilde{\xi},\mu))\right].
    \end{aligned}
  \end{cases}
\end{equation}
Taking an initial state $\xi^0\in\mathbb{L} ^2(\mathcal{F} _0)$ such that $\cL_{\xi^0}=\mu^0$, we consider the processes 
$(X_t^{\xi^0},Y_t^{\xi^0},Z_t^{\xi^0})$ and $(X_t^{x^0,\xi^0},Y_t^{x^0,\xi^0},Z_t^{x^0,\xi^0})$ which are solutions to FBSDEs (\ref{eq: mf1}) and (\ref{eq: mf2}).
We denote $\rho_t\triangleq \mathcal{L} _{X_t^{\xi^0}}$.
For some $T>0$, let $\tau$ denote the stopping time
\begin{equation}
  \tau\triangleq \inf\{t>0|(X_t^{x^0,\xi^0}, \rho_t)\notin O\}\land T.
\end{equation}
By the flow property, we have that
\begin{equation}
  \begin{aligned}
    \Psi(x^0,\mu^0) &= V(x^0,\mu^0) \\
    &= \mathbb{E}\bigg[ \int_{0}^{\tau} \mathrm{e}^{-rt} f\big( X_{t}^{x^0,\xi^0}, \rho_t, \hat{\alpha}(X_{t}^{x^0,\xi^0},\rho_t, Y_{t}^{x^0,\xi^0}) \big)\mathrm{d}t \\
    &\quad + \mathrm{e}^{-r\tau} \int_{\tau}^\infty \mathrm{e}^{-r(t-\tau)} f\big( X_{t}^{x^0,\xi^0}, \rho_t, \hat{\alpha}(X_{t}^{x^0,\xi^0},\rho_t, Y_{t}^{x^0,\xi^0}) \big)\rmd t  \bigg] \\
    &= \mathbb{E}\bigg[ \int_{0}^{\tau} \mathrm{e}^{-rt} f\big( X_{t}^{x^0,\xi^0}, \rho_t, \hat{\alpha}(X_{t}^{x^0,\xi^0},\rho_t, Y_{t}^{x^0,\xi^0}) \big) \mathrm{d}t
     + \mathrm{e}^{-r\tau} V(X_{\tau}^{x^0,\xi^0}, \rho_\tau) \bigg] \\
    &\leq \mathbb{E}\bigg[ \int_{0}^{\tau} \mathrm{e}^{-rt} f\big( X_{t}^{x^0,\xi^0}, \rho_t, \hat{\alpha}(X_{t}^{x^0,\xi^0},\rho_t, Y_{t}^{x^0,\xi^0}) \big) \mathrm{d}t
     + \mathrm{e}^{-r\tau} \Psi(X_{\tau}^{x^0,\xi^0}, \rho_\tau ) \bigg].
  \end{aligned}
\end{equation}
By the definition of the test function \(\Psi\), we can apply
It\^{o}'s formula \eqref{eq:ito} to
\(\rme^{-rt}\Psi(X_t^{x^0,\xi^0},\rho_t)\).
The linearly growing drift is treated by truncation; the finite-time
second-moment estimates and \eqref{eq: main test integrability}
justify passage to the limit. Then we get that
\begin{equation}
\begin{aligned}
0 \leq &\ \mathbb{E}\bigg[ \int_{0}^{\tau} \mathrm{e}^{-rt} \bigg( 
    f\Big( X_{t}^{x^0,\xi^0}, \rho_t, \hat{\alpha}\big(X_{t}^{x^0,\xi^0},\rho_t, Y_{t}^{x^0,\xi^0}\big) \Big) \\
&\quad - r\Psi(X_{t}^{x^0,\xi^0}, \rho_t) 
    + \partial_x\Psi(X_{t}^{x^0,\xi^0}, \rho_t) \cdot \partial_y H(X_{t}^{x^0,\xi^0}, \rho_t, Y_{t}^{x^0,\xi^0}) \\
&\quad + \frac{1}{2}\partial_{xx}\Psi(X_{t}^{x^0,\xi^0}, \rho_t) \\
&\quad + \tilde{\mathbb{E}}_{\mathcal{F}_t} \Big[ 
        \frac{1}{2}\partial_{\tilde{x}}\partial_\mu \Psi(X_{t}^{x^0,\xi^0}, \rho_t,\tilde{X}_{t}^{\xi^0}) \\
&\qquad + \partial_\mu \Psi(X_{t}^{x^0,\xi^0}, \rho_t,\tilde{X}_{t}^{\xi^0}) 
        \cdot \partial_y H(\tilde{X}_{t}^{\xi^0},\rho_t,\tilde{Y}_{t}^{\xi^0}) 
    \Big] 
\bigg) \mathrm{d}t \bigg].
\end{aligned}
\end{equation}
Then, by the relationship (\ref{eq: flow}), we have that
\begin{equation}
  \begin{aligned}
    0\leq
&\ \mathbb{E}\bigg[ \int_{0}^{\tau} \mathrm{e}^{-rt} \bigg( 
    -r\Psi(X_{t}^{x^0,\xi^0}, \rho_t) \\
&\quad + f\Big(X_{t}^{x^0,\xi^0}, \rho_t, \hat{\alpha}\big(X_{t}^{x^0,\xi^0},\rho_t, \partial_x V(X_{t}^{x^0,\xi^0}, \rho_t)\big)\Big) \\
&\quad + \partial_x \Psi(X_{t}^{x^0,\xi^0}, \rho_t) \cdot \partial_y H\Big(X_{t}^{x^0,\xi^0}, \rho_t, \partial_x V(X_{t}^{x^0,\xi^0}, \rho_t)\Big) \\
&\quad + \frac{1}{2}\partial_{xx}\Psi(X_{t}^{x^0,\xi^0}, \rho_t) \\
&\quad + \tilde{\mathbb{E}}_{\mathcal{F}_t} \bigg[ 
        \frac{1}{2}\partial_{\tilde{x}}\partial_\mu \Psi(X_{t}^{x^0,\xi^0}, \rho_t, \tilde{X}_{t}^{\xi^0}) \\
&\qquad + \partial_\mu \Psi(X_{t}^{x^0,\xi^0}, \rho_t, \tilde{X}_{t}^{\xi^0}) 
        \cdot \partial_y H\Big(\tilde{X}_{t}^{\xi^0}, \rho_t, \partial_x V(\tilde{X}_{t}^{\xi^0}, \rho_t)\Big)
    \bigg] 
\bigg) \mathrm{d}t \bigg],
  \end{aligned}
\end{equation}
which contradicts  (\ref{eq: contradic}).
Now, we finish the proof.
\qed

In Theorem \ref{thm: deri} and \ref{thm: deri eq Y}, we prove the existence of both $\pa_x V(x,\mu)$
and $\pa_{xx}V(x,\mu)$ 
for the value function $V(x,\mu)$. 
This indicates that the proof of the viscosity solution deals exclusively 
with the distributional component.
For the distribution-independent case, in which  $b(x,\mu,a)$ and $f(x,\mu,a)$
reduce to $b(x,a)$ and $f(x,a)$, and consequently the Hamiltonian $H(x,\mu,y)$ becomes $H(x,y)$, 
we can obtain a stronger conclusion.
In this case, we write \(\hat{\alpha}(x,y)\) for the minimizing feedback.
Considering the unique $L_r^2$-solution to the reduced mean field FBSDE:
\begin{equation}
  \begin{cases}
    \rmd X_{t}^{x} = \partial_yH(X_{t}^{x},Y_{t}^{x}) \rmd t +\rmd B_{t}, \\
    \rmd Y_{t}^{x} = -\left[\partial_x {H} (X_{t}^{x}, Y_{t}^{x})-rY_{t}^{x}\right] \rmd t + Z_{t}^{x}\rmd B_{t}, \\
    X_0^x=x,
    \end{cases}
\end{equation}
we define the value function
\begin{equation}
\label{eq: value-reduce}
V(x)=  \mathbb{E}\bigg[\int_{0}^{+\infty}\rme^{-rt}f\big(X_{t}^{x},
\hat{\alpha}(X_{t}^{x}, Y_{t}^{x})        \big)\rmd t\bigg],
\end{equation}
where $f(x,\hat{\alpha}(x,y))$ can be identified with $H(x,y)-\pa_yH(x,y)\cdot y$.
From the preceding derivation, we have the following proposition.
\begin{prop}
  \label{prop: classical}
Suppose that the reduced Hamiltonian $H(x,y)$ satisfies Assumption \ref{assum: H} and \ref{assum: Ha}.
We have that $V(x)$ defined in (\ref{eq: value-reduce}) is a classical solution to the equation
\begin{equation}
  rU(x)=H(x,\pa_xU(x))+\frac{1}{2}\pa_{xx}U(x).
\end{equation}
\end{prop}

\section{Higher-order regularity of the value function}
\label{sec: main higher regularity}
\label{subsec: value classical}
\setcounter{equation}{0}

This section establishes the regularity of the value function
\(V(x,\mu)\) needed for the elliptic master equation. We first
recall from \cite{songyangregularity} the variational equations and
convergence results for \(\pa_xV(x,\mu)\). The construction of
\(\pa_\mu\pa_xV(x,\mu,\tilde x)\) is carried out in \(L_r^2\),
by passing from discrete to general initial laws. We then introduce
the \(L^p\) estimates used to construct the higher-order variations.
Using these results and the discounted representation of
\(V(x,\mu)\), we construct \(\pa_\mu V(x,\mu,\tilde x)\) from
directional derivatives, prove its continuity and differentiability
in \(x\) and \(\tilde x\), and establish the growth bounds for the
derivatives appearing in the master equation. These bounds will
allow us to verify the classical solution property in
Section~\ref{sec: main classical master}.

\subsection{Estimates and variational equations}
\label{subsec: main regularity assumptions}
\label{subsec: main weighted estimates}

In this subsection, we recall the derivative constructions, estimates,
and convergence results of \cite[Sections~3--5]{songyangregularity}
that will be used to prove the regularity of \(V(x,\mu)\), and describe
the Riccati construction underlying the higher-moment estimates.
Using the notation of that paper, we write
\[
 \mcv(x,\mu)=Y_0^{x,\xi},\qquad \xi\in\dbL^2(\cF_0;\mu).
\]
Theorem~\ref{thm: deri eq Y} gives
\(\pa_xV(x,\mu)=Y_0^{x,\xi}\), so that
\(\mcv(x,\mu)=\pa_xV(x,\mu)\).

\subsubsection{\texorpdfstring{Measure derivatives in \(L_r^2\)}{Measure derivatives in L2}}
\label{subsec: main L2 construction}

\begin{assum}
\label{assum: value classical}
Assumptions~\ref{assum: H} and~\ref{assum: Ha} hold. In addition:
\begin{enumerate}[label=\textnormal{(\roman*)}]
\item The derivatives \(\pa_\mu H,\pa_{x\mu}H,
\pa_{y\mu}H\) are globally Lipschitz in all variables.
There are constants \(b_0,c_0>r/2\) and \(\lambda_\mu\geq0\)
such that
\[
 -\pa_{yy}H\geq b_0,\qquad \pa_{xx}H\geq c_0,\qquad
 |\pa_{x\mu}H|\vee|\pa_{y\mu}H|\leq\lambda_\mu,
\]
and \(\min\{b_0,c_0\}-2\lambda_\mu>r/2\).
\item All derivatives \(\pa_x^i\pa_y^jH\) with \(i+j=3\)
and the derivatives
\(\pa_{\tilde x}\pa_{q\mu}H(x,\mu,y,\tilde x)\),
\(q=x,y\), exist, are bounded, and are globally Lipschitz.
\item The derivative
\(\pa_{\tilde x}\pa_\mu H(x,\mu,y,\tilde x)\)
exists, is bounded, and is globally Lipschitz.
\end{enumerate}
Here Lipschitz continuity in the measure variable is understood
with respect to \(\cW_2\).
\end{assum}

We first recall the construction of
\(\pa_\mu\mcv(x,\mu,\tilde x)\). This part uses only Assumption~\ref{assum: value classical}(i)
and takes place in \(L_r^2\). For
\(\xi,\eta\in\dbL^2(\cF_0)\) and \(x\in\dbR\), set
\(\mu=\cL_\xi\) and \(\mu_t=\cL_{X_t^\xi}\).
Tildes denote independent copies, and
\(\widetilde{\dbE}_{\cF_t}\) denotes expectation over the copied variables.
We introduce the following FBSDEs:
\begin{equation}
\label{eq: main population directional FBSDE}
\begin{cases}
\begin{aligned}
 \rmd\delta X_t^{\xi,\eta}={}&\bigl\{ \delta X_t^{\xi,\eta}\pa_{xy}H(X_t^\xi,\mu_t,Y_t^\xi) +\delta Y_t^{\xi,\eta}\pa_{yy}H(X_t^\xi,\mu_t,Y_t^\xi)\\
 &+\widetilde{\dbE}_{\cF_t}\bigl[
 \pa_{y\mu}H(X_t^\xi,\mu_t,Y_t^\xi,\widetilde X_t^\xi)
                 \delta\widetilde X_t^{\xi,\eta}\bigr]
 \bigr\}\,\rmd t,
\end{aligned}\\
\begin{aligned}
 \rmd\delta Y_t^{\xi,\eta}=-{}&\bigl\{ \delta X_t^{\xi,\eta}\pa_{xx}H(X_t^\xi,\mu_t,Y_t^\xi) +\delta Y_t^{\xi,\eta}\pa_{xy}H(X_t^\xi,\mu_t,Y_t^\xi) -r\delta Y_t^{\xi,\eta}\\
 &+\widetilde{\dbE}_{\cF_t}\bigl[
 \pa_{x\mu}H(X_t^\xi,\mu_t,Y_t^\xi,\widetilde X_t^\xi)
                 \delta\widetilde X_t^{\xi,\eta}\bigr]
 \bigr\}\,\rmd t+\delta Z_t^{\xi,\eta}\,\rmd B_t,
\end{aligned}\\
\delta X_0^{\xi,\eta}=\eta.
\end{cases}
\end{equation}
\begin{equation}
\label{eq: main tagged directional FBSDE}
\begin{cases}
\begin{aligned}
 \rmd\delta X_t^{x,\xi,\eta}={}&\bigl\{ \delta X_t^{x,\xi,\eta}\pa_{xy}H(X_t^{x,\xi},\mu_t,Y_t^{x,\xi}) +\delta Y_t^{x,\xi,\eta}\pa_{yy}H(X_t^{x,\xi},\mu_t,Y_t^{x,\xi})\\
 &+\widetilde{\dbE}_{\cF_t}\bigl[
 \pa_{y\mu}H(X_t^{x,\xi},\mu_t,Y_t^{x,\xi},\widetilde X_t^\xi)
                 \delta\widetilde X_t^{\xi,\eta}\bigr]
 \bigr\}\,\rmd t,
\end{aligned}\\
\begin{aligned}
 \rmd\delta Y_t^{x,\xi,\eta}=-{}&\bigl\{ \delta X_t^{x,\xi,\eta}\pa_{xx}H(X_t^{x,\xi},\mu_t,Y_t^{x,\xi}) +\delta Y_t^{x,\xi,\eta}\pa_{xy}H(X_t^{x,\xi},\mu_t,Y_t^{x,\xi}) -r\delta Y_t^{x,\xi,\eta}\\
 &+\widetilde{\dbE}_{\cF_t}\bigl[
 \pa_{x\mu}H(X_t^{x,\xi},\mu_t,Y_t^{x,\xi},\widetilde X_t^\xi)
                 \delta\widetilde X_t^{\xi,\eta}\bigr]
 \bigr\}\,\rmd t+\delta Z_t^{x,\xi,\eta}\,\rmd B_t,
\end{aligned}\\
\delta X_0^{x,\xi,\eta}=0.
\end{cases}
\end{equation}
Under Assumption~\ref{assum: value classical}(i), their unique
solutions in \(L_r^2\) are denoted by
\((\delta X^{\xi,\eta},\delta Y^{\xi,\eta},\delta Z^{\xi,\eta})\)
and
\((\delta X^{x,\xi,\eta},\delta Y^{x,\xi,\eta},\delta Z^{x,\xi,\eta})\),
respectively. For \(\Phi\in\{X,Y,Z\}\), we write
\(\delta\Phi^{\xi,\eta}\) and \(\delta\Phi^{x,\xi,\eta}\)
for their corresponding components. They satisfy
\begin{equation}
\label{eq: main directional process limits}
\begin{aligned}
 &\lim_{\delta\to0}
 \left\|\frac{\Phi^{\xi+\delta\eta}-\Phi^\xi}{\delta}
               -\delta\Phi^{\xi,\eta}\right\|_{L_r^2}=0,\\
 &\lim_{\delta\to0}
 \left\|\frac{\Phi^{x,\xi+\delta\eta}-\Phi^{x,\xi}}{\delta}
               -\delta\Phi^{x,\xi,\eta}\right\|_{L_r^2}=0.
\end{aligned}
\end{equation}
As \(\delta\to0\), the backward initial values satisfy
\[
\begin{aligned}
 &\dbE\left|
   \frac{Y_0^{\xi+\delta\eta}-Y_0^\xi}{\delta}
       -\delta Y_0^{\xi,\eta}\right|^2\longrightarrow0,\\
 &\left|
   \frac{\mcv(x,\cL_{\xi+\delta\eta})-\mcv(x,\cL_\xi)}{\delta}
       -\delta Y_0^{x,\xi,\eta}\right|\longrightarrow0.
\end{aligned}
\]
These variations are bounded linear functions of \(\eta\), and are
jointly continuous in \((x,\xi,\eta)\), in \(L_r^2\) for the
processes and in \(\dbL^2\) for the backward initial values.

Following \cite{gangbo2022mean,mou2024wellposedness}, we construct
the Lions-derivative of \(\mcv\). First consider a discrete random
variable \(\xi\), with
\begin{equation}
\label{eq: main discrete law}
 \dbP(\xi=x_i)=p_i>0,\qquad i=1,\ldots,n,
 \qquad\sum_{i=1}^np_i=1.
\end{equation}
For each \(i\), introduce the following FBSDEs. The first is a
coupled system whose solution enters the second:
\begin{equation}
\label{eq: main discrete part}
\begin{cases}
\begin{aligned}
 \rmd \nabla X^{\xi,x_i}_t={}&\bigl[ \nabla X^{\xi,x_i}_t\pa_{xy}H(X_t^{x_i,\xi},\mu_t,Y_t^{x_i,\xi}) +\nabla Y^{\xi,x_i}_t\pa_{yy}H(X_t^{x_i,\xi},\mu_t,Y_t^{x_i,\xi})\\
 &+p_i\widetilde{\dbE}_{\cF_t}\bigl[
 \nabla\widetilde X_t^{\xi,x_i}\pa_{y\mu}H(X_t^{x_i,\xi},\mu_t,Y_t^{x_i,\xi},\widetilde X_t^{x_i,\xi})\\
 &\quad +\nabla\widetilde X_t^{\xi,x_i,\star}\pa_{y\mu}H(X_t^{x_i,\xi},\mu_t,Y_t^{x_i,\xi},\widetilde X_t^{\xi})\bigr]\bigr]\,\rmd t,
\end{aligned}\\
\begin{aligned}
 \rmd \nabla X^{\xi,x_i,\star}_t={}&\bigl[ \nabla X^{\xi,x_i,\star}_t\pa_{xy}H(X_t^{\xi},\mu_t,Y_t^{\xi}) +\nabla Y^{\xi,x_i,\star}_t\pa_{yy}H(X_t^{\xi},\mu_t,Y_t^{\xi})\\
 &+\widetilde{\dbE}_{\cF_t}\bigl[
 \nabla\widetilde X_t^{\xi,x_i}\pa_{y\mu}H(X_t^{\xi},\mu_t,Y_t^{\xi},\widetilde X_t^{x_i,\xi})\\
 &\quad +\nabla\widetilde X_t^{\xi,x_i,\star}\pa_{y\mu}H(X_t^{\xi},\mu_t,Y_t^{\xi},\widetilde X_t^{\xi})\bigr]I_{\{\xi\ne x_i\}}\bigr]\,\rmd t,
\end{aligned}\\
\begin{aligned}
 \rmd \nabla Y^{\xi,x_i}_t=-{}&\bigl[ \nabla X^{\xi,x_i}_t\pa_{xx}H(X_t^{x_i,\xi},\mu_t,Y_t^{x_i,\xi}) +\nabla Y^{\xi,x_i}_t\pa_{xy}H(X_t^{x_i,\xi},\mu_t,Y_t^{x_i,\xi})-r\nabla Y^{\xi,x_i}_t\\
 &+p_i\widetilde{\dbE}_{\cF_t}\bigl[
 \nabla\widetilde X_t^{\xi,x_i}\pa_{x\mu}H(X_t^{x_i,\xi},\mu_t,Y_t^{x_i,\xi},\widetilde X_t^{x_i,\xi})\\
 &\quad +\nabla\widetilde X_t^{\xi,x_i,\star}\pa_{x\mu}H(X_t^{x_i,\xi},\mu_t,Y_t^{x_i,\xi},\widetilde X_t^{\xi})\bigr]\bigr]\,\rmd t+\nabla Z^{\xi,x_i}_t\,\rmd B_t,
\end{aligned}\\
\begin{aligned}
 \rmd \nabla Y^{\xi,x_i,\star}_t=-{}&\bigl[ \nabla X^{\xi,x_i,\star}_t\pa_{xx}H(X_t^{\xi},\mu_t,Y_t^{\xi}) +\nabla Y^{\xi,x_i,\star}_t\pa_{xy}H(X_t^{\xi},\mu_t,Y_t^{\xi})-r\nabla Y^{\xi,x_i,\star}_t\\
 &+\widetilde{\dbE}_{\cF_t}\bigl[
 \nabla\widetilde X_t^{\xi,x_i}\pa_{x\mu}H(X_t^{\xi},\mu_t,Y_t^{\xi},\widetilde X_t^{x_i,\xi})\\
 &\quad +\nabla\widetilde X_t^{\xi,x_i,\star}\pa_{x\mu}H(X_t^{\xi},\mu_t,Y_t^{\xi},\widetilde X_t^{\xi})\bigr]I_{\{\xi\ne x_i\}}\bigr]\,\rmd t+\nabla Z^{\xi,x_i,\star}_t\,\rmd B_t,
\end{aligned}\\
\nabla X^{\xi,x_i}_0=1,\quad \nabla X^{\xi,x_i,\star}_0=0.
\end{cases}
\end{equation}

\begin{equation}
\label{eq: main discrete total}
\begin{cases}
\begin{aligned}
 \rmd \nabla_\mu X^{x,\xi,x_i}_t={}&\bigl[ \nabla_\mu X^{x,\xi,x_i}_t\pa_{xy}H(X_t^{x,\xi},\mu_t,Y_t^{x,\xi}) +\nabla_\mu Y^{x,\xi,x_i}_t\pa_{yy}H(X_t^{x,\xi},\mu_t,Y_t^{x,\xi})\\
 &+\widetilde{\dbE}_{\cF_t}\bigl[
 \nabla\widetilde X_t^{\xi,x_i}\pa_{y\mu}H(X_t^{x,\xi},\mu_t,Y_t^{x,\xi},\widetilde X_t^{x_i,\xi})\\
 &\quad +\nabla\widetilde X_t^{\xi,x_i,\star}\pa_{y\mu}H(X_t^{x,\xi},\mu_t,Y_t^{x,\xi},\widetilde X_t^{\xi})\bigr]\bigr]\,\rmd t,
\end{aligned}\\
\begin{aligned}
 \rmd \nabla_\mu Y^{x,\xi,x_i}_t=-{}&\bigl[ \nabla_\mu X^{x,\xi,x_i}_t\pa_{xx}H(X_t^{x,\xi},\mu_t,Y_t^{x,\xi}) +\nabla_\mu Y^{x,\xi,x_i}_t\pa_{xy}H(X_t^{x,\xi},\mu_t,Y_t^{x,\xi})\\
 &-r\nabla_\mu Y^{x,\xi,x_i}_t+\widetilde{\dbE}_{\cF_t}\bigl[
 \nabla\widetilde X_t^{\xi,x_i}\pa_{x\mu}H(X_t^{x,\xi},\mu_t,Y_t^{x,\xi},\widetilde X_t^{x_i,\xi})\\
 &\quad +\nabla\widetilde X_t^{\xi,x_i,\star}\pa_{x\mu}H(X_t^{x,\xi},\mu_t,Y_t^{x,\xi},\widetilde X_t^{\xi})\bigr]\bigr]\,\rmd t+\nabla_\mu Z^{x,\xi,x_i}_t\,\rmd B_t,
\end{aligned}\\
\nabla_\mu X^{x,\xi,x_i}_0=0.
\end{cases}
\end{equation}

These equations admit unique \(L_r^2\) solutions. For
\(\Phi\in\{X,Y,Z\}\), the solutions satisfy
\begin{equation}
\label{eq: main discrete decomposition}
\begin{aligned}
 \delta\Phi^{\xi,I_{\{\xi=x_i\}}}
 &=I_{\{\xi=x_i\}}\nabla\Phi^{\xi,x_i}
   +p_i\nabla\Phi^{\xi,x_i,\star},\\
 \delta\Phi^{x,\xi,I_{\{\xi=x_i\}}}
 &=p_i\nabla_\mu\Phi^{x,\xi,x_i}.
\end{aligned}
\end{equation}
In particular,
\[
 \nabla_\mu Y_0^{x,\xi,x_i}
 =\frac{1}{p_i}\delta Y_0^{x,\xi,I_{\{\xi=x_i\}}}.
\]
This gives the candidate values of
\(\pa_\mu\mcv(x,\cL_\xi,\cdot)\) on the support of \(\cL_\xi\).
For a label outside this support, the same FBSDEs define
the corresponding processes with \(p_i=0\).

Next, for any
\((x,\xi,\tilde x)\in\dbR\times\dbL^2(\cF_0)\times\dbR\),
introduce the following three FBSDEs. The first is
\eqref{eq: dx}, recalled here together with the other two equations:
\begin{equation}
\label{eq: main spatial first variation}
\begin{cases}
\begin{aligned}
 \rmd\nabla X_t^{x,\xi}={}&\bigl[ \nabla X_t^{x,\xi}\pa_{xy}H(X_t^{x,\xi},\mu_t,Y_t^{x,\xi}) +\nabla Y_t^{x,\xi}\pa_{yy}H(X_t^{x,\xi},\mu_t,Y_t^{x,\xi}) \bigr]\,\rmd t,
\end{aligned}\\
\begin{aligned}
 \rmd\nabla Y_t^{x,\xi}=-{}&\bigl[ \nabla X_t^{x,\xi}\pa_{xx}H(X_t^{x,\xi},\mu_t,Y_t^{x,\xi}) +\nabla Y_t^{x,\xi}\pa_{xy}H(X_t^{x,\xi},\mu_t,Y_t^{x,\xi}) -r\nabla Y_t^{x,\xi}\bigr]\,\rmd t\\
 &+\nabla Z_t^{x,\xi}\,\rmd B_t,
\end{aligned}\\
\nabla X_0^{x,\xi}=1.
\end{cases}
\end{equation}

\begin{equation}
\label{eq: main population first variation}
\begin{cases}
\begin{aligned}
 \rmd\nabla X_t^{\xi,\tilde x}={}&\bigl[ \pa_{xy}H(X_t^{\xi},\mu_t,Y_t^{\xi})\nabla X_t^{\xi,\tilde x} +\pa_{yy}H(X_t^{\xi},\mu_t,Y_t^{\xi})\nabla Y_t^{\xi,\tilde x}\\
 &+\widetilde{\dbE}_{\cF_t}\bigl[\nabla\widetilde X_t^{\tilde x,\xi}\pa_{y\mu}H(X_t^{\xi},\mu_t,Y_t^{\xi},\widetilde X_t^{\tilde x,\xi})\\
 &\quad +\nabla\widetilde X_t^{\xi,\tilde x}\pa_{y\mu}H(X_t^{\xi},\mu_t,Y_t^{\xi},\widetilde X_t^{\xi})\bigr]\bigr]\,\rmd t,
\end{aligned}\\
\begin{aligned}
 \rmd\nabla Y_t^{\xi,\tilde x}=-{}&\bigl[ \pa_{xx}H(X_t^{\xi},\mu_t,Y_t^{\xi})\nabla X_t^{\xi,\tilde x} +\pa_{xy}H(X_t^{\xi},\mu_t,Y_t^{\xi})\nabla Y_t^{\xi,\tilde x}-r\nabla Y_t^{\xi,\tilde x}\\
 &+\widetilde{\dbE}_{\cF_t}\bigl[\nabla\widetilde X_t^{\tilde x,\xi}\pa_{x\mu}H(X_t^{\xi},\mu_t,Y_t^{\xi},\widetilde X_t^{\tilde x,\xi})\\
 &\quad +\nabla\widetilde X_t^{\xi,\tilde x}\pa_{x\mu}H(X_t^{\xi},\mu_t,Y_t^{\xi},\widetilde X_t^{\xi})\bigr]\bigr]\,\rmd t+\nabla Z_t^{\xi,\tilde x}\,\rmd B_t,
\end{aligned}\\
\nabla X_0^{\xi,\tilde x}=0.
\end{cases}
\end{equation}

\begin{equation}
\label{eq: main tagged first variation}
\begin{cases}
\begin{aligned}
 \rmd\nabla_\mu X_t^{x,\xi,\tilde x}={}&\bigl[ \pa_{xy}H(X_t^{x,\xi},\mu_t,Y_t^{x,\xi})\nabla_\mu X_t^{x,\xi,\tilde x} +\pa_{yy}H(X_t^{x,\xi},\mu_t,Y_t^{x,\xi})\nabla_\mu Y_t^{x,\xi,\tilde x}\\
 &+\widetilde{\dbE}_{\cF_t}\bigl[\nabla\widetilde X_t^{\tilde x,\xi}\pa_{y\mu}H(X_t^{x,\xi},\mu_t,Y_t^{x,\xi},\widetilde X_t^{\tilde x,\xi})\\
 &\quad +\nabla\widetilde X_t^{\xi,\tilde x}\pa_{y\mu}H(X_t^{x,\xi},\mu_t,Y_t^{x,\xi},\widetilde X_t^{\xi})\bigr]\bigr]\,\rmd t,
\end{aligned}\\
\begin{aligned}
 \rmd\nabla_\mu Y_t^{x,\xi,\tilde x}=-{}&\bigl[ \pa_{xx}H(X_t^{x,\xi},\mu_t,Y_t^{x,\xi})\nabla_\mu X_t^{x,\xi,\tilde x} +\pa_{xy}H(X_t^{x,\xi},\mu_t,Y_t^{x,\xi})\nabla_\mu Y_t^{x,\xi,\tilde x}-r\nabla_\mu Y_t^{x,\xi,\tilde x}\\
 &+\widetilde{\dbE}_{\cF_t}\bigl[\nabla\widetilde X_t^{\tilde x,\xi}\pa_{x\mu}H(X_t^{x,\xi},\mu_t,Y_t^{x,\xi},\widetilde X_t^{\tilde x,\xi})\\
 &\quad +\nabla\widetilde X_t^{\xi,\tilde x}\pa_{x\mu}H(X_t^{x,\xi},\mu_t,Y_t^{x,\xi},\widetilde X_t^{\xi})\bigr]\bigr]\,\rmd t+\nabla_\mu Z_t^{x,\xi,\tilde x}\,\rmd B_t,
\end{aligned}\\
\nabla_\mu X_0^{x,\xi,\tilde x}=0.
\end{cases}
\end{equation}

Equations \eqref{eq: main spatial first variation}--\eqref{eq: main tagged first variation}
have unique solutions in \(L_r^2\), uniformly bounded and jointly
continuous in their parameters. They provide the continuous
representation obtained from the discrete construction above.
In particular, \(\nabla\Phi^{\xi,x_i}\) in
\eqref{eq: main discrete part} is defined by the equation with the
factor \(p_i\), whereas \(\nabla\Phi^{\xi,\tilde x}\) in
\eqref{eq: main population first variation} is defined by the
second equation of this three-equation system.

We recall explicitly the convergence linking these two constructions.
Suppose first that \(\cL_\xi\) is absolutely continuous. Set
\[
 x_i^n=\frac{i}{n},\qquad
 \triangle_i^n=\left[\frac{i}{n},\frac{i+1}{n}\right),
 \qquad i\in\dbZ,
\]
and denote by \(i^n(\tilde x)\) the index for which
\(\tilde x\in\triangle_{i^n(\tilde x)}^n\). Define
\begin{equation}
\label{eq: main discrete approximation}
 \xi_n=\sum_{i=-n^2}^{n^2-1}x_i^n I_{\triangle_i^n}(\xi)
       -nI_{(-\infty,-n)}(\xi)+nI_{[n,\infty)}(\xi).
\end{equation}
Then \(\xi_n\to\xi\) in \(\dbL^2\) and
\(\cW_2(\cL_{\xi_n},\cL_\xi)\to0\). For each fixed
\(\tilde x\), we have
\[
 x_{i^n(\tilde x)}^n\longrightarrow\tilde x,
 \qquad
 p_{i^n(\tilde x)}^n
 =\dbP\bigl(\xi_n=x_{i^n(\tilde x)}^n\bigr)
 =\dbP\bigl(\xi\in\triangle_{i^n(\tilde x)}^n\bigr)
 \longrightarrow0,
\]
where the equality holds for all sufficiently large \(n\).
For \(\Phi=X,Y\), we have
\begin{equation}
\label{eq: main discrete continuous limits}
\begin{aligned}
 \nabla\Phi^{\xi_n,x_{i^n(\tilde x)}^n}
 &\longrightarrow\nabla\Phi^{\tilde x,\xi},\\
 \nabla\Phi^{\xi_n,x_{i^n(\tilde x)}^n,\star}
 &\longrightarrow\nabla\Phi^{\xi,\tilde x},\\
 \nabla_\mu\Phi^{x,\xi_n,x_{i^n(\tilde x)}^n}
 &\longrightarrow\nabla_\mu\Phi^{x,\xi,\tilde x}
\end{aligned}
\end{equation}
in \(L_r^2\), and the corresponding backward initial values
converge in \(\dbL^2\). Here the left-hand sides are the solutions
of \eqref{eq: main discrete part}--\eqref{eq: main discrete total}
with \(\xi_n\) in place of \(\xi\); the right-hand sides are the
solutions of \eqref{eq: main spatial first variation}--\eqref{eq: main tagged first variation}.
In particular,
\[
 \dbE\left|\nabla_\mu Y_0^{x,\xi_n,\xi_n}
               -\nabla_\mu Y_0^{x,\xi,\xi}\right|^2
 \longrightarrow0.
\]

\begin{thm}
\label{thm:main-Lions-field}
\label{prop: main directional construction}
Under Assumption~\ref{assum: value classical}(i),
\eqref{eq: main spatial first variation},
\eqref{eq: main population first variation}, and
\eqref{eq: main tagged first variation} have unique \(L_r^2\)
solutions. These triples and their backward initial values are
uniformly bounded in \(L_r^2\) and \(\dbL^2\), respectively, and
continuous in their parameters, with \(\xi\) varying in \(\dbL^2\).
The function \(\mcv(x,\cdot)\) is Lions differentiable, and
\[
 \pa_\mu\mcv(x,\mu,\tilde x)=
 \nabla_\mu Y_0^{x,\xi,\tilde x}.
\]
For every \(\xi,\eta\in\dbL^2(\cF_0)\),
\[
 \mcv(x,\cL_{\xi+\eta})-\mcv(x,\cL_\xi)
 =\dbE[\pa_\mu\mcv(x,\cL_\xi,\xi)\eta]
  +o(\|\eta\|_2).
\]
This version of the Lions-derivative is bounded and jointly
continuous in \((x,\mu,\tilde x)\).
\end{thm}

\subsubsection{\texorpdfstring{\(L^p_r\) estimates and higher-order derivatives}{Lp estimates and higher-order derivatives}}
\label{subsec: main Lp higher derivatives}

Further differentiation of the variational equations in \(x\) or
\(\tilde x\) produces products of first variations. To control these
terms in \(L_r^2\), we use \(L_r^4\) bounds obtained from the
following \(L^p\) estimates.

For \(p\geq2\), set
\[
 \|U\|_{L_r^p}^p=\dbE\int_0^\infty \rme^{-rt}|U_t|^p\,\rmd t,
 \qquad
 \|Z\|_{\mathbb H_r^p}^p
 =\dbE\left(\int_0^\infty \rme^{-2rt/p}|Z_t|^2\,\rmd t\right)^{p/2}.
\]

Let \(\alpha,\beta,\gamma\) be progressively measurable adapted
real-valued processes, and consider the scalar linear FBSDE
\begin{equation}
\label{eq: main homogeneous linear}
\begin{cases}
 \rmd X_t=(\alpha_tX_t+\beta_tY_t)\,\rmd t,\\
 \rmd Y_t=-(\gamma_tX_t+\alpha_tY_t-rY_t)\,\rmd t+Z_t\,\rmd B_t,\\
 X_0=\eta.
\end{cases}
\end{equation}

\begin{lem}
\label{lem:main-linear-Lp}
Let \(a\geq0\), \(b_0,c_0>r/2\), and suppose that
\begin{equation}
\label{eq: main dissipativity relation}
 b_0c_0>a(r+3a).
\end{equation}
There are upper bounds \(b_1\geq b_0\) and \(c_1\geq c_0\) such that,
for any progressively measurable coefficients satisfying
\begin{equation}
\label{eq: main linear coefficient bounds}
 |\alpha_t|\leq a,\qquad b_0\leq-\beta_t\leq b_1,\qquad
 c_0\leq \gamma_t\leq c_1,
\end{equation}
and any \(\eta\in\dbL^p(\cF_0)\), \(2\leq p<\infty\), the unique
\(L_r^2\) solution of \eqref{eq: main homogeneous linear} belongs to
\(L_r^p\times L_r^p\times\mathbb H_r^p\).
\end{lem}

\begin{proof}
By the Riccati construction of Guatteri and Tessitore
\cite{guatteri2008backward}, set
\[
 F_t(q)=\gamma_t-(r-2\alpha_t)q+\beta_tq^2.
\]
For each \(T>0\), consider the equation
\begin{equation}
\label{eq: main finite Riccati}
 \rmd P_t^T=-F_t(P_t^T)\,\rmd t+Q_t^T\,\rmd B_t,
 \qquad P_T^T=0.
\end{equation}
For \(q\geq0\),
\[
 c_0-(r+2a)q-b_1q^2\leq F_t(q)
 \leq c_1-(r-2a)q-b_0q^2.
\]
The positive roots of the bounding quadratics are
\begin{equation}
\label{eq: main Riccati barriers}
 \underline\pi
 =\frac{\sqrt{(r+2a)^2+4b_1c_0}-(r+2a)}{2b_1},
 \qquad
 \overline\pi
 =\frac{\sqrt{(r-2a)^2+4b_0c_1}-(r-2a)}{2b_0}.
\end{equation}
Comparison gives \(0\leq P_t^T\leq\overline\pi\) and
\(P_t^T\geq v(T-t)\), where
\[
 v'(s)=c_0-(r+2a)v(s)-b_1v(s)^2,\qquad v(0)=0.
\]
Letting \(T\to\infty\), we obtain an infinite-time solution \((P,Q)\)
satisfying
\begin{equation}
\label{eq: main infinite Riccati}
 \rmd P_t=-\bigl[\gamma_t-(r-2\alpha_t)P_t+\beta_tP_t^2\bigr]\,\rmd t
            +Q_t\,\rmd B_t,
 \qquad \underline\pi\leq P_t\leq\overline\pi,
\end{equation}
The solution of \eqref{eq: main homogeneous linear} is given by
\begin{equation}
\label{eq: main Riccati construction}
 X_t=\eta\exp\left(\int_0^t(\alpha_s+\beta_sP_s)\,\rmd s\right),
 \qquad Y_t=P_tX_t,\qquad Z_t=Q_tX_t.
\end{equation}
Set \(\delta_*\triangleq b_0\underline\pi-a\).
For \(a>0\), choosing
\begin{equation}
\label{eq: main admissible b1}
 b_1\in\left[b_0,\,
 \frac{b_0\bigl[b_0c_0-a(r+2a)\bigr]}{a^2}\right)
\end{equation}
gives \(\delta_*>0\); the interval is nonempty
by \eqref{eq: main dissipativity relation}. If \(a=0\), any finite
\(b_1\geq b_0\) gives \(\delta_*>0\). In both cases, \(c_1\)
can be any finite constant with \(c_1\geq c_0\).
Since \(\alpha_t+\beta_tP_t\leq a-b_0\underline\pi=-\delta_*\),
the construction gives, for every \(t\geq0\),
\begin{equation}
\label{eq: main linear path decay}
 |X_t|\leq|\eta|\rme^{-\delta_*t},
 \qquad |Y_t|\leq\overline\pi|\eta|\rme^{-\delta_*t}.
\end{equation}
A simple estimate gives
\begin{equation}
\label{eq: main linear moment bounds}
 \|X\|_{L_r^p}^p+\|Y\|_{L_r^p}^p
 \leq\frac{1+\overline\pi^p}{r+p\delta_*}\dbE|\eta|^p,
 \qquad
 \|Z\|_{\mathbb H_r^p}^p\leq C_p\dbE|\eta|^p,
\end{equation}
where \(C_p\) depends only on \(p,r,a,b_0,b_1,c_0,c_1\).
The detailed proof is given in \cite[Section~3.3]{songyangregularity}.
\end{proof}

The preceding lemma and its proof lead to the following assumption on
the Hamiltonian coefficients.
\begin{assum}
\label{assum: main spatial stability}
There exist constants \(a\geq0\), \(b_1\geq b_0\), and
\(c_1\geq c_0\) such that
\begin{equation}
\label{eq: main Hamiltonian coefficient bounds}
 |\pa_{xy}H|\leq a,\qquad
 b_0\leq-\pa_{yy}H\leq b_1,\qquad
 c_0\leq\pa_{xx}H\leq c_1,
\end{equation}
and
\begin{equation}
\label{eq: main stability margin}
 \delta_*\triangleq b_0\underline\pi-a>0.
\end{equation}
Here \(b_0,c_0\) are as in
Assumption~\ref{assum: value classical}(i), and
\(\underline\pi,\overline\pi\) are given by
\eqref{eq: main Riccati barriers}.
\end{assum}

\begin{prop}
\label{prop: main homogeneous estimates}
Under Assumption~\ref{assum: value classical}(i) and
Assumption~\ref{assum: main spatial stability}, for every
\(2\leq p<\infty\), the spatial variations satisfy
\begin{equation}
\label{eq: main spatial moment bounds}
 \|\nabla X^{x,\xi}\|_{L_r^p}^p+
 \|\nabla Y^{x,\xi}\|_{L_r^p}^p+
 \|\nabla Z^{x,\xi}\|_{\mathbb H_r^p}^p\leq C_p.
\end{equation}
The same bound holds uniformly for the spatial difference quotients
\((\Phi^{x+h,\xi}-\Phi^{x,\xi})/h\), \(0<|h|\leq1\).
For \(\Phi=X,Y\),
\begin{equation}
\label{eq: main first quotient L4 convergence}
 \left\|\frac{\Phi^{x+h,\xi}-\Phi^{x,\xi}}{h}
             -\nabla\Phi^{x,\xi}\right\|_{L_r^4}
 \xrightarrow{h\to0}0.
\end{equation}
The spatial first variations are jointly continuous in \((x,\xi)\)
in \(L_r^4\) for \(X,Y\), with \(\xi\) varying in \(\dbL^2\).
All constants are independent of \(x,\xi\).
\end{prop}

To obtain higher moments for the variations containing measure terms,
we further restrict their common bound \(\lambda_\mu\) from
Assumption~\ref{assum: value classical}(i).

\begin{assum}
\label{assum: main measure stability}
With the constants in \eqref{eq: main Riccati barriers} and
\eqref{eq: main stability margin},
\begin{equation}
\label{eq: main measure contraction}
 \lambda_\mu<
 \frac{\delta_*(r+\delta_*)}
 {r+\delta_*+b_1(1+\overline\pi)}.
\end{equation}
\end{assum}

This condition makes the contraction coefficient
\[
 \vartheta_\mu\triangleq
 \frac{\lambda_\mu}{\delta_*}
 \left(1+\frac{b_1(1+\overline\pi)}{r+\delta_*}\right)
\]
strictly less than one. It ensures that the solutions of
\eqref{eq: main population first variation} and
\eqref{eq: main tagged first variation} belong to
\(L_r^p\times L_r^p\times\mathbb H_r^p\).

\begin{prop}
\label{prop:main-first-moments}
Under Assumption~\ref{assum: value classical}(i) and
Assumptions~\ref{assum: main spatial stability} and
\ref{assum: main measure stability}, for every \(2\leq p<\infty\),
\begin{equation}
\label{eq: main first moment bounds}
\begin{aligned}
 &\sup_{t\geq0}\dbE\bigl[
  |\nabla X_t^{\xi,\tilde x}|^p+
  |\nabla Y_t^{\xi,\tilde x}|^p\bigr]
  +\|\nabla Z^{\xi,\tilde x}\|_{\mathbb H_r^p}^p\leq C_p,\\
 &\|\nabla_\mu X^{x,\xi,\tilde x}\|_{L_r^p}^p+
  \|\nabla_\mu Y^{x,\xi,\tilde x}\|_{L_r^p}^p+
  \|\nabla_\mu Z^{x,\xi,\tilde x}\|_{\mathbb H_r^p}^p\leq C_p.
\end{aligned}
\end{equation}
If \((\xi_n,\tilde x_n)\to(\xi,\tilde x)\) in
\(\dbL^2\times\dbR\), then
\[
 \|\nabla\Phi^{\xi_n,\tilde x_n}
             -\nabla\Phi^{\xi,\tilde x}\|_{L_r^4}
 \longrightarrow0,\qquad \Phi=X,Y.
\]
All constants are independent of \(x,\xi,\tilde x\).
\end{prop}

For \((x,\xi,\tilde x)\in\dbR\times\dbL^2(\cF_0)\times\dbR\),
consider the following three FBSDEs, obtained by differentiating
\eqref{eq: main spatial first variation} in \(x\), and
\eqref{eq: main population first variation} and
\eqref{eq: main tagged first variation} in \(\tilde x\), respectively:
\begin{equation}
\label{eq: main second spatial variation}
\begin{cases}
\begin{aligned}
 \rmd\nabla^2 X_t^{x,\xi}={}&\bigl[ \pa_{xy}H(X_t^{x,\xi},\mu_t,Y_t^{x,\xi})\nabla^2 X_t^{x,\xi} +\pa_{yy}H(X_t^{x,\xi},\mu_t,Y_t^{x,\xi})\nabla^2 Y_t^{x,\xi}\\
 &+\pa_{xxy}H(X_t^{x,\xi},\mu_t,Y_t^{x,\xi})(\nabla X_t^{x,\xi})^2\\
 &+2\pa_{xyy}H(X_t^{x,\xi},\mu_t,Y_t^{x,\xi})\nabla X_t^{x,\xi}\nabla Y_t^{x,\xi}\\
 &+\pa_{yyy}H(X_t^{x,\xi},\mu_t,Y_t^{x,\xi})(\nabla Y_t^{x,\xi})^2\bigr]\,\rmd t,
\end{aligned}\\
\begin{aligned}
 \rmd\nabla^2 Y_t^{x,\xi}=-{}&\bigl[ \pa_{xx}H(X_t^{x,\xi},\mu_t,Y_t^{x,\xi})\nabla^2 X_t^{x,\xi} +\pa_{xy}H(X_t^{x,\xi},\mu_t,Y_t^{x,\xi})\nabla^2 Y_t^{x,\xi}-r\nabla^2 Y_t^{x,\xi}\\
 &+\pa_{xxx}H(X_t^{x,\xi},\mu_t,Y_t^{x,\xi})(\nabla X_t^{x,\xi})^2\\
 &+2\pa_{xxy}H(X_t^{x,\xi},\mu_t,Y_t^{x,\xi})\nabla X_t^{x,\xi}\nabla Y_t^{x,\xi}\\
 &+\pa_{xyy}H(X_t^{x,\xi},\mu_t,Y_t^{x,\xi})(\nabla Y_t^{x,\xi})^2\bigr]\,\rmd t+\nabla^2 Z_t^{x,\xi}\,\rmd B_t,
\end{aligned}\\
\nabla^2 X_0^{x,\xi}=0.
\end{cases}
\end{equation}

\begin{equation}
\label{eq: main population label variation}
\begin{cases}
\begin{aligned}
 \rmd\nabla^2 X_t^{\xi,\tilde x}={}&\bigl[ \pa_{xy}H(X_t^{\xi},\mu_t,Y_t^{\xi})\nabla^2 X_t^{\xi,\tilde x} +\pa_{yy}H(X_t^{\xi},\mu_t,Y_t^{\xi})\nabla^2 Y_t^{\xi,\tilde x}\\
 &+\widetilde{\dbE}_{\cF_t}\bigl[\nabla^2\widetilde X_t^{\tilde x,\xi}\pa_{y\mu}H(X_t^{\xi},\mu_t,Y_t^{\xi},\widetilde X_t^{\tilde x,\xi})\\
 &\quad +(\nabla\widetilde X_t^{\tilde x,\xi})^2\pa_{\tilde x}\pa_{y\mu}H(X_t^{\xi},\mu_t,Y_t^{\xi},\widetilde X_t^{\tilde x,\xi})\\
 &\quad +\nabla^2\widetilde X_t^{\xi,\tilde x}\pa_{y\mu}H(X_t^{\xi},\mu_t,Y_t^{\xi},\widetilde X_t^{\xi})\bigr]\bigr]\,\rmd t,
\end{aligned}\\
\begin{aligned}
 \rmd\nabla^2 Y_t^{\xi,\tilde x}=-{}&\bigl[ \pa_{xx}H(X_t^{\xi},\mu_t,Y_t^{\xi})\nabla^2 X_t^{\xi,\tilde x} +\pa_{xy}H(X_t^{\xi},\mu_t,Y_t^{\xi})\nabla^2 Y_t^{\xi,\tilde x}-r\nabla^2 Y_t^{\xi,\tilde x}\\
 &+\widetilde{\dbE}_{\cF_t}\bigl[\nabla^2\widetilde X_t^{\tilde x,\xi}\pa_{x\mu}H(X_t^{\xi},\mu_t,Y_t^{\xi},\widetilde X_t^{\tilde x,\xi})\\
 &\quad +(\nabla\widetilde X_t^{\tilde x,\xi})^2\pa_{\tilde x}\pa_{x\mu}H(X_t^{\xi},\mu_t,Y_t^{\xi},\widetilde X_t^{\tilde x,\xi})\\
 &\quad +\nabla^2\widetilde X_t^{\xi,\tilde x}\pa_{x\mu}H(X_t^{\xi},\mu_t,Y_t^{\xi},\widetilde X_t^{\xi})\bigr]\bigr]\,\rmd t+\nabla^2 Z_t^{\xi,\tilde x}\,\rmd B_t,
\end{aligned}\\
\nabla^2 X_0^{\xi,\tilde x}=0.
\end{cases}
\end{equation}

\begin{equation}
\label{eq: main tagged label variation}
\begin{cases}
\begin{aligned}
 \rmd\nabla_{\tilde x\mu} X_t^{x,\xi,\tilde x}={}&\bigl[ \pa_{xy}H(X_t^{x,\xi},\mu_t,Y_t^{x,\xi})\nabla_{\tilde x\mu} X_t^{x,\xi,\tilde x} +\pa_{yy}H(X_t^{x,\xi},\mu_t,Y_t^{x,\xi})\nabla_{\tilde x\mu} Y_t^{x,\xi,\tilde x}\\
 &+\widetilde{\dbE}_{\cF_t}\bigl[\nabla^2\widetilde X_t^{\tilde x,\xi}\pa_{y\mu}H(X_t^{x,\xi},\mu_t,Y_t^{x,\xi},\widetilde X_t^{\tilde x,\xi})\\
 &\quad +(\nabla\widetilde X_t^{\tilde x,\xi})^2\pa_{\tilde x}\pa_{y\mu}H(X_t^{x,\xi},\mu_t,Y_t^{x,\xi},\widetilde X_t^{\tilde x,\xi})\\
 &\quad +\nabla^2\widetilde X_t^{\xi,\tilde x}\pa_{y\mu}H(X_t^{x,\xi},\mu_t,Y_t^{x,\xi},\widetilde X_t^{\xi})\bigr]\bigr]\,\rmd t,
\end{aligned}\\
\begin{aligned}
 \rmd\nabla_{\tilde x\mu} Y_t^{x,\xi,\tilde x}=-{}&\bigl[ \pa_{xx}H(X_t^{x,\xi},\mu_t,Y_t^{x,\xi})\nabla_{\tilde x\mu} X_t^{x,\xi,\tilde x} +\pa_{xy}H(X_t^{x,\xi},\mu_t,Y_t^{x,\xi})\nabla_{\tilde x\mu} Y_t^{x,\xi,\tilde x}\\
 &-r\nabla_{\tilde x\mu} Y_t^{x,\xi,\tilde x}+\widetilde{\dbE}_{\cF_t}\bigl[\nabla^2\widetilde X_t^{\tilde x,\xi}\pa_{x\mu}H(X_t^{x,\xi},\mu_t,Y_t^{x,\xi},\widetilde X_t^{\tilde x,\xi})\\
 &\quad +(\nabla\widetilde X_t^{\tilde x,\xi})^2\pa_{\tilde x}\pa_{x\mu}H(X_t^{x,\xi},\mu_t,Y_t^{x,\xi},\widetilde X_t^{\tilde x,\xi})\\
 &\quad +\nabla^2\widetilde X_t^{\xi,\tilde x}\pa_{x\mu}H(X_t^{x,\xi},\mu_t,Y_t^{x,\xi},\widetilde X_t^{\xi})\bigr]\bigr]\,\rmd t+\nabla_{\tilde x\mu} Z_t^{x,\xi,\tilde x}\,\rmd B_t,
\end{aligned}\\
\nabla_{\tilde x\mu} X_0^{x,\xi,\tilde x}=0.
\end{cases}
\end{equation}

\begin{thm}
\label{thm:main-second-spatial}
\label{thm:main-population-label}
\label{thm:main-Lions-label}
Under Assumption~\ref{assum: value classical}(i)--(ii) and
Assumptions~\ref{assum: main spatial stability} and
\ref{assum: main measure stability},
\eqref{eq: main second spatial variation},
\eqref{eq: main population label variation}, and
\eqref{eq: main tagged label variation} have unique \(L_r^2\)
solutions. For \(\Phi=X,Y,Z\),
\begin{equation}
\label{eq: main second variation convergence}
\begin{aligned}
 \frac{\nabla\Phi^{x+h,\xi}-\nabla\Phi^{x,\xi}}h
 &\longrightarrow\nabla^2\Phi^{x,\xi},\\
 \frac{\nabla\Phi^{\xi,\tilde x+h}-\nabla\Phi^{\xi,\tilde x}}h
 &\longrightarrow\nabla^2\Phi^{\xi,\tilde x},\\
 \frac{\nabla_\mu\Phi^{x,\xi,\tilde x+h}
             -\nabla_\mu\Phi^{x,\xi,\tilde x}}h
 &\longrightarrow\nabla_{\tilde x\mu}\Phi^{x,\xi,\tilde x}
\end{aligned}
\end{equation}
in \(L_r^2\) as \(h\to0\). Their backward initial values satisfy
\begin{equation}
\label{eq: main second initial convergence}
\begin{aligned}
 &\left|\frac{\nabla Y_0^{x+h,\xi}-\nabla Y_0^{x,\xi}}{h}
                    -\nabla^2Y_0^{x,\xi}\right|\longrightarrow0,\\
 &\dbE\left|\frac{\nabla Y_0^{\xi,\tilde x+h}
                    -\nabla Y_0^{\xi,\tilde x}}{h}
                    -\nabla^2Y_0^{\xi,\tilde x}\right|^2
                    \longrightarrow0,\\
 &\left|\frac{\nabla_\mu Y_0^{x,\xi,\tilde x+h}
                    -\nabla_\mu Y_0^{x,\xi,\tilde x}}{h}
                    -\nabla_{\tilde x\mu}Y_0^{x,\xi,\tilde x}\right|
                    \longrightarrow0.
\end{aligned}
\end{equation}
All three solution triples and their backward initial values are
uniformly bounded in \(L_r^2\) and \(\dbL^2\), respectively.
If \((x_n,\xi_n,\tilde x_n)\to(x,\xi,\tilde x)\) in
\(\dbR\times\dbL^2\times\dbR\), then, for
\(\Phi\in\{X,Y,Z\}\),
\begin{equation}
\label{eq: main higher derivative continuity}
\begin{aligned}
 \nabla^2\Phi^{x_n,\xi_n}&\longrightarrow\nabla^2\Phi^{x,\xi},\\
 \nabla^2\Phi^{\xi_n,\tilde x_n}&\longrightarrow
                      \nabla^2\Phi^{\xi,\tilde x},\\
 \nabla_{\tilde x\mu}\Phi^{x_n,\xi_n,\tilde x_n}
 &\longrightarrow\nabla_{\tilde x\mu}\Phi^{x,\xi,\tilde x}
\end{aligned}
\end{equation}
in \(L_r^2\). The corresponding backward initial values converge
in \(\dbL^2\); the initial values in the first and third equations
are deterministic.
In particular,
\begin{equation}
\label{eq:main-field-derivative-representations}
\begin{aligned}
 \pa_x\mcv(x,\mu)&=\nabla Y_0^{x,\xi},\\
 \pa_{xx}\mcv(x,\mu)&=\nabla^2Y_0^{x,\xi},\\
 \pa_\mu\mcv(x,\mu,\tilde x)&=\nabla_\mu Y_0^{x,\xi,\tilde x},\\
 \pa_{\tilde x}\pa_\mu\mcv(x,\mu,\tilde x)
   &=\nabla_{\tilde x\mu}Y_0^{x,\xi,\tilde x}.
\end{aligned}
\end{equation}
These four derivatives are bounded and jointly continuous.
In particular, \(\mcv\in C^{2,1}(\dbR\times\cP_2)\).
\end{thm}

\subsection{Regularity of the value function}
\label{subsec: main value measure regularity}

We now apply the preceding variational results to the representation
\eqref{eq: value} of \(V(x,\mu)\). For \(\xi\in\dbL^2(\cF_0)\), set
\(\mu=\cL_\xi\) and \(\mu_t=\cL_{X_t^\xi}\).
At the minimizing feedback, the running cost is
\[
 F(x,\mu,y)=f(x,\mu,\hat{\alpha}(x,\mu,y))
 =H(x,\mu,y)-y\partial_yH(x,\mu,y).
\]
Expectations involving independent copies are taken over all the
displayed probability spaces.

For \(\mu\in\cP_2\), write
\[
 m_2(\mu)=\int_{\dbR}|\tilde x|^2\,\mu(\rmd\tilde x).
\]
The Lipschitz condition on \(\partial_\mu H\) in
Assumption~\ref{assum: value classical}(i) implies
\begin{equation}
\label{eq: value Hamiltonian growth}
 |\partial_\mu H(x,\mu,y,\tilde x)|
 \leq C\bigl(1+|x|+|y|+m_2(\mu)^{1/2}+|\tilde x|\bigr).
\end{equation}
The same linear growth bound holds for \(\partial_{x}F,\partial_{y}F,\partial_{\mu}F\), because
\begin{equation}
\label{eq: value cost derivatives}
 \partial_{x}F=\partial_{x}H-y\partial_{xy}H,\qquad \partial_{y}F=-y\partial_{yy}H,\qquad
 \partial_{\mu}F=\partial_{\mu}H-y\partial_{y\mu}H,
\end{equation}
and the second derivatives appearing here are bounded. The usual
weighted stability estimate for \eqref{eq: mf1}--\eqref{eq: mf2}, applied
also with the zero initial condition, gives
\begin{equation}
\label{eq: value base moment bounds}
\begin{aligned}
 \|X^\xi\|_{L_r^2}+\|Y^\xi\|_{L_r^2}
 &\leq C\bigl(1+\|\xi\|_2\bigr),\\
 \|X^{x,\xi}\|_{L_r^2}+\|Y^{x,\xi}\|_{L_r^2}
 &\leq C\bigl(1+|x|+\|\xi\|_2\bigr).
\end{aligned}
\end{equation}
The first estimate follows from the weighted energy inequality,
since the coefficient values at \((0,\delta_0,0)\) have finite
discounted square integrals. For the second estimate, the population
law is estimated using the first one. The same argument gives
stability with respect to the initial state and its law.

For \(\xi,\eta\in\dbL^2(\cF_0)\), denote the population and representative
player variations by
\(\delta X^{\xi,\eta}\) and \(\delta X^{x,\xi,\eta}\), respectively,
as defined in \eqref{eq: main population directional FBSDE}
and \eqref{eq: main tagged directional FBSDE}.

\begin{prop}
\label{prop: value directional differentiability}
Under Assumption~\ref{assum: value classical}(i), for every
\(x\in\dbR\) and \(\xi,\eta\in\dbL^2(\cF_0)\), the limit
\[
 \delta V^{x,\xi,\eta}
 \triangleq\lim_{\delta\to0}
 \frac{V(x,\cL_{\xi+\delta\eta})-V(x,\cL_\xi)}{\delta}
\]
exists and is given by
\begin{equation}
\label{eq: main value envelope}
 \delta V^{x,\xi,\eta}
 =\dbE\int_0^\infty\rme^{-rt}\widetilde{\dbE}_{\cF_t}
 \left[\partial_\mu H(X_t^{x,\xi},\mu_t,Y_t^{x,\xi},\widetilde X_t^\xi)
                \delta\widetilde X_t^{\xi,\eta}\right]\rmd t.
\end{equation}
For fixed \((x,\xi)\), the map
\(\eta\mapsto\delta V^{x,\xi,\eta}\) is bounded and linear, with
\begin{equation}
\label{eq: value directional bound}
 |\delta V^{x,\xi,\eta}|
 \leq C\bigl(1+|x|+\|\xi\|_2\bigr)\|\eta\|_2.
\end{equation}
Moreover, if \(x_n\to x\), \(\xi_n\to\xi\), and \(\eta_n\to\eta\)
in their respective norms, then
\[
 \delta V^{x_n,\xi_n,\eta_n}
 \longrightarrow\delta V^{x,\xi,\eta}.
\]
\end{prop}

\begin{proof}
Fix \(x,\xi,\eta\), and put \(\xi^\delta=\xi+\delta\eta\). For
\(\theta\in[0,1]\), set
\[
\begin{aligned}
 X_t^{\delta,\theta}&=X_t^\xi+\theta(X_t^{\xi^\delta}-X_t^\xi),
 &\mu_t^{\delta,\theta}&=\cL_{X_t^{\delta,\theta}},\\
 X_t^{x,\delta,\theta}&=X_t^{x,\xi}
                 +\theta(X_t^{x,\xi^\delta}-X_t^{x,\xi}),
 &Y_t^{x,\delta,\theta}&=Y_t^{x,\xi}
                 +\theta(Y_t^{x,\xi^\delta}-Y_t^{x,\xi}).
\end{aligned}
\]
The fundamental theorem of calculus, applied to the spatial
variables and to the measure argument of \(F\), gives
\begin{equation}
\label{eq: value payoff difference quotient}
\begin{aligned}
 &\frac{V(x,\cL_{\xi+\delta\eta})-V(x,\cL_\xi)}{\delta}\\
 ={}&\int_0^1\dbE\int_0^\infty\rme^{-rt}\biggl\{
 \partial_xF(X_t^{x,\delta,\theta},\mu_t^{\delta,\theta},Y_t^{x,\delta,\theta})
           \frac{X_t^{x,\xi^\delta}-X_t^{x,\xi}}{\delta}\\
 &\qquad+\partial_yF(X_t^{x,\delta,\theta},\mu_t^{\delta,\theta},Y_t^{x,\delta,\theta})
           \frac{Y_t^{x,\xi^\delta}-Y_t^{x,\xi}}{\delta}\\
 &\qquad+\widetilde{\dbE}_{\cF_t}\biggl[
 \partial_\mu F(X_t^{x,\delta,\theta},\mu_t^{\delta,\theta},Y_t^{x,\delta,\theta},
                                      \widetilde X_t^{\delta,\theta})
           \frac{\widetilde X_t^{\xi^\delta}-\widetilde X_t^\xi}{\delta}
                        \biggr]\biggr\}\rmd t\,\rmd\theta.
\end{aligned}
\end{equation}
By the linear growth of the derivatives of \(F\),
\eqref{eq: value base moment bounds}, and the convergence in
\eqref{eq: main directional process limits}, H\"{o}lder's inequality
allows us to pass to the limit in
\eqref{eq: value payoff difference quotient}, as in the proof of
Theorem~\ref{thm: deri eq Y}. Hence
\begin{equation}
\label{eq: value uncontracted variation}
\begin{aligned}
 \delta V^{x,\xi,\eta}
 =\dbE\int_0^\infty\rme^{-rt}\bigl\{
 &\partial_xF(X_t^{x,\xi},\mu_t,Y_t^{x,\xi})
                         \delta X_t^{x,\xi,\eta}\\
 &+\partial_yF(X_t^{x,\xi},\mu_t,Y_t^{x,\xi})
                         \delta Y_t^{x,\xi,\eta}\\
 &+\widetilde{\dbE}_{\cF_t}
 [\partial_\mu F(X_t^{x,\xi},\mu_t,Y_t^{x,\xi},\widetilde X_t^\xi)
                       \delta\widetilde X_t^{\xi,\eta}]\bigr\}\rmd t.
\end{aligned}
\end{equation}

Applying It\^{o}'s formula to
\(\rme^{-rt}Y_t^{x,\xi}\delta X_t^{x,\xi,\eta}\), we have
\begin{equation}
\label{eq: value adjoint product}
\begin{aligned}
 \rmd\bigl(\rme^{-rt}Y_t^{x,\xi}\delta X_t^{x,\xi,\eta}\bigr)
 =\rme^{-rt}\bigl\{&
 -\partial_xF(X_t^{x,\xi},\mu_t,Y_t^{x,\xi})
                         \delta X_t^{x,\xi,\eta}\\
 &-\partial_yF(X_t^{x,\xi},\mu_t,Y_t^{x,\xi})
                         \delta Y_t^{x,\xi,\eta}\\
 &+Y_t^{x,\xi}\widetilde{\dbE}_{\cF_t}\bigl[
 \partial_{y\mu}H(X_t^{x,\xi},\mu_t,Y_t^{x,\xi},\widetilde X_t^\xi)
                   \delta\widetilde X_t^{\xi,\eta}\bigr]\bigr\}\rmd t\\
 &+\rme^{-rt}\delta X_t^{x,\xi,\eta}Z_t^{x,\xi}\rmd B_t.
\end{aligned}
\end{equation}
Since \(\delta X_0^{x,\xi,\eta}=0\), integrating and taking
expectations as in the proof of Theorem~\ref{thm: deri eq Y},
and using \(\partial_\mu F=\partial_\mu H-y\partial_{y\mu}H\)
in \eqref{eq: value uncontracted variation}, we obtain
\eqref{eq: main value envelope}.

The variations are bounded linear functions of \(\eta\).
Their \(L_r^2\) estimates and the linear growth of the coefficients
give \eqref{eq: value directional bound}. Continuity follows from
state stability and the continuity of the directional variations
recalled in Section~\ref{subsec: main L2 construction}.
\end{proof}

For fixed \((x,\xi)\), the map
\[
 \eta\longmapsto\delta V^{x,\xi,\eta}
\]
is a bounded linear functional on \(\dbL^2(\cF_0)\). By the Riesz
representation theorem, there is a unique
\(D^{x,\xi}\in\dbL^2(\cF_0)\) such that
\begin{equation}
\label{eq: value Riesz representation}
 \delta V^{x,\xi,\eta}=\dbE[D^{x,\xi}\eta],
 \qquad \eta\in\dbL^2(\cF_0).
\end{equation}

\begin{prop}
\label{prop: value measure directional derivative}
Under Assumption~\ref{assum: value classical}(i), for every
\((x,\mu)\in\dbR\times\cP_2\), there exists a Borel
function \(\psi_0(x,\mu,\cdot)\) such that, for every
\(\xi\in\dbL^2(\cF_0;\mu)\),
\begin{equation}
\label{eq: value abstract kernel}
 D^{x,\xi}=\psi_0(x,\mu,\xi)\quad\dbP\hbox{-a.s.}
\end{equation}
Consequently, for every \(\eta\in\dbL^2(\cF_0)\),
\begin{equation}
\label{eq: value abstract directional representation}
 \delta V^{x,\xi,\eta}
 =\dbE[\psi_0(x,\mu,\xi)\eta].
\end{equation}
Moreover,
\begin{equation}
\label{eq: value abstract kernel bound}
 \|\psi_0(x,\mu,\xi)\|_2
 \leq C\bigl(1+|x|+m_2(\mu)^{1/2}\bigr).
\end{equation}
\end{prop}

\begin{proof}
Fix \(x\) and \(\xi\). We first prove that \(D^{x,\xi}\) is
\(\sigma(\xi)\)-measurable. Suppose otherwise and set
\[
 \eta_1=D^{x,\xi}-\dbE[D^{x,\xi}\mid\xi].
\]
Then
\begin{equation}
\label{eq: value Riesz residual}
 \delta V^{x,\xi,\eta_1}
 =\dbE[D^{x,\xi}\eta_1]=\dbE[|\eta_1|^2]>0.
\end{equation}
By the richness of \(\cF_0\), construct \(\eta_2\) such that,
conditional on \(\xi\), \(\eta_1\) and \(\eta_2\) are independent
and identically distributed. On the one hand,
\eqref{eq: value Riesz residual} holds. On the other hand,
\[
\begin{aligned}
 \delta V^{x,\xi,\eta_2}
 &=\dbE[D^{x,\xi}\eta_2]\\
 &=\dbE[\eta_1\eta_2]
   +\dbE\bigl[\dbE[D^{x,\xi}\mid\xi]\eta_2\bigr]\\
 &=\dbE\bigl[\dbE[\eta_1\mid\xi]\dbE[\eta_2\mid\xi]\bigr]
   +\dbE\bigl[\dbE[D^{x,\xi}\mid\xi]
               \dbE[\eta_2\mid\xi]\bigr]=0.
\end{aligned}
\]
However, \((\xi,\eta_1)\) and \((\xi,\eta_2)\) have the same law.
The definition of the directional derivative therefore gives
\(\delta V^{x,\xi,\eta_1}=\delta V^{x,\xi,\eta_2}\),
which contradicts
\eqref{eq: value Riesz residual}. Hence
\[
 D^{x,\xi}=\dbE[D^{x,\xi}\mid\xi]\quad\dbP\hbox{-a.s.}
\]
The Doob--Dynkin lemma yields a Borel function \(f_{x,\xi}\)
such that \(D^{x,\xi}=f_{x,\xi}(\xi)\).

Let \(\xi_1,\xi_2\) have the same law \(\mu\), and write
\[
 D^{x,\xi_i}=f_i^{x,\mu}(\xi_i),\qquad i=1,2.
\]
For \(g\in C_b(\dbR)\), take \(\eta_i=g(\xi_i)\), \(i=1,2\).
The pairs \((\xi_1,\eta_1)\) and \((\xi_2,\eta_2)\) have the same
law. The definition of the directional derivative gives
\(\delta V^{x,\xi_1,\eta_1}=\delta V^{x,\xi_2,\eta_2}\). Hence
\eqref{eq: value Riesz representation} gives
\[
 \int_{\dbR}f_1^{x,\mu}(\tilde x)g(\tilde x)\mu(\rmd\tilde x)
 =\int_{\dbR}f_2^{x,\mu}(\tilde x)g(\tilde x)\mu(\rmd\tilde x).
\]
Hence \(f_1^{x,\mu}=f_2^{x,\mu}\) \(\mu\)-a.e. We can therefore
choose a Borel function \(\psi_0(x,\mu,\cdot)\) satisfying
\eqref{eq: value abstract kernel} for every
\(\xi\in\dbL^2(\cF_0;\mu)\). Combining this identity with
\eqref{eq: value Riesz representation} gives
\eqref{eq: value abstract directional representation}.
The bound \eqref{eq: value abstract kernel bound} follows from
\eqref{eq: value directional bound} and the Riesz representation.
\end{proof}

Proposition~\ref{prop: value measure directional derivative}
gives a Borel representation \(\psi_0\) of the directional derivative.
To prove Lions differentiability, we construct a jointly continuous
version of \(\psi_0\), following the representation argument in
\cite[Section~4]{songyangregularity}.

The variational bounds in Theorem~\ref{thm:main-Lions-field} give
\begin{equation}
\label{eq: value first variation bounds}
 \sup_{\xi,\tilde x}\bigl(
 \|\nabla X^{\tilde x,\xi}\|_{L_r^2}
 +\|\nabla X^{\xi,\tilde x}\|_{L_r^2}\bigr)\leq C.
\end{equation}
Both variations are jointly continuous in \(L_r^2\).
For \((x,\mu,\tilde x)\in\dbR\times\cP_2\times\dbR\), choose any
\(\xi\in\dbL^2(\cF_0;\mu)\) and define
\begin{equation}
\label{eq: main value kernel}
\begin{aligned}
 \psi(x,\mu,\tilde x)
 =\dbE\int_0^\infty\rme^{-rt}\widetilde{\dbE}_{\cF_t}\bigl[
 &\partial_\mu H(X_t^{x,\xi},\mu_t,Y_t^{x,\xi},\widetilde X_t^{\tilde x,\xi})
                   \nabla\widetilde X_t^{\tilde x,\xi}\\
 &+\partial_\mu H(X_t^{x,\xi},\mu_t,Y_t^{x,\xi},\widetilde X_t^\xi)
                   \nabla\widetilde X_t^{\xi,\tilde x}\bigr]\rmd t.
\end{aligned}
\end{equation}
All copied processes are independent of the representative
player's Brownian motion. Distributional uniqueness of the state
and first-variation equations shows that the right-hand side of
\eqref{eq: main value kernel} does not depend on the particular
choice of \(\xi\) with law \(\mu\). Thus \(\psi\) is well defined.

The bounds \eqref{eq: value Hamiltonian growth},
\eqref{eq: value base moment bounds}, and
\eqref{eq: value first variation bounds} imply that
\eqref{eq: main value kernel} is absolutely integrable and
\begin{equation}
\label{eq: value candidate growth}
 |\psi(x,\mu,\tilde x)|
 \leq C\bigl(1+|x|+m_2(\mu)^{1/2}+|\tilde x|\bigr).
\end{equation}
The map \((x,\mu,\tilde x)\mapsto\psi(x,\mu,\tilde x)\) is
continuous on \(\dbR\times\cP_2\times\dbR\).
Indeed, if \(x_n\to x\), \(\cW_2(\mu_n,\mu)\to0\), and
\(\tilde x_n\to\tilde x\), realize \(\mu_n\) and \(\mu\) by
\(\xi_n\) and \(\xi\) on the rich initial space so that
\(\xi_n\to\xi\) in \(\dbL^2\). State stability and the Lipschitz
condition on \(\partial_\mu H\) give convergence of both
coefficient factors in \eqref{eq: main value kernel} in
\(L_r^2\) on the product space. The corresponding variations
converge in \(L_r^2\) by Theorem~\ref{thm:main-Lions-field}.
Cauchy--Schwarz then gives convergence of the products in
\(L_r^1\), and hence the asserted continuity.
If \(x_n\to x\) and \(\xi_n\to\xi\) in \(\dbL^2\), then
\begin{equation}
\label{eq: value evaluated candidate continuity}
 \psi(x_n,\cL_{\xi_n},\xi_n)
 \longrightarrow\psi(x,\cL_\xi,\xi)
 \quad\hbox{in }\dbL^2.
\end{equation}

\begin{thm}
\label{thm: value Lions regularity}
\label{thm: main value Lions}
Under Assumption~\ref{assum: value classical}(i), for every
\((x,\mu)\in\dbR\times\cP_2\), the function \(\psi(x,\mu,\cdot)\)
is a version of \(\psi_0(x,\mu,\cdot)\) from
Proposition~\ref{prop: value measure directional derivative}.
For each \(x\in\dbR\), \(V(x,\cdot)\) is Lions differentiable, with
\begin{equation}
\label{eq: value Lions kernel identification}
 \partial_\mu V(x,\mu,\tilde x)=\psi(x,\mu,\tilde x),
\end{equation}
that is, for every \(\xi,\eta\in\dbL^2(\cF_0)\),
\begin{equation}
\label{eq: value Lions expansion}
 V(x,\cL_{\xi+\eta})-V(x,\cL_\xi)
 =\dbE[\partial_\mu V(x,\cL_\xi,\xi)\eta]
  +o(\|\eta\|_2).
\end{equation}
The version in \eqref{eq: value Lions kernel identification}
is jointly continuous and satisfies
\begin{equation}
\label{eq: value Lions derivative growth}
 |\partial_\mu V(x,\mu,\tilde x)|
 \leq C\bigl(1+|x|+m_2(\mu)^{1/2}+|\tilde x|\bigr).
\end{equation}
\end{thm}

\begin{proof}
We first consider a discrete random variable \(\xi\) as in
\eqref{eq: main discrete law}. For each atom \(x_i\), the abstract
representation \eqref{eq: value abstract directional representation}
gives
\[
 \delta V^{x,\xi,I_{\{\xi=x_i\}}}
 =p_i\psi_0(x,\cL_\xi,x_i).
\]
By \eqref{eq: main discrete decomposition},
\begin{equation}
\label{eq: value atomic decomposition}
 \delta X_t^{\xi,I_{\{\xi=x_i\}}}
 =I_{\{\xi=x_i\}}\nabla X_t^{\xi,x_i}
       +p_i\nabla X_t^{\xi,x_i,\star},
\end{equation}
where the two variations on the right-hand side are the forward
components of the solution to \eqref{eq: main discrete part}.
On the copied event
\(\{\tilde\xi=x_i\}\), we have
\(\widetilde X^\xi=\widetilde X^{x_i,\xi}\).
By uniqueness, the processes \(X^{x_i,\xi}\) and
\(\nabla X^{\xi,x_i}\) depend on the Brownian motion and the
deterministic law \(\cL_\xi\), and are independent of the initial
variable \(\xi\).
Substituting \eqref{eq: value atomic decomposition} into
\eqref{eq: main value envelope} and taking expectation over
\(\{\tilde\xi=x_i\}\) in the first term, we obtain
\begin{equation}
\label{eq: value discrete kernel identity}
\begin{aligned}
 \psi_0(x,\cL_\xi,x_i)
 ={}&\frac{1}{p_i}\delta V^{x,\xi,I_{\{\xi=x_i\}}}\\
 ={}&\dbE\int_0^\infty\rme^{-rt}
 \widetilde{\dbE}_{\cF_t}\bigl[
 \pa_\mu H(X_t^{x,\xi},\mu_t,Y_t^{x,\xi},
                         \widetilde X_t^{x_i,\xi})
                         \nabla\widetilde X_t^{\xi,x_i}\\
 &\qquad+
 \pa_\mu H(X_t^{x,\xi},\mu_t,Y_t^{x,\xi},
                         \widetilde X_t^\xi)
                         \nabla\widetilde X_t^{\xi,x_i,\star}
 \bigr]\,\rmd t.
\end{aligned}
\end{equation}
Thus the discrete equations give the values of \(\psi_0\) at every
atom of \(\cL_\xi\). Their uniform \(L_r^2\) bounds, together with
\eqref{eq: value Hamiltonian growth} and
\eqref{eq: value base moment bounds}, imply
\begin{equation}
\label{eq: value discrete kernel bound}
 |\psi_0(x,\cL_\xi,x_i)|
 \leq C\bigl(1+|x|+\|\xi\|_2+|x_i|\bigr),
\end{equation}
where \(C\) is independent of the discrete law and of the chosen atom.

We now consider the case where \(\mu=\cL_\xi\) is absolutely
continuous. Let \(\xi_n\) be the discrete approximation in
\eqref{eq: main discrete approximation}, with the notation
\(x_i^n\) and \(i^n(\tilde x)\) introduced there.
For each fixed \(\tilde x\), we have
\(x_{i^n(\tilde x)}^n\to\tilde x\). The convergence
\eqref{eq: main discrete continuous limits} gives
\[
\begin{aligned}
 \nabla X^{\xi_n,x_{i^n(\tilde x)}^n}
 &\longrightarrow\nabla X^{\tilde x,\xi},\\
 \nabla X^{\xi_n,x_{i^n(\tilde x)}^n,\star}
 &\longrightarrow\nabla X^{\xi,\tilde x}
\end{aligned}
\qquad\hbox{in }L_r^2.
\]
State stability and the Lipschitz continuity of \(\pa_\mu H\)
also give \(L_r^2\) convergence of the two coefficient factors
in \eqref{eq: value discrete kernel identity}. Hence each product
converges in \(L_r^1\) by Cauchy--Schwarz.
For \(\mu\)-a.e. \(\tilde x\), the point
\(x_{i^n(\tilde x)}^n\) is an atom of \(\cL_{\xi_n}\) with positive
mass for all sufficiently large \(n\). Passing to the limit in
\eqref{eq: value discrete kernel identity} therefore yields
\begin{equation}
\label{eq: value discrete kernel limit}
 \psi_0(x,\cL_{\xi_n},x_{i^n(\tilde x)}^n)
 \longrightarrow\psi(x,\mu,\tilde x)
 \quad\hbox{for }\mu\hbox{-a.e. }\tilde x.
\end{equation}
For \(\dbP\)-a.e. \(\omega\), we have
\(\xi_n(\omega)=x_{i^n(\xi(\omega))}^n\) for all sufficiently large
\(n\). Moreover, \(|\xi_n|\leq|\xi|+1\). Thus
\eqref{eq: value discrete kernel bound},
\eqref{eq: value candidate growth}, and dominated convergence give
\begin{equation}
\label{eq: value evaluated discrete kernel convergence}
 \psi_0(x,\cL_{\xi_n},\xi_n)
 \longrightarrow\psi(x,\mu,\xi)
 \quad\hbox{in }\dbL^2.
\end{equation}
For any bounded Lipschitz function \(\varphi\),
\(\varphi(\xi_n)\to\varphi(\xi)\) in \(\dbL^2\).
Writing \(p_i^n=\dbP(\xi_n=x_i^n)\), we obtain
\begin{equation}
\label{eq: value kernel test identity}
\begin{aligned}
 \dbE[\psi_0(x,\mu,\xi)\varphi(\xi)]
 &=\delta V^{x,\xi,\varphi(\xi)}\\
 &=\lim_{n\to\infty}\delta V^{x,\xi_n,\varphi(\xi_n)}\\
 &=\lim_{n\to\infty}\sum_{i=-n^2}^{n^2}
   p_i^n\varphi(x_i^n)\psi_0(x,\cL_{\xi_n},x_i^n)\\
 &=\lim_{n\to\infty}
   \dbE[\psi_0(x,\cL_{\xi_n},\xi_n)\varphi(\xi_n)]\\
 &=\dbE[\psi(x,\mu,\xi)\varphi(\xi)].
\end{aligned}
\end{equation}
Here atoms of zero mass are omitted from the sum. The first limit
uses Proposition~\ref{prop: value directional differentiability},
and the last uses
\eqref{eq: value evaluated discrete kernel convergence}.
Since this holds for every bounded Lipschitz function \(\varphi\),
\[
 \psi_0(x,\mu,\cdot)=\psi(x,\mu,\cdot)
 \quad\mu\hbox{-a.e.}
\]
whenever \(\mu\) is absolutely continuous.

For a general \(\xi\in\dbL^2(\cF_0)\), choose \(\xi_n\to\xi\)
in \(\dbL^2\) with absolutely continuous laws. By the case just
proved, the continuity of the directional derivatives, and
\eqref{eq: value evaluated candidate continuity}, we have
\[
\begin{aligned}
 \dbE[\psi_0(x,\cL_\xi,\xi)\varphi(\xi)]
 &=\delta V^{x,\xi,\varphi(\xi)}\\
 &=\lim_{n\to\infty}\delta V^{x,\xi_n,\varphi(\xi_n)}\\
 &=\lim_{n\to\infty}
   \dbE[\psi(x,\cL_{\xi_n},\xi_n)\varphi(\xi_n)]\\
 &=\dbE[\psi(x,\cL_\xi,\xi)\varphi(\xi)].
\end{aligned}
\]
The same argument gives
\begin{equation}
\label{eq: value abstract concrete identification}
 \psi_0(x,\cL_\xi,\cdot)=\psi(x,\cL_\xi,\cdot)
 \quad \cL_\xi\hbox{-a.e.},
 \qquad
 D^{x,\xi}=\psi(x,\cL_\xi,\xi)\quad\dbP\hbox{-a.s.}
\end{equation}
The second identity follows from
\eqref{eq: value abstract kernel}. Thus, for every
\(\eta\in\dbL^2(\cF_0)\),
\begin{equation}
\label{eq: value directional kernel identity}
 \delta V^{x,\xi,\eta}
 =\dbE[\psi(x,\cL_\xi,\xi)\eta].
\end{equation}

By \eqref{eq: value evaluated candidate continuity}, we can
strengthen the G\^{a}teaux derivative in
\eqref{eq: value directional kernel identity} to the Fr\'{e}chet
derivative. Thus \(V(x,\cdot)\) is Lions differentiable, with
\[
 \partial_\mu V(x,\mu,\tilde x)=\psi(x,\mu,\tilde x),
\]
where \(\psi\) is given by \eqref{eq: main value kernel}.
\end{proof}

\begin{prop}
\label{prop: value mixed derivatives}
Under Assumption~\ref{assum: value classical}(i), for every
\((x,\mu,\tilde x)\in\dbR\times\cP_2\times\dbR\),
\begin{equation}
\label{eq: value mixed derivatives}
\begin{aligned}
 \pa_x\pa_\mu V(x,\mu,\tilde x)
 &=\pa_\mu\pa_xV(x,\mu,\tilde x)\\
 &=\pa_\mu\mcv(x,\mu,\tilde x)
 =\nabla_\mu Y_0^{x,\xi,\tilde x},
 \qquad \cL_\xi=\mu.
\end{aligned}
\end{equation}
These derivatives are bounded and jointly continuous.
\end{prop}

\proof
Let \(\xi\in\dbL^2(\cF_0)\), \(\cL_\xi=\mu\), and
\(\mu_t=\cL_{X_t^\xi}\). By Theorem~\ref{thm: deri eq Y},
\(\pa_xV(x,\mu)=Y_0^{x,\xi}\). The measure-variation construction
in Theorem~\ref{thm:main-Lions-field} therefore gives
\[
 \pa_\mu\pa_xV(x,\mu,\tilde x)
 =\nabla_\mu Y_0^{x,\xi,\tilde x}.
\]
We now differentiate the representation
\eqref{eq: main value kernel} of \(\pa_\mu V\) with respect to \(x\).
The population flow and all copied processes in this representation
are independent of \(x\). Moreover,
\(\pa_x(\pa_\mu H)=\pa_{x\mu}H\) and
\(\pa_y(\pa_\mu H)=\pa_{y\mu}H\); these identities follow by
taking the Lions-derivative of the spatial increment formulas for
\(H\), using the boundedness and joint continuity of its mixed
derivatives.

By Theorem~\ref{thm: deri},
\[
 \frac{X^{x+h,\xi}-X^{x,\xi}}h\longrightarrow\nabla X^{x,\xi},
 \qquad
 \frac{Y^{x+h,\xi}-Y^{x,\xi}}h\longrightarrow\nabla Y^{x,\xi}
 \quad\hbox{in }L_r^2.
\]
The mean value formula and the boundedness and continuity of
\(\pa_{x\mu}H,\pa_{y\mu}H\) give \(L_r^2\) convergence of the
corresponding difference quotients of the two \(\pa_\mu H\)
coefficients, on the product space. Indeed, the bounded coefficients
control the errors in the state difference quotients, while dominated
convergence applies to the coefficient differences multiplied by the
fixed limiting variations. Multiplication by the copied variations
then gives convergence in \(L_r^1\) by Cauchy--Schwarz. Thus the
infinite integral can be differentiated, and
\begin{equation}
\label{eq: value mixed derivative construction}
\begin{aligned}
 &\pa_x\pa_\mu V(x,\mu,\tilde x)\\
 ={}&\dbE\int_0^\infty\rme^{-rt}\widetilde{\dbE}_{\cF_t}\Big[
 \Big(
 \pa_{x\mu}H(X_t^{x,\xi},\mu_t,Y_t^{x,\xi},\widetilde X_t^{\tilde x,\xi})
                  \nabla X_t^{x,\xi}\\
 &\qquad+
 \pa_{y\mu}H(X_t^{x,\xi},\mu_t,Y_t^{x,\xi},\widetilde X_t^{\tilde x,\xi})
                  \nabla Y_t^{x,\xi}
 \Big)\nabla\widetilde X_t^{\tilde x,\xi}\\
 &\qquad+\Big(
 \pa_{x\mu}H(X_t^{x,\xi},\mu_t,Y_t^{x,\xi},\widetilde X_t^\xi)
                  \nabla X_t^{x,\xi}\\
 &\qquad+
 \pa_{y\mu}H(X_t^{x,\xi},\mu_t,Y_t^{x,\xi},\widetilde X_t^\xi)
                  \nabla Y_t^{x,\xi}
 \Big)\nabla\widetilde X_t^{\xi,\tilde x}
 \Big]\,\rmd t.
\end{aligned}
\end{equation}
The uniform \(L_r^2\) bounds and the boundedness of the coefficients
show that this derivative is bounded. Joint continuity follows from
the parameter continuity of the variations and the same product
estimates.

By \eqref{eq: main spatial first variation} and
\eqref{eq: main tagged first variation}, It\^{o}'s formula gives
\begin{equation}
\label{eq: value mixed variation duality}
\begin{aligned}
 &\rmd\Bigl[\rme^{-rt}\bigl(
 \nabla X_t^{x,\xi}\nabla_\mu Y_t^{x,\xi,\tilde x}
 -\nabla Y_t^{x,\xi}\nabla_\mu X_t^{x,\xi,\tilde x}
 \bigr)\Bigr]\\
 ={}&-\rme^{-rt}\widetilde{\dbE}_{\cF_t}\Big[
 \Big(
 \pa_{x\mu}H(X_t^{x,\xi},\mu_t,Y_t^{x,\xi},\widetilde X_t^{\tilde x,\xi})
                  \nabla X_t^{x,\xi}\\
 &\qquad+
 \pa_{y\mu}H(X_t^{x,\xi},\mu_t,Y_t^{x,\xi},\widetilde X_t^{\tilde x,\xi})
                  \nabla Y_t^{x,\xi}
 \Big)\nabla\widetilde X_t^{\tilde x,\xi}\\
 &\qquad+\Big(
 \pa_{x\mu}H(X_t^{x,\xi},\mu_t,Y_t^{x,\xi},\widetilde X_t^\xi)
                  \nabla X_t^{x,\xi}\\
 &\qquad+
 \pa_{y\mu}H(X_t^{x,\xi},\mu_t,Y_t^{x,\xi},\widetilde X_t^\xi)
                  \nabla Y_t^{x,\xi}
 \Big)\nabla\widetilde X_t^{\xi,\tilde x}
 \Big]\,\rmd t\\
 &+\rme^{-rt}\bigl(
 \nabla X_t^{x,\xi}\nabla_\mu Z_t^{x,\xi,\tilde x}
 -\nabla_\mu X_t^{x,\xi,\tilde x}\nabla Z_t^{x,\xi}
 \bigr)\,\rmd B_t.
\end{aligned}
\end{equation}
Since \(\nabla X_0^{x,\xi}=1\) and
\(\nabla_\mu X_0^{x,\xi,\tilde x}=0\), integrating and taking
expectations as in the proof of Theorem~\ref{thm: deri eq Y},
and comparing with \eqref{eq: value mixed derivative construction},
we obtain \eqref{eq: value mixed derivatives}. \qed

We next differentiate \(\pa_\mu V\) with respect to \(\tilde x\).
This step uses the higher-moment estimates and second variations
recalled in Section~\ref{subsec: main Lp higher derivatives}, under
Assumption~\ref{assum: value classical}(i)--(ii) and
Assumptions~\ref{assum: main spatial stability} and
\ref{assum: main measure stability}.
Assumption~\ref{assum: value classical}(iii) additionally gives boundedness and continuity of
\(\partial_{\tilde x}\partial_\mu H\).
By Theorem~\ref{thm:main-second-spatial}, the first variations
\(\nabla X^{\tilde x,\xi}\) and \(\nabla X^{\xi,\tilde x}\)
are differentiable with respect to \(\tilde x\) in \(L_r^2\), with
derivatives \(\nabla^2X^{\tilde x,\xi}\) and
\(\nabla^2X^{\xi,\tilde x}\), respectively.
Together with Propositions~\ref{prop: main homogeneous estimates}
and~\ref{prop:main-first-moments}, this gives
\begin{equation}
\label{eq: value higher variation bounds}
 \sup_{\xi,\tilde x}\left(
 \|\nabla X^{\tilde x,\xi}\|_{L_r^4}^2
 +\|\nabla^2X^{\tilde x,\xi}\|_{L_r^2}
 +\|\nabla^2X^{\xi,\tilde x}\|_{L_r^2}\right)\leq C.
\end{equation}
The first variations are continuous in \(L_r^4\), and the second
variations are continuous in \(L_r^2\).

\begin{thm}
\label{thm: value label regularity}
\label{thm: main value label}
Under Assumptions~\ref{assum: value classical},
\ref{assum: main spatial stability}, and
\ref{assum: main measure stability}, the continuous
version \(\partial_\mu V\) is continuously differentiable in \(\tilde x\), with
\begin{equation}
\label{eq: main value label kernel}
\begin{aligned}
 \partial_{\tilde x}\partial_\mu V(x,\mu,\tilde x)
 =\dbE\int_0^\infty\rme^{-rt}\widetilde{\dbE}_{\cF_t}\bigl[
 &\partial_{\tilde x}\partial_\mu H(X_t^{x,\xi},\mu_t,Y_t^{x,\xi},\widetilde X_t^{\tilde x,\xi})
                         (\nabla\widetilde X_t^{\tilde x,\xi})^2\\
 &+\partial_{\mu}H(X_t^{x,\xi},\mu_t,Y_t^{x,\xi},\widetilde X_t^{\tilde x,\xi})
                          \nabla^2\widetilde X_t^{\tilde x,\xi}\\
 &+\partial_{\mu}H(X_t^{x,\xi},\mu_t,Y_t^{x,\xi},\widetilde X_t^\xi)
                          \nabla^2\widetilde X_t^{\xi,\tilde x}\bigr]\rmd t.
\end{aligned}
\end{equation}
Both derivatives are jointly continuous and satisfy
\begin{equation}
\label{eq: main value kernel growth}
 |\partial_\mu V(x,\mu,\tilde x)|+|\partial_{\tilde x}\partial_\mu V(x,\mu,\tilde x)|
 \leq C\bigl(1+|x|+m_2(\mu)^{1/2}+|\tilde x|\bigr).
\end{equation}
If \(x_n\to x\) and \(\xi_n\to\xi\) in \(\dbL^2\), then both
evaluated kernels converge in \(\dbL^2\):
\begin{equation}
\label{eq: value evaluated kernels continuity}
 (\partial_\mu V,\partial_{\tilde x}\partial_\mu V)(x_n,\cL_{\xi_n},\xi_n)
 \longrightarrow(\partial_\mu V,\partial_{\tilde x}\partial_\mu V)(x,\cL_\xi,\xi).
\end{equation}
\end{thm}

\begin{proof}
Only \(X^{\tilde x,\xi}\), \(\nabla X^{\tilde x,\xi}\), and \(\nabla X^{\xi,\tilde x}\) depend on
\(\tilde x\) in \eqref{eq: main value kernel}. Put
\[
 \nabla X^{\tilde x,\xi,h}
 =\frac{X^{\tilde x+h,\xi}-X^{\tilde x,\xi}}h.
\]
For the direct summand its
difference quotient is the sum of
\[
 \partial_{\mu}H(X^{x,\xi},\mu_\cdot,Y^{x,\xi},\widetilde X^{\tilde x+h,\xi})
       \frac{\nabla\widetilde X^{\tilde x+h,\xi}-\nabla\widetilde X^{\tilde x,\xi}}h
\]
and
\[
 \left(\int_0^1 \partial_{\tilde x}\partial_\mu H
 (X^{x,\xi},\mu_\cdot,Y^{x,\xi},\widetilde X^{\tilde x,\xi}
       +s(\widetilde X^{\tilde x+h,\xi}-\widetilde X^{\tilde x,\xi}))\rmd s\right)
                         \nabla\widetilde X^{\tilde x,\xi,h}
                         \nabla\widetilde X^{\tilde x,\xi}.
\]
The first expression converges in \(L_r^1\) to
\(\partial_{\mu}H(X^{x,\xi},\mu_\cdot,Y^{x,\xi},\widetilde X^{\tilde x,\xi})\nabla^2\widetilde X^{\tilde x,\xi}\),
by the \(L_r^2\) convergence of the two factors and
Cauchy--Schwarz. In the second expression,
\(\nabla X^{\tilde x,\xi,h}\to\nabla X^{\tilde x,\xi}\)
in \(L_r^4\), and the coefficient is bounded and converges in
measure. Hence the expression converges in \(L_r^1\) to
\(\partial_{\tilde x}\partial_\mu H(X^{x,\xi},\mu_\cdot,Y^{x,\xi},\widetilde X^{\tilde x,\xi})
 (\nabla\widetilde X^{\tilde x,\xi})^2\). Indeed, the product convergence in \(L_r^2\) follows from
\[
 \|\nabla X^{\tilde x,\xi,h}\nabla X^{\tilde x,\xi}
              -(\nabla X^{\tilde x,\xi})^2\|_{L_r^2}
 \leq\|\nabla X^{\tilde x,\xi,h}-\nabla X^{\tilde x,\xi}\|_{L_r^4}
                               \|\nabla X^{\tilde x,\xi}\|_{L_r^4}.
\]
The remaining coefficient difference is handled by bounded
convergence against the fixed integrable function \(|\nabla X^{\tilde x,\xi}|^2\).

The feedback coefficient
\(\partial_{\mu}H(X^{x,\xi},\mu_\cdot,Y^{x,\xi},\widetilde X^\xi)\) is independent of \(\tilde x\).
Theorem~\ref{thm:main-second-spatial} therefore gives convergence
of its product with the difference quotient of \(\nabla X^{\xi,\tilde x}\) in \(L_r^1\),
again by Cauchy--Schwarz. These are global weighted convergences,
so they justify differentiating the infinite integral and prove
\eqref{eq: main value label kernel}.

For joint continuity, the first variation converges in \(L_r^4\)
under parameter convergence, whereas the two second variations
converge in \(L_r^2\). The \(\partial_{\mu}H\) coefficients converge in
\(L_r^2\), and the \(\partial_{\tilde x}\partial_\mu H\) coefficient is bounded and
converges in measure. The same decompositions prove continuity of
each integral in \eqref{eq: main value label kernel}. Moreover its first term
is bounded by \(\|\partial_{\tilde x}\partial_\mu H\|_\infty\|\nabla X^{\tilde x,\xi}\|_{L_r^2}^2\). The other
two terms are bounded by Cauchy--Schwarz and
\eqref{eq: value higher variation bounds}, giving the linear growth
in \eqref{eq: main value kernel growth}.
The joint continuity and linear growth also give
\eqref{eq: value evaluated kernels continuity}.
\end{proof}

We summarize the regularity of \(V\) needed for the master equation.

\begin{prop}
\label{prop: value regularity input}
Under Assumptions~\ref{assum: value classical},
\ref{assum: main spatial stability}, and
\ref{assum: main measure stability}, the derivatives
\[
 \pa_xV(x,\mu),\quad \pa_{xx}V(x,\mu),\quad
 \pa_\mu V(x,\mu,\tilde x),\quad
 \pa_{\tilde x}\pa_\mu V(x,\mu,\tilde x)
\]
exist and are jointly continuous. For \(\xi\in\dbL^2(\cF_0)\)
with \(\cL_\xi=\mu\),
\begin{equation}
\label{eq: main field derivative summary}
 \pa_xV(x,\mu)=Y_0^{x,\xi},\qquad
 \pa_{xx}V(x,\mu)=\nabla Y_0^{x,\xi},
\end{equation}
and \(\pa_\mu V\), \(\pa_{\tilde x}\pa_\mu V\) are given by
\eqref{eq: main value kernel} and
\eqref{eq: main value label kernel}, respectively.
For some constant \(C>0\),
\begin{equation}
\label{eq: main value growth}
\begin{aligned}
 |V(x,\mu)|&\leq C\bigl(1+|x|^2+m_2(\mu)\bigr),\\
 |\pa_xV(x,\mu)|&\leq C\bigl(1+|x|+m_2(\mu)^{1/2}\bigr),\\
 |\pa_{xx}V(x,\mu)|&\leq C,\\
 |\pa_\mu V(x,\mu,\tilde x)|
 +|\pa_{\tilde x}\pa_\mu V(x,\mu,\tilde x)|
 &\leq C\bigl(1+|x|+m_2(\mu)^{1/2}+|\tilde x|\bigr).
\end{aligned}
\end{equation}
Consequently, \(V\in C^{2,1}(\dbR\times\cP_2)\).
\end{prop}

\section{Classical solvability of the elliptic master equations}
\label{sec: main classical master}
\setcounter{equation}{0}

The regularity established in Section~\ref{sec: main higher regularity}
allows us to use \(V(x,\mu)\) itself as a test function in the viscosity
inequalities and hence verify the master equation pointwise.
We also establish the equation satisfied by \(\pa_xV(x,\mu)\).
We then prove uniqueness in the class of regular solutions and
verify that \(V(x,\mu)\) belongs to this class.

\subsection{Classical solution to the elliptic master equation}

\begin{thm}
\label{thm: main classical existence}
Under Assumptions~\ref{assum: value classical},
\ref{assum: main spatial stability}, and
\ref{assum: main measure stability}, the value function
\(V\) defined in \eqref{eq: value} belongs to
\(C^{2,1}(\dbR\times\cP_2)\) and is a classical solution to the
elliptic master equation \eqref{eq: master}.
\end{thm}

\proof
By Proposition~\ref{prop: value regularity input},
\(V\in C^{2,1}(\dbR\times\cP_2)\). The linear growth of its
measure derivatives in \eqref{eq: main value growth} shows that
\(V\) is a test function.
Theorem~\ref{thm: value viscosity} shows that \(V\) is a viscosity
solution. Taking \(\Psi=V\), every \((x,\mu)\in\dbR\times\cP_2\)
is both a local maximum and a local minimum of \(V-\Psi\).
The two viscosity inequalities therefore give
\begin{equation}
\begin{split}
 rV(x,\mu)={}&H(x,\mu,\pa_xV(x,\mu))
               +\frac12\pa_{xx}V(x,\mu)\\
 &+\tilde{\dbE}\left[
   \frac12\pa_{\tilde x}\pa_\mu V(x,\mu,\tilde\xi)
   +\pa_\mu V(x,\mu,\tilde\xi)
       \pa_yH(\tilde\xi,\mu,\pa_xV(\tilde\xi,\mu))\right],
 \qquad \cL_{\tilde\xi}=\mu.
\end{split}
\end{equation}
This is the desired equation. \qed

To verify that \(V\) also belongs to the class used in the uniqueness
theorem, we next establish the equation satisfied by \(\pa_xV\),
using the backward equation.

\begin{prop}
\label{prop: main derivative equation}
Under Assumption~\ref{assum: value classical}(i)--(ii) and
Assumptions~\ref{assum: main spatial stability} and
\ref{assum: main measure stability},
\(\mcv=\pa_xV\) belongs to \(C^{2,1}(\dbR\times\cP_2)\) and
is a classical solution to \eqref{eq: intro-pa}.
\end{prop}

\proof
Let \(\xi\in\dbL^2(\cF_0)\), \(\cL_\xi=\mu\), and
\(\mu_t=\cL_{X_t^\xi}\).
By Theorem~\ref{thm:main-second-spatial}, the four derivatives in
\eqref{eq:main-field-derivative-representations} are bounded and
jointly continuous, so that \(\mcv\in C^{2,1}(\dbR\times\cP_2)\).
Applying It\^{o}'s formula to
\(\mcv(X_t^{x,\xi},\mu_t)\), comparing with the backward equation
in \eqref{eq: mf2}, and using \eqref{eq: flow}, we obtain
\begin{equation}
\begin{aligned}
 0={}&\dbE\int_0^T\bigg\{
 rY_t^{x,\xi}-\pa_xH(X_t^{x,\xi},\mu_t,Y_t^{x,\xi})\\
 &\quad-\pa_x\mcv(X_t^{x,\xi},\mu_t)
             \pa_yH(X_t^{x,\xi},\mu_t,Y_t^{x,\xi})
       -\frac12\pa_{xx}\mcv(X_t^{x,\xi},\mu_t)\\
 &\quad-\tilde{\dbE}_{\cF_t}\bigg[
 \pa_\mu\mcv(X_t^{x,\xi},\mu_t,\tilde X_t^\xi)
             \pa_yH(\tilde X_t^\xi,\mu_t,\tilde Y_t^\xi)\\
 &\hspace{85pt}+\frac12\pa_{\tilde x}\pa_\mu\mcv
                     (X_t^{x,\xi},\mu_t,\tilde X_t^\xi)
 \bigg]\bigg\}\rmd t.
\end{aligned}
\end{equation}
Dividing by \(T\) and letting \(T\downarrow0\), we have
\begin{equation}
\begin{split}
 r\mcv(x,\mu)={}&\pa_xH(x,\mu,\mcv(x,\mu))
 +\pa_yH(x,\mu,\mcv(x,\mu))\pa_x\mcv(x,\mu)\\
 &+\frac12\pa_{xx}\mcv(x,\mu)
 +\tilde{\dbE}\left[
 \frac12\pa_{\tilde x}\pa_\mu\mcv(x,\mu,\tilde\xi)
 +\pa_\mu\mcv(x,\mu,\tilde\xi)
       \pa_yH(\tilde\xi,\mu,\mcv(\tilde\xi,\mu))\right].
\end{split}
\end{equation}
Here \(\cL_{\tilde\xi}=\mu\). This proves the assertion. \qed

\subsection{Uniqueness of regular solutions}
We now turn to uniqueness of classical solutions to the elliptic master
equation.
In \cite{yang2025discounted}, we provide an example of an elliptic master equation that admits multiple  classical solutions.
Considering a case where $r=2, H(x,\mu,y)=2x^2-y^2$, the PDE (\ref{eq: master}) becomes
\begin{equation}
\begin{aligned}
  2U(x,\mu)=&2x^2-\left(\pa_x U(x,\mu)\right)^2+\frac{1}{2}\pa_{xx}U(x,\mu)\\
  &+\tilde{\mathbb{E} }
  \left[\frac{1}{2}\partial_{\tilde{x}}\partial_\mu U(x,\mu,\tilde{\xi })-2\partial_\mu U(x,\mu,\tilde{\xi })
  \cdot \partial_{x}U(\tilde{\xi },\mu) \right].
\end{aligned}
\end{equation}
By solving for the solutions with  quadratic form, we obtain the following four solutions:
\begin{equation}
  \begin{aligned}
    &U_1(x,\mu)=\frac{1}{2}x^2+\frac{1}{4},\\
    &U_2(x,\mu)=\frac{1}{2}x^2-3x\bar{\mu}+\frac{3}{2}(\bar{\mu})^2+\frac{1}{4},\\
    &U_3(x,\mu)=-x^2-\frac{1}{2},\\
    &U_4(x,\mu)=-x^2+3x\bar{\mu}-\frac{3}{2}(\bar{\mu})^2-\frac{1}{2},
  \end{aligned}
\end{equation}
where we denote by $\bar{\mu}$
the expectation of the measure $\mu$.
This
demonstrates that the classical solutions to the elliptic master equation are not unique in general. 
We therefore introduce the following class of regular solutions.

\begin{defn}
    A function $U: \mathbb{R} \times \mathcal{P}_2 \to \mathbb{R}$ is called a \emph{regular solution} to the elliptic master equation \ref{eq: master} if it satisfies the following conditions:
    \begin{enumerate}
        \item $U \in C^{2,1}$ and $U$ is a classical solution to the elliptic master equation (\ref{eq: master}).

        \item $\partial_x U \in C^{2,1}$ and $\pa_x U$ is a classical solution to PDE (\ref{eq: intro-pa}).

        \item $U$ is \emph{displacement monotone}, i.e., for any $\xi_1, \xi_2 \in \dbL^2(\mathcal{F})$, the following inequality holds:
        \begin{equation}
            \mathbb{E}\left[ 
                (\xi_1 - \xi_2) \left( 
                    \partial_x U(\xi_1, \mathcal{L}_{\xi_1}) - \partial_x U(\xi_2, \mathcal{L}_{\xi_2}) 
                \right) 
            \right] \geq 0.
        \end{equation}

        \item $U$ is of at most quadratic growth. $\pa_xU$ is Lipschitz continuous, and its derivatives,
         $\pa_{xx}U$ and $\pa_{\mu}\pa_xU$  are bounded,  $\pa_{xxx}U$ and $\pa_{\tilde{x}\mu}\pa_xU$
        are of at most linear growth. The derivatives \(\pa_\mu U\) and
        \(\pa_{\tilde{x}\mu}U\) are also of at most linear growth.
    \end{enumerate}
\end{defn}

In addition, the requirement that $U$ be displacement monotone stems from the inherent properties of 
the value function $V(x,\mu)$. From (\ref{eq: 4.2}) and Theorem \ref{thm: deri eq Y}, we can obtain 
        \begin{equation}
            \mathbb{E}\left[ 
                (\xi_1 - \xi_2) \left( 
                    \partial_x V(\xi_1, \mathcal{L}_{\xi_1}) - \partial_x V(\xi_2, \mathcal{L}_{\xi_2}) 
                \right) 
            \right] \geq 0
        \end{equation}
for any $\xi_1, \xi_2 \in \dbL^2(\mathcal{F}_0)$,
which shows that the value function is displacement monotone.
We particularly remark that for the distribution-independent case, as shown in Proposition \ref{prop: classical}, 
the displacement monotonicity condition reduces to the convexity condition.

\begin{thm}
Under Assumption \ref{assum: H},  there exists at most one regular solution to the elliptic master equation (\ref{eq: master}).
\end{thm}
\proof
Suppose that $\bU$ is a regular solution. 
For any $\xi\in \dbL^2(\cF_0)$, consider the following SDE:
\begin{equation}
  \begin{cases}
    \rmd \bX_t=\pa_y H(\bX_t,\cL_{\bX_t},\pa_x\bU( \bX_t,\cL_{\bX_t} ))\rmd t+\rmd B_t,\\
\bX_0=\xi.
  \end{cases}
\end{equation}
Then, we have, for any $T>0$,
\begin{equation}
  \dbE\left[\sup_{t\in[0,T]}|\bX_t|^2\right]<\infty.
\end{equation}

Let $\bY_t\triangleq \pa_x\bU(\bX_t,\cL_{\bX_t})$ and $\bZ_t\triangleq \pa_{xx}\bU(\bX_t,\cL_{\bX_t})$.
Applying It\^{o}'s formula to $\bY_t$ on any finite interval $[0,T]$,
we can derive the following BSDE:
\begin{equation}
   \rmd \bY_{t} = -\left[\partial_x {H} (\bX_{t},\mathcal{L} _{\bX_{t}}, 
   \bY_{t})-r\bY_{t}\right] \rmd t + \bZ_{t}\rmd B_{t}.
\end{equation}
Thus, we deduce that $(\bX,\bY,\bZ)$ satisfies the infinite-time FBSDE:
\begin{equation}
    \begin{cases}
      \rmd \bX_t=\pa_y H(\bX_t,\cL_{\bX_t},\bY_t)\rmd t+\rmd B_t, \\
    \rmd \bY_{t} = -\left[\partial_x {H} (\bX_{t},\mathcal{L} _{\bX_{t}}, \bY_{t})-r\bY_{t}\right] \rmd t + \bZ_{t}\rmd B_{t},\\
    \bX_0=\xi.
    \end{cases}
\end{equation}
We must emphasize that although the triple $(\bX,\bY,\bZ)$ satisfies the FBSDE (\ref{eq: mf1}), 
it is not necessarily an $L_r^2$ solution. Consequently, we cannot utilize the uniqueness of solutions for
 this equation.
Let $(X^{\xi},Y^{\xi},Z^{\xi})$ be the unique $L_r^2$ solution of FBSDE (\ref{eq: mf1}).
We will now utilize the process $(X^{\xi},Y^{\xi},Z^{\xi})$ to prove that the $L_r^2$-norm of 
$(\bX,\bY,\bZ)$ is finite, thereby establishing that $(\bX,\bY,\bZ)$ is, in fact, equal to $(X^{\xi},Y^{\xi},Z^{\xi})$.

Denote
\begin{equation}
  I_t\triangleq \rme^{-rt}(\bX_t-X^{\xi}_t)\left(\pa_x\bU(\bX_t,\cL_{\bX_t})-\pa_x\bU(X^{\xi}_t,\cL_{X^{\xi}_t})   \right).
\end{equation}
Since $\bU$ is displacement monotone, we have $\dbE[I_T]\geq 0$, for any $T>0$.
To apply It\^{o}'s formula to $I_t$, we denote $I_t\triangleq N_t+M_t$, where
\begin{equation}
  \begin{cases}
    N_t\triangleq \rme^{-rt}(\bX_t-X^{\xi}_t)(\bY_t-Y^{\xi}_t);\\
    M_t\triangleq \rme^{-rt}(\bX_t-X^{\xi}_t)(Y^{\xi}_t- \pa_x\bU(X^{\xi}_t,\cL_{X^{\xi}_t}) ).
  \end{cases}
\end{equation}
Applying It\^{o}'s formula to $N_t$ on $[0,T]$ and taking expectation, we derive that
\begin{equation}
  \dbE[N_T]\leq -\frac{r}{2}\dbE\left[
    \int_0^T \rme^{-rt}\left( (\bX_t-X^{\xi}_t)^2+ (\bY_t-Y^{\xi}_t)^2 \right)\rmd t
  \right].
\end{equation}
Next, we apply It\^{o}'s formula to $M_t$ on $[0,T]$ and take expectation. We get that
\begin{equation}
\begin{aligned}
\dbE[M_T] = \dbE\Bigg[ & \int_0^T \rme^{-rt} \Bigg( -r (\bX_t - X^{\xi}_t)(Y^{\xi}_t - \pa_x\bU(X^{\xi}_t, \cL_{X^{\xi}_t}) ) \\
& + \bigl( \pa_y H(\bX_t, \cL_{\bX_t}, \bY_t) - \pa_y H(X_t^{\xi}, \cL_{X_t^{\xi}}, Y_t^{\xi}) \bigr) (Y^{\xi}_t - \pa_x\bU(X^{\xi}_t, \cL_{X^{\xi}_t}) ) \\
& + (\bX_t - X^{\xi}_t) \bigl( -\partial_x H(X^{\xi}_t, \mathcal{L}_{X^{\xi}_t}, Y^{\xi}_t) + rY^{\xi}_t \bigr) \\
& - (\bX_t - X^{\xi}_t) \bigg( \pa_{xx}\bU(X^{\xi}_t, \cL_{X^{\xi}_t}) \pa_y H(X_t^\xi, \cL_{X_t^\xi}, Y_t^\xi) + \frac{1}{2} \pa_{xxx}\bU(X^{\xi}_t, \cL_{X^{\xi}_t}) \\
& \qquad + \tilde{\dbE}_{\cF_t} \Big[ \pa_\mu\pa_x\bU(X^{\xi}_t, \cL_{X^{\xi}_t}, \tilde{X}^{\xi}_t) \pa_y H(\tilde{X}_t^\xi, \cL_{X_t^\xi}, \tilde{Y}_t^\xi) \\
& \qquad\qquad + \frac{1}{2} \pa_{\tilde{x}\mu}\pa_x\bU(X^{\xi}_t, \cL_{X^{\xi}_t}, \tilde{X}^{\xi}_t) \Big] \bigg) \Bigg) \rmd t \Bigg].
\end{aligned}
\end{equation}
Noticing that 
$\pa_yH(x,\mu,y)$ is Lipschitz continuous and considering the growth condition of 
$\pa_x\bU$, we get that
\begin{equation}
\begin{aligned}
   \dbE[M_T]\leq& \epsilon\dbE\left[
    \int_0^T \rme^{-rt}\left( (\bX_t-X^{\xi}_t)^2+ (\bY_t-Y^{\xi}_t)^2 \right)\rmd t
  \right] \\ & +C(\epsilon) \dbE\left[
    \int_0^T \rme^{-rt}\left(1+ (X^{\xi}_t)^2+ (Y^{\xi}_t)^2 \right)\rmd t
  \right].
\end{aligned}
\end{equation}
Here, $\epsilon$ is a positive constant that is less than $r/2$ and $C(\epsilon)$ depends only
 on $\epsilon$.
Then, we get that there exists a constant $C>0$ (independent of $T$), such that
\begin{equation}
\dbE\left[
    \int_0^T \rme^{-rt}\left( (\bX_t-X^{\xi}_t)^2+ (\bY_t-Y^{\xi}_t)^2 \right)\rmd t
  \right]
\leq C   \dbE\left[
    \int_0^T \rme^{-rt}\left(1+ (X^{\xi}_t)^2+ (Y^{\xi}_t)^2 \right)\rmd t
  \right].
\end{equation}
Letting $T\to\infty$, we have that $(\bX,\bY,\bZ)$ belongs to $L_r^2$.
This means $\pa_x \bU(x,\mu)$ is unique.

Now, consider another regular solution $\bU'$, we know that $\pa_x \bU=\pa_x\bU'$.
Denote $F(\mu)\triangleq \bU(x,\mu)-\bU'(x,\mu)$. It satisfies the following PDE:
\begin{equation}
\begin{split}
  r F(\mu)=\tilde{\mathbb{E} }
  \left[\frac{1}{2}\partial_{\tilde{x}}\partial_\mu F(\mu,\tilde{\xi })+\partial_\mu F(\mu,\tilde{\xi })
  \partial_y H(\tilde{\xi },\mu,\partial_x\bU(\tilde{\xi},\mu))\right].
\end{split}
\end{equation}
The linear growth of \(\pa_\mu F\) and \(\pa_{\tilde{x}\mu}F\),
together with the finite-time second-moment estimate for \(\bX\),
allows us to apply It\^{o}'s formula to \(F(\cL_{\bX_t})\). We obtain
\begin{equation}
\begin{aligned}
  \rmd F(\cL_{\bX_t})&={\mathbb{E} }
  \left[\frac{1}{2}\partial_{\tilde{x}}\partial_\mu F(\cL_{\bX_t},\bX_t)+\partial_\mu F(\cL_{\bX_t},\bX_t)
  \partial_y H(\bX_t, \cL_{\bX_t},\partial_x\bU(\bX_t, \cL_{\bX_t}))\right]\rmd t\\
  &=rF(\cL_{\bX_t})\rmd t.
\end{aligned}
\end{equation}
Thus, we know
\begin{equation}
  F(\cL_{\bX_t})=F(\cL_{\xi})\rme^{rt}.
\end{equation}
Recalling that $F$ is of at most quadratic growth and $\bX\in L_r^2$, we can choose a sequence $T_i\to\infty$
such that $\rme^{-rT_i}\dbE[|\bX_{T_i}|^2]\to 0$. Then we know $F(\cL_{\xi}) =0$ for any $\xi\in \dbL^2(\cF_0)$.
Now, we finish the proof of the uniqueness.

\qed

\begin{cor}
\label{cor: main unique regular}
Under Assumptions~\ref{assum: value classical},
\ref{assum: main spatial stability}, and
\ref{assum: main measure stability}, \(V\) is the unique regular
solution to \eqref{eq: master} in the sense of the preceding definition.
\end{cor}

\section*{Acknowledgements}
Song Y. is financially supported by National Key R\&D Program of China (No. 2024YFA1013503 \& No. 2020YFA0712700) and the National Natural Science Foundation of China (No. 12431017).

\end{document}